\documentclass{amsart}
\usepackage[normalem]{ulem}
\usepackage{float}
\usepackage{todonotes}
\usepackage{xcolor}
\usepackage[all]{xy}
\usepackage{extarrows}
\usetikzlibrary{arrows}
\usepackage{microtype, mathtools, mathrsfs, enumerate, dsfont, esvect}
\usepackage{latexsym,amsxtra,amscd,ifthen}
\usepackage{amsfonts}
\usepackage{graphicx}
\usepackage{cleveref}
\usepackage{verbatim}
\usepackage{amsmath}
\usepackage{amsthm}
\usepackage{amssymb}
\usepackage{url}
\usepackage[toc]{appendix}

\usepackage{fullpage}
\usepackage{tikz}
\usetikzlibrary{arrows,shapes,chains}
\usepackage{tikz}
\usetikzlibrary{shapes,arrows,positioning,fit,calc}
\tikzstyle{block1} = [rectangle, draw, thick,
align=center]
\tikzstyle{block2} = [rectangle, draw, thick, 
align=center]
\usepackage{tikz,lipsum,lmodern}
\usepackage[all]{xy}
\usepackage[english]{babel}
\theoremstyle{plain}

\newtheorem{theorem}{Theorem}
\newtheorem{lemma}[theorem]{Lemma}
\newtheorem{proposition}[theorem]{Proposition}
\newtheorem{definition-lemma}[theorem]{Definition-Lemma}
\newtheorem{corollary}[theorem]{Corollary}

\newtheorem{hypothesis}[theorem]{Hypothesis}
\newtheorem{notation}[theorem]{Notation}

\numberwithin{theorem}{section}
\numberwithin{equation}{theorem}

\theoremstyle{definition}
\newtheorem{definition}[theorem]{Definition}
\newtheorem{example}[theorem]{Example}
\newtheorem{remark}[theorem]{Remark}
\newtheorem{question}[theorem]{Question}
\newtheorem*{question*}{Question}

\newcommand{\Z}{\mathbb{Z}}

\newcommand{\im}{\text{\upshape im}}
\DeclareMathOperator{\divv}{div}

\usepackage[utf8]{inputenc}
 \usepackage{multicol}
\usepackage{multicol}

\newcommand*\circled[1]
{\tikz[baseline=(char.base)]{
        \node[shape=circle,draw,inner sep=2pt] (char) {\vphantom{1g}{#1}};}}

\DeclareMathOperator{\Aut}{Aut}
\DeclareMathOperator{\gr}{gr}

\DeclareMathOperator{\GKdim}{GKdim}

\newcommand\kk{{\Bbbk}}

\begin{document}

\title[Weighted Poisson polynomial rings]
{Weighted Poisson polynomial rings in dimension three}

\author{Hongdi Huang, Xin Tang, Xingting Wang, and James J. Zhang}

\address{Huang: Department of Mathematics, Shanghai University,
Shanghai, 200444, China}
\email{hdhuang@shu.edu.cn}

\address{Tang: Department of Mathematics \& Computer Science, 
Fayetteville State University, Fayetteville, NC 28301,
USA}

\email{xtang@uncfsu.edu}

\address{Wang: Department of Mathematics, 
Louisiana State University, Baton Rouge, Louisiana 70803, USA}
\email{xingtingwang@math.lsu.edu}

\address{Zhang: Department of Mathematics, Box 354350,
University of Washington, Seattle, Washington 98195, USA}

\email{zhang@math.washington.edu}

\begin{abstract}
We discuss Poisson structures on a weighted 
polynomial algebra $A:=\Bbbk[x, y, z]$ defined by a homogeneous element $\Omega\in A$, called a {\it potential}. We start with classifying potentials $\Omega$ of degree $\deg(x)+\deg(y)+\deg(z)$ with any positive weight $(\deg(x),\deg(y),\deg(z))$ and list all with isolated singularity. Based on the classification, we study the rigidity of $A$ in terms of graded twistings and classify Poisson fraction fields of $A/(\Omega)$ for irreducible potentials. Using Poisson valuations, we characterize the Poisson automorphism group of $A$ when $\Omega$ has an isolated singularity extending a nice result of Makar-Limanov-Turusbekova-Umirbaev. Finally, Poisson cohomology groups are computed for new classes of Poisson polynomial algebras.
\end{abstract}

\makeatletter
\@namedef{subjclassname@2020}{\textup{2020} Mathematics Subject Classification}
\makeatother




\maketitle

\section*{Introduction}
Poisson algebras are used in classical mechanics to describe observable evolution in Hamiltonian systems. They have been widely studied concerning topics such as  (twisted) Poincar{\'e} duality and modular derivations \cite{LuWW1, LvWZ2, Wa}, Poisson Dixmier-Moeglin equivalence \cite{BLSM, Go1, GLa, LaS, LuWW2}, Poisson enveloping algebras \cite{Ba1, Ba2, Ba3, LvWZ1, LvWZ3}, noncommutative discriminant \cite{BY1, BY2, LY1, NTY} and so on. They have also been utilized to study the representation theory of PI Sklyanin algebras \cite{WWY1, WWY2}. Additionally, Poisson algebras have been investigated in the context of the isomorphism problem, invariant theory, and the cancellation problem \cite{GVW, GW, GWY, Ma}. 

Let $\kk$ be an algebraically closed base field of characteristic zero throughout. Quadratic Poisson structures with $\deg(x_i)=1$ for all $i$ on $\Bbbk[x_1,\ldots, x_n]$ have been applied in various fields, as discussed in papers \cite{Bo, Go2, GLe, LX, Py} and their references. Note that the deformation quantization of such a Poisson structure is the homogeneous coordinate ring of quantum ${\mathbb P}^{n-1}$. It is Calabi-Yau if the corresponding formal Poisson bracket on $\Bbbk[x_1,\ldots,x_n]$ is unimodular \cite{Do}. The modular derivations can be considered a Poisson analog of the Nakayama automorphisms of skew Calabi-Yau algebras. For more information on skew Calabi-Yau algebras, see \cite{RRZ1, RRZ2} and related references. A notable family of quadratic Poisson structures is the elliptic Poisson algebras. These were independently introduced by Feigin and Odesskii \cite{FO98} and Polishchuk \cite{Po}. Elliptic Poisson algebras can be viewed as semi-classical limits of the elliptic Sklyanin algebras studied by Feigin and Odesskii \cite{FO89}. 

A unimodular Poisson structure on $\Bbbk[x, y, z]$ is determined by a potential $\Omega \in \Bbbk[x, y, z]$. Elliptic Poisson algebras in $3$ variables are defined as a particular case by a homogeneous potential $\Omega$ of degree $3$ with an isolated singularity at the origin. Van den Bergh earlier considered these elliptic Poisson algebras in his work on Hochschild homology of 3-dimensional Sklyanin algebras \cite{VdB}. Makar-Limanov-Turusbekova-Umirbaev computed Poisson automorphism groups of these algebras \cite{MTU} when all the generators have degree one. Recently, it has been proved that every connected graded Poisson polynomial algebra is a twist of a unimodular Poisson polynomial algebra \cite{TWZ}.

In associative algebra, Stephenson \cite{St1, St2} has classified and studied the weighted version of connected graded Artin-Schelter regular (or skew Calabi-Yau) algebras of global dimension three. However, there is little knowledge about the Poisson analog of these algebras.

\subsection{General setup}
\label{zzsec0.1}
In this paper, the graded Poisson polynomial algebras $A$ in dimension three are given by a weighted homogeneous potential $\Omega$. We relax two assumptions made in the elliptic ones above: (a) generators being in degree one and (b) isolated singularities of potential $\Omega$. We study those $A$ that exhibit similar Poisson cohomological behaviors to elliptic Poisson algebras. Additionally, we are interested in the Poisson automorphism groups of $A$ and some of its quotients $A/(\Omega-\xi)$ where $\xi\in \Bbbk$. 

Set $\mathbb N=\{0, 1, 2,\cdots\}$ and 
$\mathbb N_+=\{1, 2, 3,\cdots\}$. An algebra $A$ is said to be 
{\it connected graded} if $A=\bigoplus_{i\ge 0} A_i$ is 
$\mathbb N$-graded and $A_0=\kk$. If so, we use $|f|$ to denote 
the degree of a homogeneous element $f\in A$. We 
say $A$ is a {\it connected $w$-graded Poisson algebra} 
(for $w\in \mathbb Z$) if $A=\bigoplus_{i\ge 0} A_i$ is 
a connected graded algebra where the Poisson bracket of 
$A$ satisfies $\{A_i,A_j\}\subseteq A_{i+j-w}$ for all 
$i,j\ge 0$. If $w=0$, we simply say $A$ is a connected graded Poisson algebra. Below is a general setup for some of the main objects in this paper. We shall make it clear in the context when such a setup or part of it is in place.

\begin{hypothesis}
\label{zzhyp0.1}
\begin{enumerate}
\item[(1)]
Let $A\colon= \kk[x, y, z]$ be a weighted polynomial algebra with 
$\deg(x)=a, \deg(y)=b, \deg (z)=c$ for $a,b,c\in 
\mathbb N_+$. 
\item[(2)]
Let $\Omega\in A$ be a nonzero homogeneous element of degree $n>0$. We call $\Omega$ a potential of $A$ and set $w:=n-a-b-c$.
\item[(3)]{\rm{[\Cref{zzdef3.1}]}}
Let $A_\Omega$ denote $\Bbbk[x, y, z]$ with a Poisson structure determined by $\Omega\in \Bbbk[x, y, z]$ as follows
\begin{equation}\notag
\{f,g\}~=~\det \begin{pmatrix} 
f_{x} & f_{y} & f_{z}\\
g_{x} & g_{y} & g_{z}\\
\Omega_{x} & \Omega_{y} & \Omega_{z}
\end{pmatrix}\quad 
\text{for all $f,g\in \Bbbk[x, y, z]$}.
\end{equation}
\item [(4)] 
$w=0$, equivalently $n=a+b+c$.
\end{enumerate}
\end{hypothesis}

\begin{remark}
\label{zzrem0.2}
To save space, the following will be implemented.
\begin{enumerate}
\item[(1)] Concerning the classification of $\Omega$, it is always assumed that $1\leq a \leq b \leq c$. 
\item[(2)]
We will use tables, such as \Cref{t:potential}, at the end of the Introduction and in the \Cref{appendix} to present results concisely.
\item[(3)]
Some computations will not be shown for \Cref{zzthm3.6}, \Cref{zzthm4.2} and \Cref{zzpro6.3} but are available from the authors.
\item[(4)]
When analyzing arguments divided into cases, authors typically provide a detailed analysis for one case and skip details for others if the proofs are similar. All necessary details can be provided upon request.
\end{enumerate}
\end{remark}

\subsection{Classification}
\label{zzsec0.2}
We classify all potentials $\Omega$ in $\kk[x,y,z]$ of degree $a+b+c$ (refer to \Cref{zzthm2.5} for details). The classification for $(a, b, c)=(1, 1, 1)$ is well-known \cite{BM, DH, DML, KM, LX}. We characterize all possible weights $(a, b, c) \in \mathbb{N}_{+}^{3}$ on $\Bbbk[x, y, z]$ that guarantee the existence of potentials $\Omega$ of degree $a+b+c$ with isolated singularities [\Cref{zzlem3.10}]. Together with \Cref{zzthm2.5}, we identify all three possible parametric families of such potentials. 

\begin{theorem}[\Cref{zzlem3.10}]
\label{zzthm0.3}
Assume (1), (2) and (4) of \Cref{zzhyp0.1} with $a\leq b\leq c$ and ${\rm gcd}(a, b, c)=1$. Below is a complete list of potentials with isolated singularities of degree $a+b+c$ (up to graded automorphisms):
\begin{itemize}
\item[(1)] 
$\Omega_1=x^3+y^3+z^3+\lambda xyz$ for $(-\lambda)^3\neq 3\cdot 3\cdot 3$ 
and $(a,b,c)=(1,1,1)$.   
\item[(2)] 
$\Omega_2=x^4+y^4+z^2+\lambda xyz$ for $(-\lambda)^4\neq 4\cdot 4\cdot 2^2$ 
and $(a,b,c)=(1,1,2)$.
\item[(3)] 
$\Omega_3=x^6+y^3+z^2+\lambda xyz$ for $(-\lambda)^6\neq 6\cdot 3^2\cdot 2^3$ 
and $(a,b,c)=(1,2,3)$.  
\end{itemize}
\end{theorem}

Note that these potentials correspond to the homogeneous coordinate rings of Veronese embeddings of elliptic curves in a weighted projective plane. The embeddings are defined by a divisor $D=kP$ with a marked point $P\in E$ and $k=3,2,1$. Our result will hopefully have independent interests in weighted projective spaces.

\subsection{Rigidities}
\label{zzsec0.3}
In \cite{TWZ}, the Poisson version of graded twists of graded associative algebras introduced by \cite{Zh} was used to define a numerical invariant $rgt(A)$ [\Cref{zzdef1.2}] for any $\Z$-graded Poisson algebra $A$. This invariant measures the size of the vector space of graded twists of $A$. If $rgt(A)=0$, then all graded twists of $A$ are isomorphic to $A$, and we call $A$ rigid. 
When $A=\kk[x,y,z]$ is a polynomial Poisson algebra generated 
in degree one, it was shown in \cite[Corollary 6.7]{TWZ} that 
any connected graded unimodular Poisson structure on $A$ is 
rigid if and only if the associated potential $\Omega$ is 
irreducible. We generalize this equivalence to the weighted 
case. 

\begin{theorem}[\Cref{zzthm4.2}]
\label{zzthm0.4}
Assume \Cref{zzhyp0.1}. Then $rgt(A_{\Omega})=0$ if and only if 
$\Omega$ is irreducible.  
\end{theorem}

In Subsections 4.2 and 4.3, we will briefly discuss two other types of rigidities.

\subsection{Automorphism problem}
\label{zzsec0.4}
One of the aims of this paper is to present universal methods for studying Poisson algebra on essential subjects such as Poisson automorphism groups [\Cref{zzthm0.5}] and Poisson cohomologies [\Cref{zzthm0.6}]. 

In \cite{HTWZ1}, Poisson valuations were introduced and used to solve rigidity, automorphism, isomorphism, and embedding problems for various Poisson algebras/fields. We use them to determine Poisson automorphism groups for $A_{\Omega}$ and its quotient $A_{\Omega}/(\Omega-\xi)$ when $\Omega$ is a potential of degree $a+b+c$ with isolated singularity. Our approach provides an alternative method for determining the automorphism groups of three-dimensional elliptic Poisson algebras when the Poisson algebra $A_{\Omega}$ is generated in degree one, differing from \cite{MTU}.  

\begin{theorem}[\Cref{zzthm3.9}] 
\label{zzthm0.5}
Assume \Cref{zzhyp0.1}. Suppose $\Omega$ has an isolated 
singularity. Denote by $P_{\Omega-\xi}=A_{\Omega}/(\Omega-\xi)$ for 
$\xi\in \kk$ and write $P_{\Omega-0}$ simply as $P_{\Omega}$. 
\begin{itemize}
\item[(1)] 
Every Poisson automorphism of $P_{\Omega}$ is graded and 
every Poisson automorphism of $P_{\Omega-\xi}$ is linear when 
$\xi\ne0$. 
\item[(2)] 
Every Poisson automorphism of $A_{\Omega}$ is graded. 
\end{itemize}
Moreover, the explicit automorphism groups of 
$P_{\Omega-\xi}$ and $A_{\Omega}$ are listed in 
Lemmas~\ref{zzlem3.11}--\ref{zzlem3.14}.
\end{theorem}

 If $\Omega$ has no isolated singularity, finding the Poisson automorphism group of $A_\Omega$ becomes challenging. According to \Cref{zzthm3.6}, the Poisson fraction fields $Q=Q(P_\Omega)$ can be divided into three families. For 
convenience, we call $\Omega$ {\it Weyl type} if $Q$ is isomorphic 
to the Poisson Weyl field $K_{Weyl}:=\kk(x,y)$ with $\{x,y\}=1$, 
and we call $\Omega$ {\it quantum type} if $Q$ is isomorphic 
to the Poisson quantum field $K_q:=\kk(x,y)$ with $\{x,y\}=qxy$ 
for some $q\in \kk^\times$. When $\Omega$ is of Weyl type, we can construct many ungraded Poisson automorphisms of $A$ (see \Cref{zzexa3.15}). However, we cannot construct any ungraded Poisson automorphisms of $A$ if $\Omega$ is of quantum type.  We are curious if the Poisson automorphisms of their Poisson algebras are similar to those of elliptic Poisson algebras, which are all graded (see \Cref{zzque3.16}).  

\subsection{Poisson cohomologies}
\label{zzsec0.5}
Poisson cohomologies can be notoriously difficult to calculate. We characterize Poisson cohomological groups for various Poisson algebras in dimension three. Inspired by the $PH^1$-minimality from \cite{TWZ}, we introduce the concept of $uPH^2$-vacancy in \Cref{zzdef5.8} to control the second Poisson cohomology. The property of being $uPH^2$-vacant was implicitly used by Pichereau in \cite[Remark 3.9]{Pi1} for $\Omega$ having an isolated singularity.  
In the theorem below, we generalize \cite[Theorem 0.6]{TWZ} to 
the weighted case. We call an irreducible potential $\Omega$ in 
$\kk[x, y, z]$ {\it balanced} if $\Omega_x \Omega_y\Omega_z\ne 0$ 
for any choice of graded generators $(x, y, z)$; otherwise, we 
call it {\it non-balanced}.

\begin{theorem}
\label{zzthm0.6}
Let $A:=\kk[x, y, z]$ be a connected graded Poisson polynomial 
algebra assuming (1) of \Cref{zzhyp0.1}. Denote by $Z$ the 
Poisson center of $A$. Then, the following statements are 
equivalent.
\begin{enumerate}
\item[(1)]
$rgt(A)=0$ and any homogeneous Poisson derivation of $A$ with 
negative degree is zero. 
\item[(2)]
Any graded twist of $A$ is isomorphic to $A$, and any homogeneous Poisson derivation of $A$ with a negative degree is zero. 
\item[(3)]
The Hilbert series of the graded vector space of Poisson
derivations of $A$ is $\frac{1}{(1-t^a)(1-t^b)(1-t^c)}$.
\item[(4)]
$h_{PH^1(A)}(t)$ is $\frac{1}{1-t^{a+b+c}}$.
\item[(5)] 
$h_{PH^1(A)}(t)$ is equal to $h_Z(t)$.
\item[(6)]

Every Poisson derivation $\phi$ of $A$ has a decomposition 
$\phi=zE+H_{f}$, where $z\in Z$ and $f\in A$. Here, $z$ is 
unique, and $f$ is unique up to a Poisson central element. 
\item[(7)]
Every Poisson derivation of $A$ that vanishes on $Z$ is Hamiltonian. 
\item[(8)]
$A$ is an unimodular Poisson algebra determined by an irreducible potential $\Omega$ that is balanced.
\item[(9)]
$h_{PH^3(A)}(t)-h_{PH^2(A)}(t)=t^{a+b+c}$. 
\item[(10)] $A$ is unimodular and $h_{PH^2(A)}(t)=\frac{1}{t^{a+b+c}}\left(\frac{(1-t^{a+b})(1-t^{a+c})
(1-t^{b+c})}{(1-t^{a+b+c})(1-t^{a})(1-t^{b})(1-t^{c})}-1\right)$.
\item[(11)] $A$ is unimodular and 
$h_{PH^3(A)}(t)=\frac{(1-t^{a+b})(1-t^{a+c})(1-t^{b+c})}
{t^{a+b+c}(1-t^{a+b+c})(1-t^{a})(1-t^{b})(1-t^{c})}$.
\item[(12)] 
$A$ is $uPH^2$-vacant.
\end{enumerate}
\end{theorem}

It is important to note that Van den Bergh \cite{VdB} already computed the Poisson (co)homology of $A$ for the case where $A$ is generated in degree one and $\Omega$ is a cubic polynomial with isolated singularity. Moreover, it was later computed by Pichereau \cite{Pi1} for an arbitrary weighted homogeneous $\Omega$ with isolated singularity. Additional computations can be found in \cite{Pe1, Pe2}. As stated in \Cref{zzthm0.6}, the calculation of Poisson cohomology for $A_{\Omega}$ is possible for both quantum and balanced Weyl types of $\Omega$, regardless of whether or not isolated singularity exists.

In the graph below, all potentials of degree $a+b+c$ listed in \Cref{zzthm2.5} are divided into irreducible potentials in the black circle and reducible ones in the complement, where irreducible potentials are subdivided into three types: isolated singularity, quantum and (balanced and non-balanced) Weyl.  

\vskip 6mm
\parbox{0.5\textwidth}
{\tikzset{
  ephemeral/.style = {opacity=1}
}
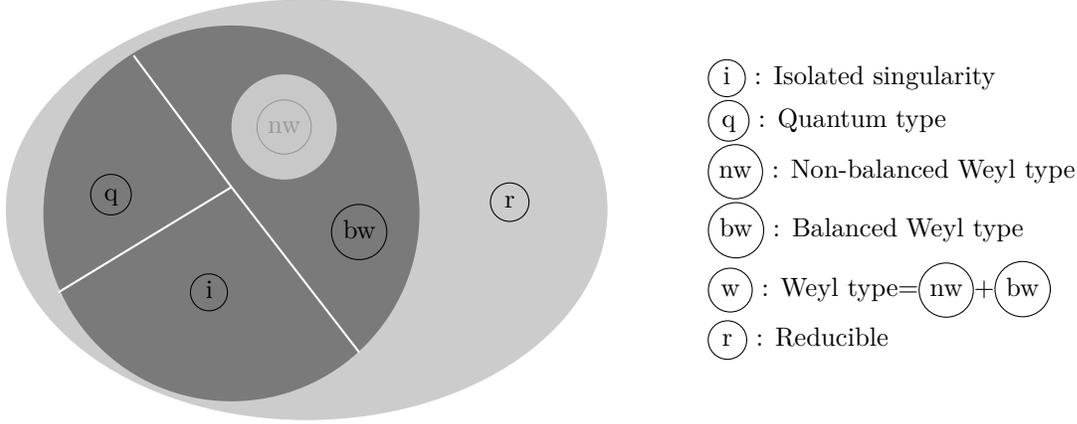
\begin{tikzpicture}
\begin{scope}[transparency group]
\begin{scope}[blend mode=screen]
\fill[black, opacity=.4](2.3,1.85) circle (2.5);
\fill[black, opacity=0.2] 
(3.3, 1.9) ellipse (4 and 2.8);
\fill[white, opacity=0.6](3,3) circle (0.7);
\draw[thick, white] (1,3.95) -- (2.3,2.2);
\draw[thick, white] (0,0.8) -- (2.3,2.2);
\draw[thick, white] (4,0) -- (2.3,2.2);
\node  at (3,3) {\textcolor{black}{\circled{nw}}};
\node  at (2,0.8) {\circled{i}};
\node  at (0.7,2.1) {\circled{q}};
\node  at (4,1.6) {\circled{bw}};
\node  at (6,2) {\circled{r}};
\end{scope}
\end{scope}
\end{tikzpicture}
}
\hskip 10mm
\parbox{0.5\textwidth}{
    \circled{i} : Isolated singularity\\
    \circled{q} : Quantum type\\
    \circled{nw} : Non-balanced Weyl type\\
     \circled{bw} : Balanced Weyl type\\
    \circled{w} : Weyl type=\circled{nw}+\circled{bw}\\
    \circled{r} : Reducible
}
\vskip 6mm

We conclude the introduction with a table summarizing the 
main results for each type of $\Omega$. 
Indeed, the table provides information concerning the smoothness of the projective curve $\Omega=0$ (see 
\Cref{zzrem3.7}), the Gelfand-Kirillov dimension (GKdim) 
of $A_{sing}\colon=A/(\Omega_{x}, \Omega_{y}, \Omega_{z})$, 
the rigidity of $A$ (see \Cref{tab:rgtanddim}), the 
Poisson automorphism group of $A$ (see \Cref{zzthm0.5} and 
\Cref{zzexa3.15}), the Poisson fraction field of $A/\Omega$ 
(see \Cref{zzthm3.6}), the $uPH^2$-vacancy (see \Cref{zzthm0.6}), 
and the $K_1$-sealedness (see \Cref{zzdef5.2}).

\begin{table}[H]
\centering
\caption{\label{t:potential} Potential $\Omega$} 
\vskip 0.6em
\begin{tabular}{|c|cccccccc|}\hline
$\Omega$-type & Proj. Curve  & GKdim  & $rgt(A)$ & $\Aut(A)$
& $Q(P_{\Omega})$ &  $uPH^2$- & $K_1$-   &$h_{PH^i(A)}$\\
& $\Omega=0$ & of $A_{sing}$ & &graded?&&vacant &sealed & computed? \\  [0.6em] 
\hline
\circled{i} &smooth&0 &0&yes& $S_{\zeta, \lambda}$ &yes& yes& yes\\  [1em]
\hline
\circled{q} & nodal singularity &1&0&?&$K_q$& yes & ? &yes\\ [1em]\hline
\circled{bw} &cusp singularity 
    &1&0&no& $K_{Weyl}$ &yes&? &yes\\  [1em]\hline
\circled{nw}  &cusp singularity & 1&0&no& $K_{Weyl}$&no&no &no\\ [1em]\hline
\circled{r} &reducible &\{1,2\}& 
$\leq -1$&some&undefined&no&no&some\\
[1em]\hline
\end{tabular}
\end{table}

This paper is divided into six sections. Section 1 provides basic notations and results for Poisson algebras and briefly describes Poisson valuations. In Section 2, we classify all homogeneous polynomials $\Omega$ in $\kk[x,y,z]$ such that $|\Omega|=|x|+|y|+|z|$ and prove Theorem 0.3. In Section 3, we prove Theorem 0.5; in Section 4, we prove results about several different rigidities, including Theorem 0.4. We study $K_1$-sealedness and $uPH^2$-vacancy in Section 5, which will be useful for the following section. In Section 6, we establish the results on Poisson cohomology for Poisson algebras $\kk[x,y,z]$ with irreducible potentials $\Omega$ of degree $|x|+|y|+|z|$ as summarized in Theorem 0.6.

\section{Preliminaries}
\label{zzsec1}

\subsection{Terminology}
\label{zzsec1.1}
Let $A=\kk[x_1,\ldots,x_n]$ be a polynomial algebra. We denote 
by $\mathfrak X^\bullet(A)=\oplus_{i=0}^{\infty} 
\mathfrak X^{i}(A)$ the set of skew-symmetric multi-derivations 
of $A$.  For $P\in \mathfrak X^p(A)$ and $Q\in \mathfrak X^q(A)$, 
their {\it wedge product} $P\wedge Q\in \mathfrak X^{p+q}(A)$ is 
the skew-symmetric $(p+q)$-derivation of $A$, defined by
\begin{align*}
(P\wedge Q)&(a_1,\ldots,a_{p+q})~:=~\sum_{\sigma\in \mathbb S_{p,q}}
{\rm sgn}(\sigma) P(a_{\sigma(1)},\cdots,a_{\sigma(p)})
\,Q(a_{\sigma(p+1)},\ldots, a_{\sigma(p+q)})
\end{align*}
for all $a_1,\ldots,a_{p+q}\in A$, where $\mathbb S_{p,q}\subset 
\mathbb S_{p+q}$ is the set of all $(p,q)$-shuffles. Note that 
$(\mathfrak X^\bullet(A),\wedge)$ is a graded commutative algebra 
\cite[Proposition 3.1]{LPV}. Recall that the {\it Schouten bracket} 
on $\mathfrak X^\bullet(A)$ is given by 
\[
[\cdot \, , \,\cdot]_S: \mathfrak X^p(A)\times \mathfrak X^q(A)\to 
\mathfrak X^{p+q-1}(A)
\]
such that  
\begin{align*}
[P,\,Q]_S(a_1,&\ldots,a_{p+q-1})
~=~\sum_{\sigma\in \mathbb S_{q,p-1}}{\rm sgn}(\sigma)P
\left(Q(a_{\sigma(1)},\ldots, a_{\sigma(q)}),
a_{\sigma(q+1)},\ldots,a_{\sigma(q+p-1)}\right)\\
&-(-1)^{(p-1)(q-1)}\sum_{\sigma\in \mathbb S_{p,q-1}}{\rm sgn}
(\sigma)Q\left(P(a_{\sigma(1)},\ldots, a_{\sigma(p)}),
a_{\sigma(p+1)},\ldots,a_{\sigma(p+q-1)}\right) 
\end{align*}
for any $P\in \mathfrak X^p(A), Q\in \mathfrak X^q(A) $ and 
$p,q\in \mathbb N$. Note that 
$(\mathfrak X^\bullet(A),\wedge,[-,-]_S)$ is a Gerstenhaber 
algebra \cite[Proposition 3.7]{LPV}. 

Let $\Omega^1(A)$ be the module of K{\" a}hler differentials 
over $A$ and $\Omega^p(A)=\wedge^p_A \Omega^1(A)$ for $p\ge 2$. 
The differential $d: A\to \Omega^1(A)$ extends to a 
well-defined differential of the complex $\Omega^\bullet (A)$ 
and the complex $(\Omega^\bullet ,d)$ is called the algebraic de Rham complex of $A$.

For every $P\in \mathfrak {X}^p(A)$, the {\it internal product} 
with respect to $P$, denoted by $\iota_P$, is an $A$-module map 
$$\iota_P: \Omega^{\bullet}(A)\to \Omega^{\bullet-p}(A)$$
which is determined by
\begin{equation}
\notag
\iota_P(dF_1\wedge dF_2 \wedge \cdots \wedge dF_k)~=~
\begin{cases} 0& k<p,\\
\sum\limits_{\sigma\in {\mathbb S}_{p,k-p}} {\rm sgn}(\sigma)
P(F_{\sigma(1)},\ldots, F_{\sigma(p)}) \\
\qquad\qquad dF_{\sigma(p+1)} \wedge \cdots \wedge dF_{\sigma(k)}
\in \Omega^{k-p}(A) &k\geq p
\end{cases}
\end{equation}
for all $dF_1\wedge dF_2 \wedge \cdots \wedge dF_k \in 
\Omega^k(A)$. Then the {\it Lie derivative} with respect to 
$P$ is defined to be
\begin{equation}
\notag
{\mathcal L}_{P}~=~[\iota_P,d]: \Omega^{\bullet}(A)
\to \Omega^{\bullet-p+1}(A)
\end{equation}
see \cite[(3.49)]{LPV}. Let $\delta\in \mathfrak X^1(A)$ be a 
derivation of $A$. The {\it divergence} of $\delta$, denoted by 
${\rm div}(\delta)$, is an element in $A$ defined by the equation
\begin{equation}
\notag
\mathcal L_\delta (\nu) ~= ~{\rm div}(\delta)\nu,   
\end{equation} 
where $\nu\in \Omega^n(A)$ is a fixed volume form for $A$. In 
particular, if we choose $\nu=dx_1\wedge \cdots \wedge dx_n$, 
from \cite[Lemma 1.2(1)]{TWZ} we get 
\begin{align*}
    {\rm div}(\delta)~=~\sum_{1\leq i\leq n} 
\frac{\partial \delta(x_i)}{\partial x_i}.
\end{align*}

Let $(A,\pi)$ be a Poisson algebra with $\pi\in \mathfrak X^2(A)$ 
satisfying $[\pi,\pi]_S=0$. We usually write the corresponding 
Poisson bracket on $A$ as $\{-,-\}=\pi(-,-)$. A derivation 
$\delta$ on $A$ is called a {\it Poisson derivation} if 
$[\delta,\pi]_S=0$ or $\delta(\{f,g\})=\{\delta(f),g\}
+\{f,\delta(g)\}$ for all $f,g\in A$. There is a special class of 
Poisson derivations on $A$ called {\it Hamiltonian derivations}, 
which is given by $H_f:=\{f,-\}$ for any $f\in A$. The 
{\it modular derivation} of $A$ is defined by
\begin{align*}
\mathfrak m(f)~:=~-\divv(H_f)
\end{align*}
for all $f\in A$. We call $A$ {\it unimodular} if $\mathfrak m=0$.

For each $q\ge 0$, the {\it $q$-th Poisson cohomology} of $A$ is 
defined to be the $q$th-cohomology of the cochain complex 
$(\mathfrak X^\bullet(A),d_\pi^\bullet)$ with differential 
$d_\pi=-[-,\pi]_S$. In particular, for any 
$f\in \mathfrak X^q(A)$, $d_\pi^q(f)\in \mathfrak X^{q+1}(A)$ 
is determined by 
\begin{align}
\label{E1.0.1}\tag{E1.0.1}
d_{\pi}^q(f)(a_0,\ldots,a_q)
~=~&\sum_{i=0}^q (-1)^i \{ a_i, f(a_0,\ldots,\widehat{a_i},\ldots,a_q)\}\\
&+\sum_{0\leq i< j\leq q}
(-1)^{i+j} f(\{a_i,a_j\}, a_0,\ldots,
\widehat{a_i},\ldots,\widehat{a_j},\ldots,a_q)\notag
\end{align}
for any $a_0,a_1,\ldots,a_q\in A$. We denote by 
\begin{align*}
   PH^q(A)~:=~\ker(d_\pi^q)/\im(d_\pi^{q-1}). 
\end{align*}

Let $Pd(A)$ be the Lie algebra of all Poisson derivations of $A$ 
and let $Hd(A)$ be the Lie ideal of $Pd(A)$ consisting of all 
Hamiltonian derivations. We also denote by $Z_P(A)$ the Poisson 
center of $A$. In particular, 
\begin{align*}
PH^0(A)~=~Z_P(A),\quad PH^1(A)~=~Pd(A)/Hd(A).
\end{align*}
For each $q\ge 0$, the {\it $q$-th 
Poisson homology} of $A$ is defined to be the $q$th-homology of 
the chain complex $(\Omega^\bullet(A),\partial^\pi)$, where the 
differentials are given by $\partial^\pi_q=\mathcal L_\pi=[i_{\pi}, d]: 
\Omega^q(A)\to \Omega ^{q-1}(A)$. We denote by 
\begin{align*}
  PH_q(A)~:=~\ker(\partial^\pi_q)/\im(\partial^\pi_{q+1}).  
\end{align*}
When $A$ has a unimodular Poisson structure $\pi$, a duality exists between its Poisson homology and Poisson cohomology \cite{LuWW1}.
\footnote{We may omit these sentences.}

Let's review the concepts of $H$-ozoneness and $PH^1$-minimality about a connected graded Poisson algebra.  

\begin{definition}\cite[Definition 7.1]{TWZ}
\label{zzdef1.1}
Let $A=\kk[x_1,\ldots,x_n]$ be a connected graded Poisson 
algebra with its Poisson center denoted by $Z$.
\begin{itemize}
\item[(1)] 
$\delta\in Pd(A)$ is called {\it ozone} if $\delta(Z)=0$.
\item[(2)] 
Let $Od(A)$ denote the Lie algebra of all ozone Poisson 
derivations of $A$.
\item[(3)] We say $A$ is $H$-ozone if $Od(A)=Hd(A)$, namely, 
any ozone derivation is Hamiltonian. 
\item[(4)] We say $A$ is $PH^1$-minimal if $PH^1(A)\cong ZE$ 
as graded $Z$-modules where $E$ is the Euler derivation 
\eqref{E1.1.1} below.
\end{itemize}
\end{definition}

\subsection{Twists of graded Poisson algebras}
\label{zzsec1.2}
Let $A=\kk[x_1,\ldots,x_n]$ be a graded Poisson polynomial 
algebra with Poisson bracket $\pi=\{-,-\}$ of degree zero. In \cite[\S 2]{TWZ}, 
the notion of graded twists of $A$ was introduced. For any 
homogeneous element $f\in A$, we use $|f|$ to denote its degree 
in $A$. Define the {\it Euler derivation} $E$ of $A$ by
\begin{equation}
\label{E1.1.1}\tag{E1.1.1}
E(f)~:=~|f|\, f
\end{equation}
for all homogeneous elements $f\in A$. We point out that 
$E$ is a Poisson derivation and $\divv(E)=\sum_{i=1}^n \deg(x_i)$. 
Recall that a derivation $\delta$ on $A$ is said to be a 
{\it semi-Poisson derivation} if 
\[[E\wedge \delta,\pi]_S~=~E\wedge [\delta,\pi]_S~=~0.\] 
The set of all graded semi-Poisson derivations (resp. graded 
Poisson derivations) of $A$ is denoted by $Gspd(A)$ 
(resp. $Gpd(A)$). When $A$ is a $\Z$-graded Poisson algebra, 
$Gspd(A)$ is a $\Bbbk$-vector space. For any 
$\delta\in Gspd(A)$, we can define a new Poisson algebra 
$A^\delta:=(A,\pi_{new})$, called {\it a graded twist} of $A$, 
with
\begin{equation}
\label{E1.1.2}\tag{E1.1.2}
\pi_{new}~:=~\pi\,+\, E\wedge \delta
\end{equation}
or namely, $\{f,g\}_{new}=\{f,g\}+E(f)\delta(g)-\delta(f)E(g)$ 
for all homogeneous elements $f,g\in A$.

\begin{definition}\cite[Definition 4.3]{TWZ}
\label{zzdef1.2}
Let $A=\kk[x_1,\ldots,x_n]$ be a $\Z$-graded Poisson algebra. 
The {\it rigidity} of $A$ is defined to be  
\begin{align*}
    rgt(A)~:=~1-\dim_\kk Gspd(A).
\end{align*}

In particular, we say $A$ is {\it rigid} if $rgt(A)=0$.
\end{definition}

Let $(A, \pi)$ be a Poisson algebra with Poisson bracket $\pi$. 
Let $\xi$ be any nonzero scalar. We define a new Poisson bracket 
$\pi_\xi:=\xi\pi$ or $\{-,-\}_\xi := \xi\{-,-\}$ on $A$. Then, it 
is easy to see that $A':=(A, \pi_\xi)$ is a Poisson 
algebra. The following lemma shows how Poisson structures and 
their Poisson cohomologies behave when we replace $A$ (resp. $\pi$) 
by $A'$ (resp. $\pi':=\pi_{\xi}$).

\begin{lemma}\cite[Lemma 1.5]{TWZ}
\label{zzlem1.3}
Retain the notations as above with $\xi\in \kk^\times$.
Let $d_{\pi}^q$ {\rm{(}}resp. $d_{\pi'}^q${\rm{)}} be the 
differential of ${\mathfrak X}^{\bullet}(A)$ {\rm{(}}resp. 
${\mathfrak X}^{\bullet}(A')${\rm{)}} as defined in 
\eqref{E1.0.1}. The following is true. 
\begin{enumerate}
\item[(1)]
$d_{\pi'}^q=\xi d_{\pi}^q$ for all $q$.
\item[(2)]
$\ker(d_{\pi'}^q)=\ker(d_{\pi}^q)$ for all $q$.
\item[(3)]
$\im(d_{\pi'}^q)=\im(d_{\pi}^q)$ for all $q$.
\item[(4)]
$PH^q(A)=PH^q(A')$ for all $q$.
\item[(5)] 
$rgt(A_\Omega)=rgt(A_{\xi\Omega})$.
\item[(6)] 
$A_\Omega$ is $H$-ozone if and only if $A_{\xi\Omega}$ 
is $H$-ozone.
\item[(7)] 
$A_\Omega$ is $PH^1$-minimal if and only if $A_{\xi\Omega}$ 
is $PH^1$-minimal.
\end{enumerate}
\end{lemma}

\subsection{Notations for Poisson (co)homology in dimension three}
\label{zzsec1.3}
We consider the  polynomial algebra $A=\kk[x,y,z]$ with 
grading $(\deg(x),\deg(y),\deg(z))=(a,b,c)\in (\mathbb N_+)^3$. 
Note that a connected $w$-graded unimodular Poisson structure $\pi$ 
on $A$ is determined by a homogeneous polynomial $\Omega$ of 
degree $n$ (not necessarily equal to $a+b+c$). We write 
$(A,\pi_\Omega)$ for the corresponding Poisson algebra where 
the Poisson bracket on $A$ is homogeneous of degree 
$w:=n-a-b-c$. So the cochain complex 
$(\mathfrak X^\bullet(A), d_\pi^\bullet)$ consisting of graded 
vector spaces of skew-symmetric multi-derivations are given by 
\begin{equation}
\label{E1.3.1}\tag{E1.3.1}
0\xlongrightarrow{} \mathfrak X^0(A)
\xlongrightarrow{d_\pi^0}\mathfrak X^1(A)[w]
\xlongrightarrow{d_\pi^1}\mathfrak X^2(A)[2w]
\xlongrightarrow{d_\pi^2}\mathfrak X^3(A)[3w]
\xlongrightarrow{}0.
\end{equation}
Here, we choose the natural isomorphisms of graded vector spaces as follows
\begin{equation}
\label{E1.3.2}\tag{E1.3.2}
\left \{\begin{aligned}
&\mathfrak X^1(A)\xrightarrow{\sim} 
A[a]\oplus A[b]\oplus A[c] \quad && V\mapsto (V(x),V(y), V(z))\\
&\mathfrak X^2(A)\xrightarrow{\sim} A[b+c]
\oplus A[a+c]\oplus A[a+b]\quad
&& V\mapsto (V(y,z),V(z,x),V(x,y))\\
&\mathfrak X^3(A)\xrightarrow{\sim} A[a+b+c] 
\quad && V\mapsto (V(x,y,z)).\\
\end{aligned}\right.
\end{equation}
Using these isomorphisms, it becomes convenient to compute the associated Hilbert series. The elements of $A^{\oplus3}$ are viewed as 
vector-valued functions on $A$, and we denote such an element by $\overrightarrow{f}\in A^{\oplus3}$. Let $\cdot,\times$ 
denote the usual inner and cross products, respectively; 
while $\overrightarrow{\nabla}, \overrightarrow{\nabla}\times$ 
and ${\rm Div}$ denote respectively the gradient, the curl, and 
the divergence operators. Therefore, the cochain complex 
\eqref{E1.3.1} can be identified as 
\begin{equation}
\label{E1.3.3}\tag{E1.3.3}
\centering
\begin{tabular}{ccccc}
&$A[w+a]$ & &$A[2w+b+c]$& \vspace*{-2.1mm}\\
$0\xlongrightarrow{} A \xlongrightarrow{\delta^0_\Omega}$
& \hspace*{-2.5mm}$\oplus A[w+b]$
&\hspace*{-2.5mm}$\xlongrightarrow{\delta^1_\Omega}$
& \hspace*{-2.5mm}$ \oplus A[2w+a+c]$
&\hspace*{-2.5mm}$
\xlongrightarrow{\delta^2_\Omega}A[3w+a+b+c]\xlongrightarrow{}0 $\\
&\hspace*{-2.5mm}$\oplus A[w+c]$& 
&\hspace*{-2.5mm}$\oplus A[2w+a+b]$&
\end{tabular}
\end{equation}
where the differential $\delta_\Omega$ can be written in a 
compact form
\begin{align}
\delta_\Omega^0(f)&~=~ \overrightarrow{\nabla} f\times 
\overrightarrow{\nabla} \Omega,\quad \text{for}\ f\in A
\xrightarrow{\sim}\mathfrak X^0(A),\label{E1.3.4}\tag{E1.3.4}\\
\delta_\Omega^1(\overrightarrow{f})&~=~ -\overrightarrow{\nabla}
(\overrightarrow{f}\cdot  \overrightarrow{\nabla} \Omega)
+{\rm Div}(\overrightarrow{f})\overrightarrow{\nabla} \Omega,
\quad \text{for}\ \overrightarrow{f}\in A^{\oplus3}\xrightarrow{\sim} 
\mathfrak X^1(A),\label{E1.3.5}\tag{E1.3.5}\\
\delta_\Omega^2(\overrightarrow{f})&~=~ -\overrightarrow{\nabla} 
\Omega\cdot( \overrightarrow{\nabla} \times \overrightarrow{f})
=-{\rm Div}(\overrightarrow{f}\times\overrightarrow{\nabla} \Omega),
\quad \text{for}\ \overrightarrow{f}\in A^{\oplus3}\xrightarrow{\sim} 
\mathfrak X^2(A).\label{E1.3.6}\tag{E1.3.6}
\end{align}

For any graded vector space $M=\oplus_{i\in \mathbb Z} M_i$ that 
is locally finite, we use 
\[h_M(t)~=~\sum_{i\in \mathbb Z} \dim_\kk(M_i)\, t^i\]
to denote the {\it Hilbert series} of $M$. Note that the Hilbert series of $A$ is given by
\[
h_A(t)~=~\frac{1}{(1-t^a)(1-t^b)(1-t^c)}.
\]

As a consequence of \eqref{E1.3.3},  the Poisson cohomology 
$PH^{\bullet}(A)$ (resp. Poisson homology $PH_\bullet(A)$) are graded vector spaces, and we denote their Hilbert series as 
$h_{PH^\bullet(A)}(t)$ (resp. $h_{PH_\bullet(A)}(t)$). By additivity of the Hilbert series of \eqref{E1.3.1}-\eqref{E1.3.3}, 
we have
\begin{equation}\label{E1.3.7}\tag{E1.3.7}
\sum_{i=0}^3 (-t^{-w})^{i} h_{PH^i(A)}(t)=
-\frac{1}{t^{3w+a+b+c}}\frac{(1-t^{w+a})(1-t^{w+b})(1-t^{w+c})}
{(1-t^a)(1-t^b)(1-t^c)}.
\end{equation}

\subsection{Poisson valuations and filtrations}
\label{zzsec1.4}
In \cite{HTWZ1}, the notion of Poisson valuations was 
introduced to solve rigidity, automorphism, 
isomorphism, and embedding problems for various classes 
of Poisson algebras/fields. In this subsection, we recall some 
basics of Poisson valuations. Let $w$ be any integer.

\begin{definition}\cite[Definition 1.1]{HTWZ1}
\label{zzdef1.4}
Let $K$ be a Poisson algebra (or a Poisson field) over $\Bbbk$. A {\it $w$-valuation} on $K$ is a map
$$\nu: K \to {\mathbb Z}\cup\{\infty\}$$
which satisfies the following properties: for all $f,g\in K$,
\begin{enumerate}
\item[(a)]
$\nu(f)=\infty$ if and only if $f=0$,
\item[(b)]
$\nu(f)=0$ for all $f\in \Bbbk^{\times}:=\Bbbk\setminus \{0\}$,
\item[(c)]
$\nu(fg)=\nu(f)+\nu(g)$ (assuming $n+\infty=\infty$ when $n\in
{\mathbb Z}\cup\{\infty\}$),
\item[(d)]
$\nu(f+g)\geq \min\{\nu(f),\nu(g)\}$, with equality if $\nu(f)
\neq \nu(g)$.
\item[(e)]
$\nu(\{f,g\})\geq \nu(f)+\nu(g)-w$.
\end{enumerate}
\end{definition}

Note that conditions (a)-(d) mean $\nu$ is an ordinary 
valuation on $K$. Next, we state the definition of $w$-filtration closely related to the Poisson $w$-valuation on $K$.

\begin{definition}\cite[Definition 2.2]{HTWZ1}
\label{zzdef1.5}
Let $A$ be a Poisson algebra. Let 
$\mathbb F=\{F_i\,|\, i\in \mathbb Z\}$ be a 
chain of $\kk$-subspaces of $A$. We say $\mathbb F$ is a {\it $w$-filtration} of $A$ if it satisfies 
\begin{itemize}
\item[(a)] 
$F_i\supseteq F_{i+1}$ 
for all $i\in \mathbb Z$ and $1\in F_0\setminus F_{1}$;
\item[(b)] 
$F_iF_j\subseteq F_{i+j}$ for all $i,j\in \mathbb Z$;
\item[(c)] 
$\cap_{i\in \mathbb Z} F_i=\{0\}$;
\item[(d)] 
$\cup_{i\in \mathbb Z} F_i=A$;
\item[(e)] 
$\{ F_i, F_j\}\subseteq  F_{i+j-w}$ for all $i,j\in \mathbb Z$.
\end{itemize}
\end{definition} 

Let $\mathbb F=\{ F_i\,|\,i\in \mathbb Z\}$ be a $w$-filtration 
of the Poisson algebra $A$. The {\it associated graded algebra} 
of the $w$-filtration $\mathbb F$ of $A$ is defined to be 
$${\rm gr}_\mathbb FA:=\bigoplus_{i\in \mathbb Z}\, F_i/F_{i+1}.$$
For any nonzero element $f\in F_i$, we denote $\overline{f}$ the 
element $f+F_{i+1}$ in the $i$th degree component 
$({\rm gr}_\mathbb FA)_i:=F_i/F_{i+1}$. It is clear that 
${\rm gr}_\mathbb FA$ is a graded algebra. Moreover, by 
\cite[Lemma 2.3]{HTWZ1}, ${\rm gr}_\mathbb FA$ is a $w$-graded 
Poisson algebra with the induced homogeneous Poisson bracket of 
degree $-w$ such that 
\[
\{F_i/ F_{i+1}, F_{j}/ F_{j+1}\}\subseteq  F_{i+j-w}/ F_{i+j+1-w}, 
\] 
namely, $\{({\rm gr}_\mathbb FA)_i, ({\rm gr}_\mathbb FA)_j\}
\subseteq ({\rm gr}_\mathbb FA)_{i+j-w}$ for all $i,j\in  
\mathbb Z$. We call $\mathbb F$  a {\it good} filtration if 
${\rm gr}_\mathbb FA$ is a domain. 

Give a good $w$-filtration $\mathbb F$ on $A$, we define the 
notion of a {\it degree} function, denoted by $\deg: \;
A\to \mathbb Z \cup \{\infty\}$ via 
\begin{align}
\notag
\deg(f):=i\ \text{if $f\in F_i\setminus F_{i+1}$\ and 
\ $\deg(0)=+\infty$.}
\end{align}
One can see that $\deg$ is a valuation on $A$. Conversely, 
given a valuation $\nu$ on $A$, we can define a filtration  
$\mathbb F^\nu_i:=\{F^\nu_i\,|\, i\in \mathbb Z\}$ of $A$ 
(associated to $\nu$) by
\begin{align}
\notag
F_i^\nu:=\{f\in A\,|\, \nu(f)\ge i\}.
\end{align}
The corresponding associated graded algebra of $A$ is 
denoted by ${\rm gr}_{\nu}A.$

\begin{proposition}\cite[Lemma 2.6]{HTWZ1}
\label{zzpro1.6}
Let $A$ be a Poisson algebra. There is a one-to-one 
correspondence between the set of good $w$-filtrations of $A$ 
and the set of $w$-valuations on $A$ via the above constructions. 
\end{proposition}

In this paper, we will mainly focus on the following special 
class of Poisson valuations. 

\begin{definition}\cite[Definition 3.1]{HTWZ1}
\label{zzdef1.7}
Let $K$ be a Poisson field over $\kk$. 
A $w$-valuation $\nu$ on $K$ is called a 
{\it faithful $w$-valuation}  if the following hold.
\begin{enumerate}
\item[(a)] 
The image of $\nu$ is $\mathbb Z\cup \{\infty\}$.
\item[(b)] 
$\GKdim(\gr_\nu K)=\GKdim(K)$.
\item[(c)] 
The $w$-graded Poisson bracket on $\gr_\nu K$ is nonzero. 
\end{enumerate}
A Poisson $w$-valuation on a Poisson domain $A$ is faithful if its natural extension to the Poisson fractional field $Q(A)$ is faithful. 
\end{definition}

Note that the above conditions (b), and (c) are different from 
the original ones \cite[Definition 3.1(1)]{HTWZ1}. 
But it is easy to see they are equivalent 
(see \cite[Definitions 0.2 and 2.18]{HTWZ2}).

\section{Classification of potentials $\Omega$ 
in $\Bbbk[x,y,z]$}
\label{zzsec2}

In this section, we classify all possible homogeneous polynomials $\Omega$ satisfying (1), (2) and (4) of \Cref{zzhyp0.1}, up to some graded automorphism of $A$. We will use the following definition. 

\begin{definition}
\label{zzdef2.1}
We define the {\it Jacobian quotient algebra} of $A$ with respect 
to $\Omega$ to be 
\[A_{sing}:=\kk[x,y,z]/(\Omega_x,\Omega_y,\Omega_z).\]
It is clear that $A_{sing}$ is independent of the choices of 
graded generators $(x,y,z)$.
\end{definition}

\begin{lemma}
\label{zzlem2.2}
Assume (1), (2) and (4) of \Cref{zzhyp0.1}. 
\begin{enumerate}
\item[(1)] 
The following are equivalent.
\begin{enumerate}
\item[(1a)]
There is an irreducible homogeneous potential $\Omega$ in $\kk[x,y]$.
\item[(1b)]
$2<a<b<c=\frac{ab}{\gcd(a,b)}-a-b$. 
\end{enumerate}
In this case, $\Bbbk[x,y]_{c}=0$, and, up to a graded automorphism of $A$, 
$\Omega=x^{\frac{b}{\gcd(a,b)}}+y^{\frac{a}{\gcd(a,b)}}$. 
\item[(2)]
Let $\Omega$ denote a nonzero homogeneous polynomial in $\kk[x,y]$. 
The following are equivalent.
\begin{enumerate}
\item[(2a)]
$A$ has a nonzero graded derivation $\delta$ of degree $0$ 
satisfying $\divv(\delta)=\delta(\Omega)=0$.
\item[(2b)]
$\Omega$ is reducible.
\end{enumerate}
\end{enumerate}
\end{lemma}

\begin{proof}
(1) $(1b)\Rightarrow (1a)$: This is clear by taking 
$\Omega=x^{\frac{b}{\gcd(a,b)}}+y^{\frac{a}{\gcd(a,b)}}$.

$(1a)\Rightarrow (1b)$:
Let $\Omega=h(x,y)\in \kk[x, y]$ be an irreducible polynomial 
of degree $a+b+c$. We have $\deg(h(x,y))=a+b+c=ka+lb$ for 
some $k, l\in \mathbb N$. As a result, $c=(k-1)a+(l-1)b$. If 
needed, we can divide the degrees of $x$ and $y$ by 
${\rm gcd}(a,b)$ and thus assume that ${\rm gcd}(a, b)=1$. If 
no more than one of the $x^{\frac{a+b+c}{a}}$ and 
$y^{\frac{a+b+c}{b}}$ terms appears in $h(x,y)$, then 
$h(x,y)=xf(x,y)$ or $h(x,y)=yf(x,y)$ for some non-constant 
polynomial $f(x,y)\in \Bbbk[x,y]$. In this case, 
$\Omega=h(x,y)$ is reducible. Next we consider the case where 
$h(x,y)$ contains both $x^{\frac{a+b+c}{a}}$ and 
$y^{\frac{a+b+c}{b}}$ terms. We have $a\mid b+c$ and $b\mid a+c$, 
which implies that $a\mid l$ and $b\mid k$. Say $l=m_1a$ and 
$k=m_2b$ for some $m_1,m_2 \in \mathbb N_+$. Thus $a+b+c=(m_1+m_2)ab$. So 
we have 
\[
h(x,y)=\lambda_{m_1+m_2}\,x^{(m_1+m_2)b}+\lambda_{m_1+m_2-1}\,x^{(m_1+m_2-1)b}y^a
+\cdots+\lambda_1\,x^by^{(m_1+m_2-1)a}+\lambda_0\,y^{(m_1+m_2)a}
\] 
for some $\lambda_0,\ldots,\lambda_{m_1+m_2}\in \kk$. We rewrite 
$p=x^b$ and $q=y^a$. Then, we can get 
\[
h(x,y)=h(p,q)=\lambda_{m_1+m_2}\,p^{m_1+m_2}+\lambda_{m_1+m_2-1}\,
p^{m_1+m_2-1}q+\cdots+\lambda_0\,q^{m_1+m_2}
\]
with $\deg(p)=\deg(q)=ab$. If $m_1+m_2>1$, then $h(p,q)$ is always 
reducible. Since we assume $\Omega$ is irreducible, we obtain
that $m_1+m_2=1$. Then, after a linear transformation, we can assume 
that $h(x,y)=x^b+y^a$ (which is irreducible for $\gcd(a, b)=1$). 
Since $1\leq a\leq b\leq c$ and ${\rm gcd}(a,b)=1$, we have 
$1\leq a<b$. Since $2b< a+b+c=ab$, we have $a>2$. Note that 
$c=ab-a-b=(a-1)b-a$. Thus, we have $c>b$. Thus, we obtain (1b) 
when ${\rm gcd}(a,b)=1$. Therefore (1b) holds by lifting 
to the general case when ${\rm gcd}(a,b)> 1$.

One can easily show that the conditions in (1b) imply that
$\kk[x,y]_c=0$.

(2) $(2b)\Rightarrow (2a)$: By assumption, $\Omega$ is 
reducible. By the proof of part (1), $\deg(\Omega)=ka+lb$ 
for some $k, l\in \mathbb N$. If $k, l\ge 1$, we can let 
$\delta=x^{k-1}y^{l-1}\frac{\partial}{\partial z}$. 
Otherwise, we may assume $k=0$ and $\Omega=y^l$ (as 
$\Omega$ is reducible). Then we let 
$\delta=z\frac{\partial}{\partial z}-x\frac{\partial}{\partial x}$. 
Then (2a) holds.

$(2a)\Rightarrow (2b)$: Assume to the contrary that $\Omega$ is 
irreducible. Without loss of generality, let $\Omega=x^b+y^a$ 
with $\gcd(a,b)=1$ as in part (1). By (1b), we have $m_1+m_2=1$ and 
$a<b<c=ab-a-b$. This (together with $\kk[x,y]_c=0$) implies 
that any graded derivation $\delta$ of $A$ of degree zero must 
have the form $\delta(x)=\alpha x$, $\delta(y)=\beta y$ and 
$\delta(z)=\gamma z$ for some $\alpha,\beta,\gamma\in \kk$. 
So $\divv(\delta)=\delta(\Omega)=0$ yields that $\delta=0$.
This finishes the proof.
\end{proof}

\begin{proposition}
\label{zzpro2.3}
Let $A=\kk[x, y, z]$ be a weighted polynomial algebra with 
$\deg x=1, \deg y=1, \deg z=2$. Then the nonzero homogeneous 
degree $4$ polynomials $\Omega\in A$ can be classified in 
\Cref{tab:112}. In particular, $\Omega$ has an isolated singularity 
if and only if $\Omega=z^2+xy^3+\lambda x^2y^2+x^3y$ with $\lambda\ne \pm 2$ up to graded isomorphisms of $A$.
\end{proposition}

\begin{proof}
Since $\deg(\Omega)=4$, we have 
$\Omega=l_1z^2+l_2zg(x,y)+h(x,y)$, where $\deg g(x,y)=2$, 
$\deg h(x,y)=4$ and $l_1, l_2\in \kk$. If $l_1\ne 0$, then 
after a linear transformation of $z$, we can assume that 
$\Omega=z^2+h(x,y)$. If $l_1=0$ and $l_2\ne 0$, we 
can assume that $g(x,y)=x^2$ or $xy$ after a further 
linear transformation of $x$ and $y$. So, we only need to 
consider the following cases.

\smallskip

\noindent{\bf Case 1:} $\Omega= z^2+h(x,y)$.

If $h(x,y)=0$, then $\Omega=z^2$. If $0\neq h(x,y)$ has a root of multiplicity $4$ in $\mathbb P^1$, then without loss of 
generality, we can assume that $h(x,y)=x^4$.  If $h(x,y)$ has a root of multiplicity $3$, then we can assume that 
$h(x,y)=x^3y$ due to the symmetry between $x$ and $y$. If 
$h(x,y)$ has  a root of multiplicity $2$, then we can assume that $h(x,y)=x^2y(y+\lambda x)$ for some $\lambda=0$ or $1$. 
If $h(x,y)$ has no repeated root, then we can assume that 
$h(x,y)=xy(y+x)(y+kx)$ for some $k\in \kk\setminus\{0,1\}$. 
Then $h(x,y)=xy^3+(k+1)x^2y^2+kx^3y$. By a suitable 
re-scaling of $x$ and $y$, we obtain 
$h(x,y)=xy^3+\lambda x^2y^2+x^3y$ for some $\lambda\in\kk$.

\noindent{\bf Case 2:} $\Omega=x^2z+h(x,y)$.

After a linear transformation of $z$ and re-scaling of 
$x,y$, and $z$ if necessary, we can assume that 
$\Omega=x^2z+\lambda_1xy^3+\lambda_2y^4$ for some 
$\lambda_1, \lambda_2\in \{0,1\}$.

\noindent{\bf Case 3:} $\Omega=xyz+h(x,y)$.

After a linear transformation of $z$ and re-scaling of 
$x,y$ and $z$ if necessary, we can have 
$\Omega=xyz+\lambda_1x^4+\lambda_2y^4$ for some 
$\lambda_1, \lambda_2\in \{0,1\}$. 

\noindent{\bf Case 4:} $\Omega=h(x,y)$.

If $\Omega=h(x,y)$, then by the same argument of Case 1, 
we can show that $\Omega$ is one of the following forms: 
\[
x^4, \quad x^3y, \quad x^2y^2, \quad x^2y^2+ x^3y, \quad 
xy^3+\lambda x^2y^2+x^3y\,\,\, \text{for some} \,\lambda\in \kk.
\]
By direct computation, we can verify that 
$z^2+xy^3+\lambda x^2y^2+x^3y$ with $\lambda \ne \pm 2$ has 
an isolated singularity. 
\end{proof}

\begin{proposition}
\label{zzpro2.4}
Let $A=\kk[x, y, z]$ be a weighted polynomial algebra with 
$\deg x=1, \deg y=2, \deg z=3$. Then the nonzero homogeneous 
degree $6$ polynomials $\Omega\in A$ can be classified in 
\Cref{tab:123}. In particular, $\Omega$ has an isolated singularity if and only if $\Omega=z^2+y^3+\lambda x^2y^2+x^4y$ 
with $\lambda\ne \pm2$ up to graded isomorphisms of $A$.
\end{proposition}

\begin{proof}
Since $\deg (\Omega)=6$, we have 
$\Omega=l_1z^2+l_2zg(x,y)+h(x,y)$, where $g(x,y)\in\Bbbk[x,y]$ 
and $h(x,y)\in\Bbbk[x,y]$ have degrees $3$ and $6$, 
respectively, and $l_1, l_2\in \kk$. If $l_1\ne0$, then by a 
linear transformation of $z$, we can assume that 
$\Omega=z^2+h(x,y)$. If $l_1=0$ and $l_2\ne 0$, by a possible 
linear transformation of $x$ and $y$, we can have $g(x,y)=x^3$ 
or $g(x,y)=xy$. We write, in general, 
\[h(x,y)=w_1y^3+w_2x^2y^2+w_3x^4y+w_4x^6,\]
where $w_i\in \kk$ for $1\leq i\leq 4$.
So, we only need to consider the following cases.

\noindent{\bf Case 1:} $\Omega=z^2+h(x,y)$.

\medskip

\noindent{\bf Subcase 1:} If $w_1\ne 0$, then we can write 
$h(x,y)=(y+ax^2)(y+bx^2)(y+cx^2)$ for $a, b, c\in \kk$. 
After a linear transformation of $y$, we can assume 
$h(x,y)=y(y+ax^2)(y+bx^2)$. By a possible re-scaling of $x,y$ and $z$, 
we can assume that $\Omega$ is one of the following forms: 
\[
z^2+y^3,\quad z^2+y^3+x^2y^2, \quad z^2+y^3
+\lambda x^2y^2+x^4y \,\,\, \text{for some} \, \lambda \in \kk.
\]

\noindent{\bf Subcase 2:} If $w_1=0$ and $w_2\neq 0$, then, 
similarly, we can assume that $h(x,y)=x^2(y+ax^2)(y+bx^2)$ 
for some $a,b\in \kk$. A further linear transformation of 
$x,y$ and $z$ yields $\Omega=z^2+x^2y^2$ or $\Omega=z^2+x^2y^2+x^4y$.

\noindent{\bf Subcase 3:} 
If $w_1=w_2=0$ and $w_3\neq 0$, then we can assume that 
$h(x,y)=x^4(y+ax^2)$ for some $a\in \kk$. A linear 
transformation of $y$ yields $\Omega=z^2+x^4y$.

\noindent{\bf Subcase 4:} 
If $w_1=w_2=w_3=0$ and $w_4\neq 0$, then by a re-scaling 
of $x$, we get $\Omega=z^2+x^6$.

\noindent{\bf Subcase 5:} 
Finally, if $w_1=w_2=w_3=w_4=0$, then we have $\Omega=z^2$.

\medskip

\noindent{\bf Case 2:} $\Omega=x^3z+h(x,y)$.

After a linear transformation of $z$, we can assume that 
$\Omega=x^3z+\lambda_1x^2y^2+\lambda_2y^3$ for some 
$\lambda_1, \lambda_2\in \kk$. It is easy to check that 
$\Omega$ is one of the following forms:
\[
x^3z, \quad x^3z+y^3, \quad x^3z+x^2y^2, \quad x^3z+x^2y^2+y^3.
\]

\noindent{\bf Case 3:} $\Omega=xyz+h(x,y)$. 

Again, via a linear transformation of $z$, one can assume that  
\[
\Omega=xyz+\lambda_1x^6+\lambda_2y^3
\]
for some $\lambda_1, \lambda_2\in \kk$. After re-scaling 
$x$ and $y$ as needed, we can assume that $\Omega$ is of 
one of the following forms:
\[
xyz, \quad xyz+x^6, \quad xyz+y^3, \quad xyz+x^6+y^3.
\]

\noindent{\bf Case 4:} 
If $\Omega=h(x,y)$, then by the same argument as in Case 1, 
$\Omega$ can be assumed to be one of the following forms:
\[
y^3, \quad y^3+x^2y^2,\quad y^3+\lambda x^2y^2+x^4y,\quad x^4y,
\quad x^2y^2, \quad x^2y^2+x^4y,\quad x^6
\]
where $\lambda\in\kk$.

By a direct computation, we can further verify that 
$\Omega=z^2+y^3+\lambda x^2y^2+x^4y$ has an isolated 
singularity if and only if $\lambda\ne \pm2$.
\end{proof}

\begin{theorem}
\label{zzthm2.5}
Let $A=\kk[x, y, z]$ be a weighted polynomial algebra with 
$\deg(x)=a, \deg (y)=b, \deg(z)=c$ for $1\leq a\leq b\leq c$. 
Let $\Omega$ be a nonzero homogeneous polynomial of degree 
$a+b+c$. Then, up to a graded automorphism of $A$, we have 
the following: 
\begin{enumerate}
\item[(1)] \cite{BM, DH, DML, KM, LX}
If $a=b=c$, then $\Omega$ is one of the forms listed in 
\Cref{tab:111}.
\item[(2)]  
If $a=b<c$, then $\Omega$ is one of the forms listed in 
\Cref{tab:112} and 
\Cref{tab:a=b<c}.
\item[(3)]  
If $a<b=c$, then every $\Omega$ is reducible and is one of 
the forms listed in \Cref{tab:a<b=c}. 
\item[(4)]  If $a<b<c$, then $\Omega$ is one of the forms 
listed in \Cref{tab:123} and \Cref{tab:abc2}. 
\end{enumerate}
\end{theorem}

\begin{proof}
$(1)$ Since $a=b=c$, we can reduce the classification of $\Omega$ 
to the case where the degrees of $x, y$ and $z$ are all equal to 
$1$. In this case, the classification of $\Omega$ is well-known. 
Also see \cite[Corollary 6.7]{TWZ}.
\smallskip 

$(2)$ Since $a=b<c$, then $\deg (\Omega)=2a+c<3c$. So we can 
write $\Omega= z^2f(x,y)+zg(x,y)+h(x,y)$, where 
$\deg(f(x,y))=a+b-c<a$, $\deg(g(x,y))=2a$ and 
$\deg(h(x,y))=2a+c$. 

If $a\nmid c$, then in particular, $c\ne a+b$, whence $f(x,y)=0$. 
If $h(x,y)\neq 0$, then $a\mid \deg (h(x, y))=2a+c$, we get 
$a\mid c$, yielding a contradiction. So $\Omega=zg(x,y)$ is 
reducible. By a linear transformation of $x,y$, we can assume 
that $\Omega=xyz$ or $x^2z$.  

If $c=ka$ for some integer $k\geq 2$, then $(a, b, c)=(a, a, ka)$. 
If $k=2$, then the result is given by Proposition \ref{zzpro2.3}. 
Now, assume $k>2$. Then we have $\Omega=zg(x,y)+h(x,y)$. We can 
assume $g(x,y)=x^2$ or $xy$ after a linear transformation of 
$x,y$. If $\Omega=x^2z+h(x,y)$, after a necessary linear 
transformation of $z$, we can write  
$\Omega=x^2z+\lambda_1xy^{k+1}+\lambda_2y^{k+2}$ for 
$\lambda_1, \lambda_2\in \{0,1\}$. If $\Omega=xyz+h(x,y)$, then 
similarly, we can write 
$\Omega=xyz+\lambda_1x^{k+2}+\lambda_2y^{k+2}$ for some 
$\lambda_1, \lambda_2\in \{0,1\}$. 

\smallskip 

$(3)$ Since $a<b=c$, we have that $\deg (\Omega)=a+2b<3b$. If 
$a\mid b$, then we have 
$\Omega=\lambda x^{1+2\frac{b}{a}}+ x^{1+\frac{b}{a}}g(y,z)+xf(y,z)$, 
where $\lambda\in \kk, \deg (f(y,z))=2b$ and $\deg (g(y,z))=b$. 
Thus $\Omega$ is reducible. If $a\nmid b$, then $\Omega=xf (y, z)$ or $\Omega=\lambda x^{1+b}+xf(y, z)$ when $a=2$, which is again reducible.

If $a\nmid b$ and $a\neq 2$, after a linear transformation of $y, z$, we can 
assume that $\Omega=xyz$ or $xy^2$. If $a\nmid b$ and $a=2$, we can have $\Omega=x^{1+b}+xyz$ or $\Omega=x^{1+b}+xz^2$ after a linear transformation for $x, y$ and $z$. If $a\mid b$, we can assume 
that $\Omega=x\Omega_1$ where $\Omega_1=\lambda u^2+ug(y,z)+f(y,z)$ 
with $u=x^{\frac{b}{a}}$ for some $\lambda \in \kk$.  We can 
rewrite $\Omega_{1}$ as follows
\[\Omega_1=k_1z^2+k_2zh_1(y,u)+h_2(y,u)
\]
for some $k_1, k_2\in \kk$ and $h_1(y,u), h_2(y,u)\in \kk[y,u]$. 
If $k_1\ne 0$, then we can assume that $\Omega$ is one of the 
following forms:
\[
xz^2, \quad xz^2+xy^2, \quad xz^2+x^{1+\frac{2b}{a}},\quad 
xz^2+xy^2+x^{1+\frac{2b}{a}}, \quad x^{1+\frac{b}{a}}y+xz^2.
\]
If $k_1=k_2=0$, then $\Omega=xy^2, x^{1+\frac{2b}{a}}, 
x^{1+\frac{b}{a}}y, xy^2+x^{1+\frac{2b}{a}}$. If $k_1=0$ 
and $k_2\ne 0$, then $\Omega=x^{1+\frac{b}{a}}z, 
x^{1+\frac{b}{a}}z+xy^2$, $xyz, xyz+x^{1+\frac{2b}{a}}$.

\smallskip

$(4)$ If $\Omega\in \kk[x,y]$, by \Cref{zzlem2.2}, the irreducible ones are given by $\Omega=x^{b/d}+y^{a/d}$ where $d=\gcd(a, b)$ and $2<a<b<c=ab/d-a-b$. Moreover, such irreducible $\Omega$ won't occur unless $c=ma+nb, c\neq a+b, a\nmid b$ for some $m, n\in \mathbb{Z}$. Let us assume $\Omega\not\in \kk[x,y]$ in the remaining argument. 
Since $a<b<c$, we have that $\deg (\Omega)=a+b+c<3c$. Thus we 
can assume that $\Omega=z^2f(x,y)+zg(x,y)+h(x,y)$, where 
$\deg(f(x,y))=a+b-c<a$, $\deg(g(x,y))=a+b$ and 
$\deg(h(x,y))=a+b+c$. We divide the argument into two cases.

{\bf Case 1:} $c=ma+nb$ for some integers $m$ and $n$. 

{\bf Subcase 1:} If $c=a+b$, then we have $\deg (\Omega)=2c$ 
and $f(x,y)\in \kk$. As a result, we can assume that 
$\Omega=z^2+h(x,y)$ or $zg(x,y)+h(x,y)$. Since 
$\deg(h(x,y))=2a+2b<4b$, we can write 
$h(x,y)=h_3(x)y^3+h_2(x)y^2+h_1(x)y+h_0(x)$ where $h_i(x)$ is 
a monomial in $x$ of degree $2a+(2-i)b$ for $0\leq i\leq 3$. 
Let us further assume that $a\nmid b$. After a linear 
transformation of $x,y$, we can have $g(x,y)=xy$ if 
$g(x, y)\neq 0$. If $a\nmid 2b$, then $h(x,y)=\lambda x^2y^2$ 
for some $\lambda \in \kk$.  Hence, by a linear transformation 
of $z$, $\Omega$ can be one of the following possible forms:
\[
z^2, \quad z^2+x^2y^2,\quad xyz.
\]
If $a\mid 2b$, then $h(x,y)$ can be one of the following forms:
\[
x^2y^2,\quad x^{2+\frac{2b}{a}}, \quad x^2y^2+ x^{2+\frac{2b}{a}}.
\]
By a linear transformation of $z$, $\Omega$ can be one of the 
following possible forms:
\[
z^2, \quad z^2+x^2y^2, \quad z^2+ x^{2+\frac{2b}{a}},\quad 
z^2+x^2y^2+x^{2+\frac{2b}{a}}, \quad xyz, \quad xyz+x^{2+\frac{2b}{a}}.
\]
Now we assume that $a\mid b$. If $b=2a$, then $c=3a$ and so 
see \Cref{zzpro2.4}. We next consider the case $b\neq 2a$. 
Since $a+b+c=2a+2b<4b$ and $b\neq 2a$, we can write 
$h(x,y)=\lambda_1x^2y^2+\lambda_2x^{2+\frac{b}{a}}y
+\lambda_3x^{2+2\frac{b}{a}}$ for some $\lambda_1, \lambda_2, 
\lambda_3\in\kk$. If $\Omega=z^2+h(x,y)$, then after a linear 
transformation of $x, y$, 
$\Omega$ is one of the following forms:
\[
z^2, \quad z^2+x^2y^2, \quad z^2+x^{2+\frac{b}{a}}y, \quad 
z^2+x^{2+\frac{2b}{a}}, \quad z^2+x^2y^2+x^{2+\frac{b}{a}}y.
\]
If $\Omega=zg(x,y)+h(x,y)$, then after a linear transformation 
of $x,y$, we can have $g(x,y)=xy$ or $x^{1+\frac{b}{a}}$.   
By a linear transformation of $z$, $\Omega$ is one of the following:
\[
xyz,\quad xyz+x^{2+\frac{2b}{a}}, \quad x^{1+\frac{b}{a}}z,\quad 
x^{1+\frac{b}{a}}z+x^2y^2.
\]

{\bf Subcase 2:} If $c\neq a+b$, then we can assume that 
$\Omega=g(x,y)z+h(x,y)$ where $\deg(g(x,y))=a+b$ and 
$\deg(h(x,y))=(m+1)a+(n+1)b$. Let $k, l$ be the possible 
integers such that $k=(n+1)b/a$ and $l=(m+1)a/b$. Note that 
$kl=(m+1)(n+1)$. If $a\nmid b$, then we can assume that 
$g(x,y)=xy$ if $g(x,y)\neq 0$. If $a\mid b$, then we have 
$g(x,y)=xy$ or $x^{1+\frac{b}{a}}$ if $g(x,y)\neq 0$. By a 
linear transformation of $z$, $\Omega$ can be reduced to one 
of the following forms:
\[
xyz, \quad xyz+x^{m+1+k}, \quad xyz+y^{n+1+l}, \quad 
xyz+x^{m+1+k}+y^{n+1+l},\quad x^{1+\frac{b}{a}}z
\]
and 
\[
x^{1+\frac{b}{a}}z+x^{\frac{b}{a}}y^{n+l},\quad 
x^{1+\frac{b}{a}}z+x^{\frac{b}{a}}y^{n+l}+ y^{1+n+l}, 
\quad x^{1+\frac{b}{a}}z+ y^{1+n+l}
\]
where $k$ and $l$ are assumed to be integers in the first place.

{\bf Case 2:} If $c\ne ma+nb$ for any $m,n\in \mathbb Z$, then 
$f(x,y)=0$. Suppose $h(x,y) \ne 0$. Then 
$a+b+c=\deg (h(x,y))=sa+lb$ for some $s,l\in \mathbb N$. Then 
$c=(s-1)a+(l-1)b$, which is a contradiction. So $h(x,y)=0$. Then 
$\Omega=zg(x,y)$ is reducible. After a linear transformation of 
$x, y$, we can assume that $\Omega=xyz$ for $a\nmid b$ and 
$\Omega= xyz$ or $x^{1+\frac{b}{a}}z$ for $a\mid b$.

This completes the proof.
\end{proof}

\section{Poisson fraction fields and automorphism groups}
\label{zzsec3}

In this section, we discuss Poisson fraction fields and 
Poisson automorphism groups related to unimodular Poisson 
algebras in dimension three. Throughout this section, we will assume \Cref{zzhyp0.1} unless mentioned otherwise.

\begin{definition}
\label{zzdef3.1}
We define a map $\pi: A\to \mathfrak X^2(A)$ such that, for 
each $h\in A$, $\pi_{h}\colon =\pi(h)$ is defined by
\begin{equation}
\label{E3.1.1}\tag{E3.1.1}
\pi_h(f,g)=\det \begin{pmatrix} 
f_{x} & f_{y} & f_{z}\\
g_{x} & g_{y} & g_{z}\\
h_{x} & h_{y} & h_{z}
\end{pmatrix}
(=
\det\left(\frac{\partial\left(f,\,g,\,h\right)}
{\partial\left(x,\,y,\,z\right)}\right)=\{f,g\}_h)
\quad 
\text{for all $f,g\in A:=\kk[x,y,z]$}.
\end{equation}
Note that the definition of $\pi_h$ depends on the 
generating set $\{x,y,z\}$. It is clear that if a new set of generators $\{x,y,z\}$ is used, then the corresponding $\pi_{h}$ will be a scalar multiple of the original $\pi_{h}$.  
\end{definition}

One can check that $[\pi_h,\pi_h]_S=0$. In particular, $\pi_\Omega$ is a connected graded Poisson bracket on $A$ if $\Omega$ satisfies \Cref{zzhyp0.1}. To keep things simple, we will introduce the following. 

\begin{notation} Here $w=n-(a+b+c)$ may not be necessarily zero.
\label{zznot3.2}
\
\begin{itemize}
\item 
$A_\Omega:=(A,\pi_\Omega)$: the unimodular Poisson algebra with 
potential $\Omega$.
\item 
${\rm Aut}_{gr}(A)$: the group of all graded algebra automorphisms 
of $A$.
\item 
${\rm Aut}_{P}(A)$ {\rm{(}}resp. ${\rm Aut}_{grP}(A)${\rm{)}}: 
the group of all {\rm{(}}resp. graded{\rm{)}} Poisson 
automorphisms of $A=A_\Omega$. 
\item 
$P_{\Omega-\xi}:=A/(\Omega-\xi)$ {\rm{(}}$\xi\in \kk${\rm{)}}: 
the quotient Poisson algebra.
\item 
${\rm Aut}_{P}(P_\Omega)$ {\rm{(}}resp. 
${\rm Aut}_{grP}(P_\Omega)${\rm{)}}: the group of all 
{\rm{(}}resp. graded{\rm{)}} Poisson automorphisms of $P_\Omega$.
\item 
$Q(P_{\Omega-\xi})$: the Poisson fraction field of 
$P_{\Omega-\xi}$ whenever it is an integral domain. 
\end{itemize}
\end{notation}

Note that any graded unimodular Poisson structure on $A$ is 
determined by such a potential $\Omega$. If $\Omega$ is 
homogeneous of degree $|x|+|y|+|z|$, then  
$P_{\Omega-\xi}\cong P_{\Omega-1}$ whenever $\xi\neq 0$. In 
this case, we can always assume $\xi$ to be either $0$ or $1$.

\begin{lemma}
\label{zzlem3.3}
Let $\Omega$ and $\Omega'$ be two homogeneous potentials for $A=\kk[x,y,z]$ of degrees $>\max\{|x|,|y|,|z|\}$. Then the  two Poisson algebras $A_\Omega$ and $A_{\Omega'}$ are isomorphic if and only if there is an algebra 
automorphism $\phi$ of $A$ such that 
$\Omega'=\phi(\Omega)/\det(\phi)$ where 
$\det(\phi)=\det\left(\frac{\partial\left(\phi(x),\,\phi(y),
\,\phi(z)\right)}{\partial\left(x,\,y,\,z\right)}\right)$. 
As a consequence,
\begin{align*}
{\rm Aut}_{P}(A_\Omega)~=~\{\phi\in {\rm Aut}(A)\,|
\, \phi(\Omega)=\det(\phi)\Omega\},
\end{align*}
and ${\rm Aut}_P(A_\Omega)={\rm Aut}_{P}(A_{\lambda\Omega})$ 
for any $\lambda\in \kk^\times$. 
\end{lemma}

\begin{proof}
Suppose there is a Poisson algebra 
isomorphism $\phi: A_\Omega\to A_{\Omega'}$. For any polynomials $f, g\in A$, $\phi(f)$ and $\phi(g)$ can 
be regarded as polynomials of variables $\phi(x), \phi(y)$ and 
$\phi(z)$ in $\phi(A)=A$. Therefore, according to \eqref{E3.1.1}, 
we have 
\begin{align*}
\{\phi(f),\phi(g)\}_{\Omega'}
&=\phi(\{f,g\}_{\Omega})=
\phi\left(\det(\frac{\partial(f,\,g,\,\Omega)}{\partial(x,y,z)})\right)=
\det\left(\frac{\partial\left(\phi(f),\,\phi(g),\,
\phi(\Omega)\right)}{\partial\left(\phi(x),\,\phi(y),\,\phi(z)\right)}\right)
\\
&=\det\left(\frac{\partial
\left(x,\,y,\,z\right)}
{\partial\left(\phi(x),\,\phi(y),\, \phi(z)\right)}\right)
\det\left(\frac{\partial\left(\phi(f),\,\phi(g),\,
\phi(\Omega)\right)}{\partial\left(x,\,y,\,z\right)}\right)\\
&=\det(\phi)^{-1}\{\phi(f),\phi(g)\}_{\phi(\Omega)}=\{\phi(f),\phi(g)\}_{\det(\phi)^{-1}\phi(\Omega)}. 
\end{align*}
This implies that 
\[\{-,-\}_{\Omega'}=\{-,-\}_{\det(\phi)^{-1}\phi(\Omega)}\] 
which yields the same Poisson bracket on $A$. Thus $\phi(\Omega)=\det(\phi)\Omega'+\lambda$ for some scalar $\lambda \in\kk$. Since $\Omega$ and $\Omega'$ are homogeneous of degrees $>\max\{|x|,|y|,|z|\}$, their partial derivatives are all homogeneous of positive degree. Therefore, one can check that $A/(\det(\phi)\Omega'+\lambda)$ is regular for $\lambda \neq 0$ while $A/(\Omega)$ is not regular by the Jacobian criterion. As a consequence, the induced algebra isomorphism $\phi: A/(\Omega)\to A/(\phi(\Omega))=A/(\det(\phi)\Omega'+\lambda)$ implies that $\lambda=0$. So $\Omega'=\phi(\Omega)/\det(\phi)$.   

Conversely, suppose there is an algebra isomorphism $\phi$ of $A$ such that $\Omega'=\phi(\Omega)/\det(\phi)$. Then we have
\begin{align*}
\phi(\{f,g\}_{\Omega})&=
\phi\left(\det(\frac{\partial(f,\,g,\,\Omega)}{\partial(x,y,z)})\right)=
\det\left(\frac{\partial\left(\phi(f),\,\phi(g),\,
\phi(\Omega)\right)}{\partial\left(\phi(x),\,\phi(y),\,\phi(z)\right)}\right)= \det(\phi)\det\left(\frac{\partial\left(\phi(f),\,\phi(g),\,\Omega'\right)}{\partial\left(\phi(x),\,\phi(y),\,\phi(z)\right)}\right)\\
&=\det(\phi)\det\left(\frac{\partial
\left(x,\,y,\,z\right)}
{\partial\left(\phi(x),\,\phi(y),\, \phi(z)\right)}\right)
\det\left(\frac{\partial\left(\phi(f),\,\phi(g),\,
\Omega'\right)}{\partial\left(x,\,y,\,z\right)}\right)\\
&=\det(\phi)\det(\phi)^{-1}\{\phi(f),\phi(g)\}_{\Omega'}=\{\phi(f),\phi(g)\}_{\Omega'}.
\end{align*}
So, $\phi$ is indeed a Poisson isomorphism. Finally, the 
consequences follow immediately.
\end{proof}

Now, we can classify all connected graded unimodular Poisson algebras in dimension three. 

\begin{theorem}
\label{zzthm3.4}
Any connected graded unimodular Poisson algebra $A$ is 
isomorphic to some $A_{\lambda \Omega}$, where $\Omega$ 
is listed in \Cref{zzthm2.5} and $\lambda\in \kk^\times$. 
\end{theorem}

\begin{proof}
By \cite[Theorem 5]{Pr}, any unimodular Poisson structure on $A=\kk[x,y,z]$ is given by $\pi_\Omega$ for some potential $\Omega$. Moreover, $\pi_\Omega$ being graded implies $|\Omega|=|x|+|y|+|z|$. Hence our classification follows 
from \Cref{zzdef3.1}, \Cref{zzthm2.5}, and \Cref{zzlem3.3}.
\end{proof}

\begin{example}
\label{zzexa3.5}
We introduce some examples of Poisson fields of transcendental 
degree two over the field $\kk$.
\begin{itemize}
\item[(1)] 
We define $K_{Weyl}:=\kk(x,y)$ to be the Poisson field 
$\kk(x, y)$ with the Poisson bracket $\{x,y\}=1$.
\item[(2)] 
We define $K_q:=\kk(x,y)$ to be the Poisson field $\kk(x, y)$ 
with the Poisson bracket $\{x,y\}=qxy$ for some 
$q\in \kk^\times$. It is shown in \cite[Corollary 5.4]{GLa} 
that $K_q\cong K_{p}$ as Poisson fields if and only if 
$p=\pm q$. Note that $K_{Weyl}$ is not isomorphic to any 
$K_{q}$ as Poisson fields by \cite[Corollary 5.3]{GLa}. 
\item[(3)] 
Consider the irreducible cubic polynomial 
\[\Omega_{\zeta,\lambda}:
=\zeta\left(x^3+y^3+z^3\right)+\lambda xyz\]
with two parameters $\zeta,\lambda\in \kk$ such that 
$\zeta\neq 0,\lambda^3\neq -3^3\zeta^3$. We denote the 
corresponding graded unimodular Poisson algebra by 
$A_{\Omega_{\zeta,\lambda}}$,  where the Poisson structure on 
$\kk[x,y,z]$ is defined by
\[
\{x,y\}=3\zeta z^2+\lambda xy,\ \{y,z\}
=3\zeta x^2+\lambda yz,\ \{z,x\}=3\zeta y^2+\lambda xz
\]
with $\deg(x)=\deg(y)=\deg(z)=1$. Let 
$S_{\zeta,\lambda}:=Q(P_{\Omega_{\zeta,\lambda}})$ be 
the Poisson fraction field. By \cite[Corollary 6.4]{HTWZ1}, 
$S_{\zeta, \lambda}$ is not isomorphic to $K_{Weyl}$ or 
$K_q$. However, $S_{\zeta,\lambda}=F_{\zeta,\lambda}(t)$ 
as fields, where $F_{\zeta,\lambda}$ is the function field of 
the elliptic curve  $\Omega_{\zeta,\lambda}=0$. Suppose 
$S_{\zeta,\lambda}\cong S_{\zeta',\lambda'}$ as Poisson 
fields, whence they are isomorphic as function fields. By 
\cite[Theorem 2]{De}, we have $F_{\zeta,\lambda}\cong F_{\zeta',\lambda'}$ 
as function fields. As a result, the two elliptic curves 
$\Omega_{\zeta,\lambda}=0$ and $\Omega_{\zeta',\lambda'}=0$ are 
birationally equivalent or have the same $j$-invariant 
$27\frac{(\lambda/\zeta)^3 ((\lambda/\zeta)^3+8)^3}
{((\lambda/\zeta)^3-1)}=27\frac{(\lambda'/\zeta')^3 
((\lambda'/\zeta')^3+8)^3}{((\lambda'/\zeta')^3-1)}$ 
(for the $j$-invariant of a Hesse form of a smooth elliptic 
curve, see \cite[Theorem 2.11]{Frium}). It is unclear if two elliptic curves with equations $\Omega_{\zeta,\lambda}=0$ and $\Omega_{\zeta',\lambda'}=0$ being birationally equivalent necessarily implies an isomorphism between their corresponding Poisson fields $S_{\zeta,\lambda}$ and $S_{\zeta',\lambda'}$.
\end{itemize}
\end{example}

\begin{theorem}
\label{zzthm3.6}
Let $A_\Omega$ be a connected graded unimodular Poisson algebra 
defined by some homogeneous irreducible potential $\Omega$. Let 
$Q=Q(P_\Omega)$ be the Poisson fraction field of $P_{\Omega}$. 
Then the following hold:   
\begin{itemize}
\item[(1)] 
If $\Omega$ does not have an isolated singularity, then $Q$ is 
isomorphic to either $K_{Weyl}$ or  $K_q$ for some $q \in \kk^\times$.
\item[(2)] 
If $\Omega$ has an isolated singularity, then $Q$ is isomorphic 
to $S_{\zeta,\lambda}$ for some $\zeta,\lambda\in \kk$ with 
$\zeta \neq 0$ and $\lambda^3\neq -3^3\zeta^3$.
\end{itemize}
Moreover, we label those $\Omega$ with the corresponding $Q$ that are isomorphic to 
$S_{\zeta,\lambda}$ by 
\raisebox{.5pt}{\textcircled{\raisebox{-.5pt} {i}}} and 
those that are isomorphic to $K_q$ by 
\raisebox{.5pt}{\textcircled{\raisebox{-.5pt} {q}}} and 
those that are isomorphic to $K_{Weyl}$ by 
\raisebox{.5pt}{\textcircled{\raisebox{-.5pt} {w}}}  
in \Cref{tab:111}--\Cref{tab:abc2}. 
\end{theorem}

\begin{proof}
We conduct a case-by-case verification for those irreducible 
potentials $\Omega$ listed in \Cref{tab:111}--\Cref{tab:abc2}. 
Then the result follows from \Cref{zzthm3.4}.

(1) Suppose $\Omega$ does not have an isolated singularity. 
As an illustration, we provide some details when 
$\Omega=z^2+y^3+2x^2y^2+x^4y$ ($\deg x=1, \deg y =2$ and 
$\deg z=3$) in \Cref{tab:123} and 
$\Omega=z^2+x^2y^2+x^{2+\frac{2b}{a}}$ with $a\nmid b$ 
($\deg x=a>2, \deg y =b$ and $\deg z=c$) in \Cref{tab:abc2}. 
We also check for $\Omega=z^2+x^{2+\frac{2b}{a}}$ when 
$c=a+b$ and $a\nmid b$ in \Cref{tab:abc2}.  
Let $Q=Q(P_{\lambda\Omega})$ for some $\lambda\in \kk^\times$.

If $\Omega=z^2+y^3+2x^2y^2+x^4y$, then we have $z^2+y(y+x^2)^2=0$ 
in $P_\Omega$. Let $w=\frac{z}{y+x^2}$. Then $y=-w^2$ and 
$Q=\kk(x,w)$. Consider $\{x, w\}=\{x, \frac{z}{y+x^2}\}
=-\lambda(y+x^2)=\lambda(w^2-x^2)$. After a linear 
transformation of $x$ and $w$, we have $\{x, w\}=pxw$ for 
some $p\in \kk^{\times}$. In this case, 
$Q(P_{\lambda\Omega})\cong K_{p}$.

If $\Omega=z^2+x^2y^2+x^{2+\frac{2b}{a}}$, then $\frac{2b}{a}$ 
is odd since $a\nmid b$, say $2t+1$. In $Q$, we have 
$x=-(\frac{z}{x^{1+t}})^2-(\frac{y}{x^{t}})^2$. Set 
$u=\frac{z}{x^{1+t}}$ and $v=\frac{y}{x^t}$. Then $Q=\kk(u,v)$. 
One can check that $\{u,v\}=\lambda(u^2+v^2)$. After a linear 
transformation of $u$ and $v$, we can obtain $\{u,v\}=puv$ 
for some $p\in \kk^{\times}$.

If $\Omega = z^2+x^{2+\frac{2b}{a}}$ with $a\nmid b$ and $(a<b<c)\neq (1, 2, 3)$, then $\frac{2b}{a}$ is odd. Assume that $2+\frac{2b}{a}=2l+1$. Thus we have $z^2+x^{2l+1}=0$ and $x=-(\frac{z}{x^l})^2$ in $Q$. Set $s=y$ and $t=\frac{z}{x^l}$. We have $Q=\Bbbk(s, t)$ and 
\begin{eqnarray*}
\{s, t\}&=& \{y, \frac{z}{x^l}\}= \frac{1}{x^l} \{y, z\}-\frac{lz}{x^{l+1}} \{y, x\}= \frac{(2l+1)x^{2l}}{x^l}+\frac{2lz^2}{x^{l+1}}\\ &=&\frac{(2l+1) x^{2l+1} +2lz^2}{x^{l+1}}
=\frac{x^{2l+1}}{x^{l+1}}=x^{l}=(-t^2)^{l}=(-1)^l t^{2l}.
\end{eqnarray*}
Set $u=\frac{s}{(-1)^{l} t^{2l}}$ and $v=t$. Then we have $Q=\Bbbk(u, v)$ with $\{u, v\}=1$. Thus $Q$ is of Weyl type.

(2) Suppose $\Omega$ has an isolated singularity. By \Cref{zzthm2.5}, we need to consider the 
following three cases. Up to a scalar multiple, we can assume that 
\begin{enumerate}
\item[(a)] 
$\Omega_1=z^2x+yx^2+\lambda xy^2+y^3$ with $(a,b,c)=(1,1,1)$ 
and $\lambda\ne \pm 2$;
\item[(b)] 
$\Omega_2=z^2+x^3y+\lambda x^2y^2+xy^3$ with 
$(a, b, c)=(1,1,2)$ and $\lambda\ne \pm 2$; and
\item[(c)] 
$\Omega_3=z^2+y^3+\lambda x^2y^2+x^4y$ with 
$(a, b, c)=(1,2,3)$ and $\lambda\ne \pm 2$.
\end{enumerate} 
We use an alternative form in (a) following \cite[p. 255]{DML} 
instead of the Hessian normal form. We consider 
$Q_i=Q(P_{\mu\Omega_i})$ for some $\mu\in \kk^\times$ for 
$i=1,2,3$.

{\bf Case (a):}  In $Q_1$, we have 
\begin{equation*}
\{x,y\}= 2\mu zx,  \quad \quad 
\{y,z\}= \mu(z^2+2xy+\lambda y^2), \quad {\rm and} \quad 
\{z,x\}= \mu(x^2+2\lambda xy+3y^2) 
\end{equation*}
Denote $u=\frac{z}{x}$, $v=\frac{y}{x}$ and $w=x$. We have 
\[0=\frac{\Omega_1}{x^3}=\left(\frac{z}{x}\right)^2+\frac{y}{x}
+\lambda\left(\frac{y}{x}\right)^2+\left(\frac{y}{x}\right)^3
=u^2+v+\lambda v^2+v^3.\] 
A direct computation yields that 
\[\{u,w\}=\mu(w+2\lambda wv+3wv^2),\quad \{v,w\}=-2\mu wu,
\quad \{u,v\}=0.\]

{\bf Case (b):} In $Q_2$, we have 
\begin{equation*}
\{x,y\}= 2\mu z,  \,\,
\{y,z\}=\mu( 3x^2y+2\lambda xy^2+y^3), \,\, 
{\rm and}\,\, \{z,x\}= \mu(x^3+2\lambda x^2y+3xy^2).
\end{equation*}
Denote $u=\frac{z}{x^2}$, $v=\frac{y}{x}$ and $w=x$. We have 
\[0=\frac{\Omega_2}{x^4}=\left(\frac{z}{x^2}\right)^2+\frac{y}{x}
+\lambda\left(\frac{y}{x}\right)^2+\left(\frac{y}{x}\right)^3
=u^2+v+\lambda v^2+v^3.\] 
We can easily verify that 
\[\{u,w\}=\mu(w+2\lambda wv+3wv^2),\,\, \{v,w\}=-2\mu wu,\,\, 
\{u,v\}=0.\]
 
{\bf Case (c):} In $Q_3$, we have
\begin{equation*}
\{x,y\}= 2\mu z, \,\,
\{y,z\}= \mu(2\lambda xy^2+4x^3y),\,\, 
\{z,x\}= \mu(3y^2+2\lambda x^2y+x^4).
\end{equation*}
Denote $u=\frac{z}{x^3}$, $v=\frac{y}{x^2}$ and $w=x$. Then 
\[0=\frac{\Omega_3}{x^6}=\left(\frac{z}{x^3}\right)^2+
\left(\frac{y}{x^2}\right)^3+\lambda\left(\frac{y}{x^2}\right)^2
+\left(\frac{y}{x^2}\right)=u^2+v^3+\lambda v^2+v.\] 
It is routine to check that
\[\{u,w\}=\mu(w+2\lambda wv+3wv^2),\,\, \{v,w\}=-2\mu wu, \,
\,\{u,v\}=0.\]
Therefore, $Q_1\cong Q_2\cong Q_3$ as Poisson fraction 
fields. Finally, we reuse the Hesse form for (1) and note that 
any re-scaling of $\Omega$ can be written as 
$\Omega=\zeta(x^3+y^3+z^3)+\lambda xyz$ for some 
$\zeta,\lambda\in \kk$ such that $\zeta\neq 0, 
\lambda^3\neq -3^3\zeta^3$. The result follows from 
\Cref{zzexa3.5}(3).

\end{proof}

\begin{remark}
\label{zzrem3.7}
Let $\Omega$ be irreducible and denote the Poisson fraction 
field by $Q=Q(P_\Omega)$. Let $X$ be the projective curve in 
the weighted projective space $\mathbb P(a,b,c)$ determined 
by such $\Omega$. Then the statement of \Cref{zzthm3.6} can 
be refined as: if $X$ is smooth, then $Q\cong S_{\sigma,\lambda}$; 
if $X$ has nodal singularity, then $Q\cong K_q$ for some 
$q\in \Bbbk^{\times}$; and if  $X$ has cusp singularity, 
then $Q\cong K_{Weyl}$. 
\end{remark}

\begin{remark}
\label{zzrem3.8}
We note that the isomorphisms between $Q(P_\Omega)$ depend not on grading but on projective curves $\Omega=0$ when $\Omega$ is homogeneous of degree $|x|+|y|+|z|$. The geometry of elliptic curves appears to reflect properties of these Poisson algebras.
\end{remark}

We aim to investigate the relationship between the Poisson automorphism group of $A_{\Omega}$ and the type of $\Omega$. We determine the automorphism group of every connected graded unimodular Poisson algebra $A_\Omega$ where the potential $\Omega$ has an isolated singularity. 

\begin{theorem}
\label{zzthm3.9}
Let $A_\Omega$ be a connected graded unimodular Poisson algebra. 
If $\Omega$ has an isolated singularity, then every Poisson 
automorphism of $A$ is graded. As a consequence, 
\[{\rm Aut}_P(A_{\Omega})\cong {\rm Aut}_{grP}(A_{\Omega})\cong 
{\rm Aut}_{grP}(P_\Omega)\cong {\rm Aut}_{P}(P_\Omega).\]
\end{theorem}

We will establish a few lemmas before proving \Cref{zzthm3.9}. We'll use the Hessian normal forms of these elliptic curves in the weighted projective space $\mathbb P(a,b,c)$ to prove the lemmas ahead. These forms can be obtained from the $\Omega$ listed in part (2) of \Cref{zzthm3.6} by a linear transformation, as stated in the following Lemma.

\begin{lemma}
\label{zzlem3.10} 
Let $A=\kk[x, y, z]$ be a weighted polynomial algebra with 
$\deg(x)=a, \deg(y)=b, \deg (z)=c$ for some $1\leq a\leq b\leq c$ and ${\rm gcd}(a, b, c)=1$. 
If $\Omega$ is a homogeneous polynomial of degree $a+b+c$ that 
has an isolated singularity, then $\Omega$ is one of the forms
\begin{itemize}
\item[(1)] 
$\Omega_1=x^3+y^3+z^3+\lambda xyz$ for $(-\lambda)^3\neq 27$ 
and $(a,b,c)=(1,1,1)$.   
\item[(2)] 
$\Omega_2=x^4+y^4+z^2+\lambda xyz$ for $(-\lambda)^4\neq 64$ 
and $(a,b,c)=(1,1,2)$.
\item[(3)] 
$\Omega_3=x^6+y^3+z^2+\lambda xyz$ for $(-\lambda)^6\neq 432$ 
and $(a,b,c)=(1,2,3)$.  
\end{itemize}
\end{lemma}

\begin{proof}
Case (1) follows from a classical result that a smooth cubic projective curve is projectively equivalent to $x^3+y^3+z^3=3kxyz$ with $k^3\neq -1$, whose proof can be found in \cite[Theorem 2.12]{BM}, see also \cite{DML, DH, KM, LX}. We will show case (2) in detail and obtain the third case in a similar fashion. Replacing $z+\frac{1}{2}\lambda xy$ by $z$, one can rewrite 
$\Omega_2$ as $z^2+x^4+y^4-\frac{1}{4}\lambda^2 x^2y^2$, whence 
it becomes
\[z^2+(x^2-u y^2)(x^2-v y^2)=z^2+(x+\sqrt{u}y)(x-\sqrt{u}y)
(x+\sqrt{v}y)(x-\sqrt{v}y)\] 
with $uv=1$ and $u+v=\frac{1}{4}\lambda^2$. Since 
$(-\lambda)^4\ne 64$, it implies $\sqrt{v}\ne \sqrt{u}$. Lastly, 
after a linear transformation of $x$ and $y$, we can obtain that 
$\Omega_2=z^2+xy^3+x^3y+kx^2y^2$ for some $k\in \kk$. Since 
$\Omega_2$ has isolated singularity, it forces that $k\ne \pm 2$ 
as desired.

\end{proof}

Recall that $P_{\Omega-\xi}=A/(\Omega-\xi)$ for any $\xi\in \kk$. 
Let $P=P_\Omega$. Since $P=\oplus_{i\ge 0}P_i$ is graded, $P$ 
has two natural $0$-filtrations denoted by 
$\mathbb F^{Id}=\{F_i^{Id}\,|\, i\in \mathbb Z\}$ and 
$\mathbb F^{-Id}=\{F_i^{-Id}\,|\, i\in \mathbb Z\}$ respectively
such that 
\begin{align}
\label{E3.10.1}\tag{E3.10.1}
F_i^{Id}=\bigoplus_{n\ge i} P_n\ \text{and}
\ F_i^{-Id}=\bigoplus_{n\leq -i} P_n.
\end{align}
Since $\gr_{\mathbb F^{\pm Id}} P\cong P$ (with grading flipped 
$i\leftrightarrow -i$ for $\mathbb F^{-Id}$), $P$ has two 
canonical faithful $0$-valuations. We denote them by $\nu^{Id}$
and $\nu^{-Id}$ respectively. In particular, we have
\begin{align}
\label{E3.10.2}\tag{E3.10.2}
\nu^{Id}(f)=n\ \text{and}\ \nu^{-Id}(f)=-m
\end{align}
for any $f=\sum_{i=n}^{m}f_i$ with $n\le m$, $f_i\in P_i$ and 
$f_n,f_m\neq 0$.

Now consider $P_{\Omega-\xi}=A/(\Omega-\xi)$ for any 
$\xi\in \kk^\times$. We define a filtration 
$\mathbb F^{-Id}=\{F_i^{-Id}\,|\, i\in \mathbb Z\}$ of 
$P_{\Omega-\xi}$ by 
\[F_{-i}^{-Id}P_{\Omega-\xi}=\{\sum a_jf_j \,|\, 
a_j\in \kk,\ \text{$f_j$ are monomials of degree $\le i$}\}.\]
One can check that $\gr_{\mathbb F^{-Id}}P_{\Omega-\xi}\cong P$, 
and the corresponding faithful $0$-valuation is given by 
$\nu^{-Id}$ via $\nu^{-Id}(x)=-a$, $\nu^{-Id}(y)=-b$ and 
$\nu^{-Id}(z)=-c$.  We say a Poisson algebra automorphism 
$\sigma$ of $P_{\Omega-\xi}$ is {\it linear} if it preserves the 
specific $0$-filtration $\mathbb F^{-Id}$ on $P_{\Omega-\xi}$, 
or namely $\sigma(F^{-Id}_i)\subseteq F^{-Id}_i$ for $i$. 

The Poisson valuations directly apply to Poisson automorphism groups of $P_{\Omega-\xi}$ when the homogeneous potential $\Omega$ has an isolated singularity. 

\begin{lemma}
\label{zzlem3.11}
Let $A_\Omega$ be a connected graded Poisson algebra. We have the following if the 
potential $\Omega$ has an isolated singularity. 
\begin{itemize}
\item[(1)] 
Every Poisson automorphism of $P_\Omega$ is graded.
\item[(2)] 
If $\xi\neq 0$, then every Poisson automorphism of 
$P_{\Omega-\xi}$ is linear. 
\end{itemize}
\end{lemma}

\begin{proof}
We only check for $\Omega_2$ (with 
$(\deg(x),\deg(y),\deg(z))=(1,1,2)$) and the proofs for other 
cases are similar. For simplicity, we write $\Omega=\Omega_2$, 
$P=P_{\Omega_2}$ and $P_\xi=P_{\Omega_2-\xi}$.  

(1) By \Cref{zzthm3.6}(2), the fraction ring $Q=Q(P)$ is isomorphic 
to $Q(P_{\Omega_1})$ as Poisson fields. By 
\cite[Theorem 3.8]{HTWZ1}, $Q$ has exactly two faithful $0$-valuations, 
namely $\{\nu^{\pm Id}\}$ as discussed above. Let $\phi$ be any 
Poisson automorphism of $P$. We extend $\phi$ to a Poisson field 
automorphism of $Q$, which we still denote by $\phi$. It is clear 
that $\nu^{\pm Id}\circ \phi$ are two distinct faithful 
$0$-valuations of $Q$. So 
$\{\nu^{\pm Id}\}=\{\nu^{\pm Id}\circ \phi\}$.  By 
\eqref{E3.10.2}, an element $f\in P$ is homogeneous if and only 
if $\nu^{Id}(f)+\nu^{-Id}(f)=0$. Let $f$ be any homogeneous 
element of $P$. So we have 
\[0=\nu^{Id}(f)+\nu^{-Id}(f)=(\nu^{Id}\circ \phi)(f)
+(\nu^{-Id}\circ \phi)(f)=\nu^{Id}(\phi(f))+\nu^{-Id}(\phi(f)).\]
This implies that $\phi(f)$ is again homogeneous. In particular, 
$\phi(x),\phi(y),\phi(z)$ are homogeneous and nonzero.
By recycling letters $a,b,c$, respectively, we assume they have degrees $a,b,c$. Hence, 
$\phi(\Omega)=\phi(x)^4+\phi(y)^4+\phi^2(z)
+\lambda\phi(x)\phi(y)\phi(z)$ is homogeneous. So we have 
$4a=4b=2c=a+b+c$. Since $P$ is generated by $\phi(x),\phi(y),
\phi(z)$, we must have $a=b=1$ and $c=2$. Hence, $\phi$ is graded.

(2) Let $Q=Q(P_{\xi})$ for $\xi\in \kk^\times$. We can use the 
similar argument of \cite[Theorem 3.11]{HTWZ1} to show that $Q$ 
has only one faithful $0$-valuation, namely $\nu^{-Id}$ as 
discussed above. So for any Poisson automorphism $\phi$ of $P_\xi$, 
we have $\nu^{-Id}=\nu^{-Id}\circ \phi$. Let $f\in P_\xi$. By 
definition, we have $f\in F^{-Id}_i$ if and only if 
$\nu^{-Id}(f)\leq -i$ if and only if $(\nu^{-Id}\circ \phi)(f)
=\nu^{-Id}(\phi(f))\leq -i$. Hence $\phi$ preserves the 
$0$-filtration $\mathbb F^{-Id}$ of $P_\xi$ and so is linear.
\end{proof}

Next, we explicitly compute the Poisson automorphism groups of 
$P_{\Omega-\xi}$. 

\begin{lemma}
\label{zzlem3.12}
Let $\Omega=x^3+y^3+z^3+\lambda xyz$ for $(-\lambda)^3\neq 27$ 
with $\deg(x)=1, \deg(y)=1,\deg(z)=1$.
\begin{itemize}
\item[(1)] 
There is a short exact sequence of groups:
\[1\to G\to {\rm Aut}_P(P_\Omega)\to C_3\to 1,\]
where $G=\{(\alpha_1,\alpha_2,\alpha_3)\,|\, \alpha_1^3
=\alpha_2^3=\alpha_3^3=\alpha_1\alpha_2\alpha_3\}\subset 
\kk^\times\times \kk^\times\times \kk^\times$ and 
${\rm Aut}_P(P_\Omega)\cong C_3\ltimes G$.
\item[(2)] 
There is a short exact sequence of groups:
\[
1\to G'\to {\rm Aut}_P(P_{\Omega-1})\to C_3\to 1,\]
where $G'=\{(\alpha_1,\alpha_2,\alpha_3)\,|\, \alpha_1^3
=\alpha_2^3=\alpha_3^3=\alpha_1\alpha_2\alpha_3=1\}\subset 
\kk^\times\times \kk^\times\times \kk^\times$ and 
${\rm Aut}_P(P_{\Omega-1})\cong C_3\ltimes G'$. 
\end{itemize}
\end{lemma}

\begin{proof}
(1) Note that the argument of \cite[Theorem 1]{MTU} works for 
${\rm Aut}_{grP}(P_\Omega)$ as well. So 
${\rm Aut}_{grP}(P_\Omega)$ is generated by all possible 
diagonal actions 
$G=\{(\alpha_1,\alpha_2,\alpha_3)\,|\,\alpha_1^3=\alpha_2^3
=\alpha_3^3=\alpha_1\alpha_2\alpha_3\}$ and the permutation 
$\tau(x,y,z)=(y,z,x)$. So our result follows \Cref{zzlem3.11}. 

(2) By \Cref{zzlem3.11}, every automorphism $\phi$ of 
$P_{\Omega-1}$ is linear. Note that for $0\le i\le 2$, we will 
identify  $A_i$ with $(P_{\Omega-1})_i$ as vector spaces. So we 
can write 
\begin{align*}
\phi(x)=f_1+f_0,\ \phi(y)=g_1+g_0,\ \phi(z)=h_1+h_0
\end{align*}
where $f_i,g_i,h_i$ are homogeneous polynomials of degree $i$ 
in $A_i=(P_{\Omega-1})_i$ for all possible $0\le i\le 1$. For simplicity, we will denote the images of 
$x, y, z$ in $P_{\Omega-1}$ by $x,  y, z$. We have 
\[\{x,y\}=3z^2+\lambda xy,\ 
\{y,z\}=3x^2+\lambda yz,\ \{z,x\}=3y^2+\lambda xz
\]
in $P_{\Omega-1}$. We apply $\phi$ to each of the above three 
equations. After comparing the constant terms on both sides of 
the resulting equations, we obtain
\[
3f_0^2+\lambda g_0h_0= 3g_0^2+\lambda f_0h_0
=3h_0^2+\lambda f_0g_0=0.\]
Note that we have $(a+b+c)\Omega=ax\Omega_x+by\Omega_y+cz\Omega_z$ in $A=\Bbbk[x, y, z]$. Since $(a, b, c)=(1, 1, 1)$, it follows that 
\begin{eqnarray*}
\Omega(f_0, g_0, h_0) &= &\frac{1}{a+b+c}[af_{0} \Omega_{x}(f_0, g_0, h_0)+bg_{0}\Omega_y(f_0, g_0, h_0)+ch_0\Omega_z(f_0, g_0, h_0)]\\
&=&\frac{1}{3}[f_{0} (3g_0^2+\lambda f_0h_0) + g_0(3h_0^2+\lambda f_0g_0)+h_0(3f_0^2+\lambda g_0h_0)]\\
&=&0.
\end{eqnarray*} 
So $(f_0, g_0,h_0)$ is a singular point on the surface $\Omega=0$. 
Hence $(f_0,g_0,h_0)=(0, 0, 0)$ since $\Omega$ only has an 
isolated singularity at the origin. So $\phi$ maps 
$\kk x+\kk y+\kk z$ to itself. This implies that 
$\gr(\phi): \gr (P_{\Omega-1})\to \gr (P_{\Omega-1})$ is a 
graded automorphism of $P_{\Omega}(\cong \gr (P_{\Omega-1}))$ 
which has been described in part (1).  The rest of the proof 
follows from a direct computation.  
\end{proof}

\begin{lemma}
\label{zzlem3.13}
Let $\Omega=x^4+y^4+z^2+\lambda xyz$ for $(-\lambda)^4\neq 64$ 
with $\deg(x)=\deg(y)=1$ and $\deg(z)=2$. 
\begin{itemize}
\item[(1)] 
There is a short exact sequence of groups:
\[
1\to G\to {\rm Aut}_P(P_\Omega)\to C_2\to 1
\]
where $G=\{(\alpha_1,\alpha_2)\,|\, \alpha_1^2=\alpha_2^2\}
\subset \kk^\times\times \kk^\times$ and 
${\rm Aut}_P(P_\Omega)\cong C_2\ltimes G$.
\item[(2)] 
There is a short exact sequence of groups:
\[
1\to G'\to {\rm Aut}_P(P_{\Omega-1})\to C_2\to 1
\]
where $G'=\{(\alpha_1,\alpha_2)\,|\, \alpha_1^4=\alpha^4_2=1, 
\alpha_1^2=\alpha_2^2\}\subset \kk^\times\times \kk^\times$ 
and ${\rm Aut}_P(P_{\Omega-1})\cong C_2\ltimes G'$. 
\end{itemize}
\end{lemma}

\begin{proof}
(1) By \Cref{zzlem3.11}, every Poisson automorphism $\phi$ of 
$P_\Omega$ is graded. For the convenience of this proof, we write 
$(x,y,z)=(x_1,x_2,x_3)$. In general, we can assume that $\phi$ 
is given by
\begin{equation*}
\left\{
\begin{aligned}
\phi(x_1)&=\alpha_1x_1+\alpha_2x_2\\
\phi(x_2)&=\beta_1x_1+\beta_2x_2\\
\phi(x_3)&=\gamma x_3+h(x_1, x_2)
\end{aligned}\right.
\end{equation*}
where $\gamma$ and $\alpha_1\beta_2-\alpha_2\beta_1$ are not 
zero in $\kk$ and $h(x_1,x_2)$ is a quadratic polynomial in 
$\kk[x_1,x_2]$, as 
$\det(\phi)=\gamma(\alpha_1\beta_2-\alpha_2\beta_1)\ne 0$. After 
a proper re-scaling of variables, we first assume that 
$\alpha_1\beta_2-\alpha_2\beta_1=1$. Applying $\phi$ to the 
equation $\{x_1,x_2\}=2x_3+\lambda x_1x_2$, we get 
$\gamma=1$ and $h(x_1, x_2)=\frac{\lambda}{2}
(x_1x_2-\phi(x_1)\phi(x_2))$. Applying $\phi$ to the other two 
equations $\{x_2,x_3\}=4x_1^3+\lambda x_2x_3$ and 
$\{x_3,x_1\}=4x_2^3+\lambda x_1x_3$, we further obtain the 
following relations:
\begin{equation*}
\left\{\begin{aligned}
24\alpha_1^2\alpha_2
&=\lambda^2(\beta_1+\alpha_2\beta_1^2+2\beta_1\beta_2\alpha_1)
 =3\lambda^2\alpha_1\beta_1\beta_2,\\
24\alpha_1\alpha_2^2
&=\lambda^2(-\beta_2+\alpha_1\beta_2^2+2\beta_1\beta_2\alpha_2)
 =3\lambda^2\alpha_2\beta_1\beta_2,\\
24\beta_1^2\beta_2
&=\lambda^2(-\alpha_1+\alpha_1^2\beta_2+2\alpha_1\alpha_2\beta_1)
 =3\lambda^2\alpha_1\alpha_2\beta_1,\\
24\beta_1\beta_2^2
&=\lambda^2(\alpha_2+\alpha_2^2\beta_1+2\alpha_1\alpha_2\beta_2)
 =3\lambda^2\alpha_1\alpha_2\beta_2.
\end{aligned}\right.
\end{equation*}
Since $\alpha_1\beta_2-\alpha_2\beta_1=1$,  $\alpha_1$ and 
$\alpha_2$ cannot be simultaneously equal to zero and similar 
situation happens to $\beta_1$ and $\beta_2$. After simplifying 
the above equations, we have 
$8\alpha_1\alpha_2=\lambda^2\beta_1\beta_2$ and 
$8\beta_1\beta_2=\lambda^2\alpha_1\alpha_2$ and so 
$(64-\lambda^4)\alpha_1\alpha_2=(64-\lambda^4)\beta_1\beta_2=0$. 
Since $\lambda^4\neq 64$, we get 
$\alpha_1\alpha_2=\beta_1\beta_2=0$.  Thus, we have either 
$\alpha_2=\beta_1=0$ or $\alpha_1=\beta_2=0$. Now, by taking care 
of the re-scaling of variables, we can find a permutation 
$\sigma\in S_2$ and write $\phi$ as 
\begin{align}
\label{E3.13.1}\tag{E3.13.1}
\phi(x_1)=\alpha_1x_{\sigma(1)},\,\, 
\phi(x_2)=\alpha_2x_{\sigma(2)},\,\,
\phi(x_3)=\alpha_3x_3+h.
\end{align}
for some nonzero scalars $\alpha_1,\alpha_2,\alpha_3$ and some 
quadratic polynomial $h(x_1,x_2)$. Then, it is routine to check 
that $\phi$ is a Poisson automorphism of $P_\Omega$ if and only if
\[
\alpha_3={\rm sgn}(\sigma)\alpha_1\alpha_2,\ 
\alpha_1^2=\alpha_2^2,\  
h=\frac{\lambda}{2}({\rm sgn}(\sigma)-1)\alpha_1\alpha_2x_1x_2.
\]
Consider the normal subgroup $G$ of ${\rm Aut}_{grP}(P_\Omega)$ 
given by $\phi(x)=\alpha_1x, \phi(y)=\alpha_2y, 
\phi(z)=\alpha_1\alpha_2z$ satisfying $\alpha_1^2=\alpha_2^2$. 
Finally, it is clear that ${\rm Aut}_{grP}(P_\Omega)/G\cong S_2$.

(2) By \Cref{zzlem3.11}, every automorphism $\phi$ of 
$P_{\Omega-1}$ is linear. So we can write 
\begin{align*}
\phi(x)=f_1+f_0,\,\,
\phi(y)=g_1+g_0,\,\,
\phi(z)=h_2+h_1+h_0
\end{align*}
where $f_i,g_i,h_i$ are homogeneous polynomials of degree $i$ 
in $A_i=(P_{\Omega-1})_i$ for all possible $0\le i\le 2$. We 
apply $\phi$ to the following Poisson brackets:
\[\{x,y\}=2z+\lambda xy,\ 
\{y,z\}=4x^3+\lambda yz,\ 
\{z,x\}=4y^3+\lambda xz.\]
By examining these equations at the degree zero part, we get 
\[2h_0+\lambda f_0g_0=4f_0^3+\lambda g_0h_0
=4g_0^2+\lambda f_0h_0=0.\]
This implies that $(f_0,g_0,h_0)$ is a singular point for 
$\Omega=0$. So $f_0=g_0=h_0=0$ since $\Omega$ only has an isolated 
singularity at the origin. Now since $\phi(\Omega)=1$ in $P_{\Omega-1}$, 
we get
\[f_1^4+g_1^4+(h_2+h_1)^2+\lambda f_1g_1(h_2+h_1)=1.\]
The degree two part of the above equation yields that $h_1^2=0$ 
and thus $h_1=0$. Hence, $\phi$ maps $\kk x+\kk y$ (resp.
$\kk z$) to itself, and we can write 
$\phi$ as in \eqref{E3.13.1}. Similar to (1), indeed, we have 
\begin{align*}
\phi(x_1)=\alpha_1x_{\sigma(1)},\,\,
\phi(x_2)=\alpha_2x_{\sigma(2)},\,\,
\phi(x_3)=\alpha_3x_3+\frac{\lambda}{2}
({\rm sgn}(\sigma)-1)\alpha_1\alpha_2x_1x_2,
\end{align*}
for some $\sigma\in S_2$ and $\alpha_1^4=\alpha_2^4=1$ and 
$\alpha_1^2=\alpha_2^2$ since $\phi(\Omega)=1$. So the result 
follows.
\end{proof}

\begin{lemma}
\label{zzlem3.14}
Let $\Omega=x^6+y^3+z^2+\lambda xyz$ for $(-\lambda)^6\neq 432$ 
with $\deg(x)=1, \deg(y)=2,\deg(z)=3$.
\begin{itemize}
\item[(1)] 
Every Poisson automorphism of $P_\Omega$ is of the form 
\begin{equation*}
x\mapsto \zeta x,\ y\mapsto \zeta^2y,\ z\mapsto \zeta^3z
\end{equation*}
for some $\zeta\in \kk^\times$.
\item[(2)] 
Every Poisson automorphism of $P_{\Omega-1}$ is of the form 
\begin{equation*}
x\mapsto \zeta x,\ y\mapsto \zeta^2y,\ z\mapsto \zeta^3z,
\end{equation*}
where $\zeta^6=1$.
\end{itemize}
\end{lemma}

\begin{proof}
(1) By \Cref{zzlem3.11}, every automorphism $\phi$ of 
$P_\Omega$ is graded. So we can write 
\begin{equation}\label{E3.14.1}\tag{E3.14.1}
\phi(x)=\alpha_1x,\ \phi(y)=\beta_1y+\beta_0x^2,\ 
\phi(z)=\gamma_2z+\gamma_1xy+\gamma_0x^3
\end{equation}
for some $\alpha_1,\beta_1,\gamma_2\in \kk^\times$ and 
$\beta_0,\gamma_1,\gamma_0\in \kk$. Applying $\phi$ to
$\{z,x\}=3y^2+\lambda xz$, we obtain that $\beta_0=\gamma_0
=\gamma_1=0$. So $\phi$ equals a scalar multiple when it is 
applied to the generators. Finally, from 
$\{\phi(y),\phi(z)\}=\phi(6x^5+\lambda yz)$, 
$\{ \phi(x), \phi(y)\}=\phi(2z+\lambda xy)$ and 
$\{\phi(z), \phi(x)\}=\phi(3y^2+\lambda xz)$, we obtain
$\beta_1\gamma_2=\alpha_1^5$, $\alpha_1\beta_1=\gamma_2$, 
and $\gamma_2\alpha_1=\beta_1^2$. Let $\zeta=\gamma_2/\beta_1$. 
We have $\alpha_1=\zeta,\beta_1=\zeta^2,\gamma_2=\zeta^3$. So 
the result follows. 

(2) By \Cref{zzlem3.11}, every automorphism $\phi$ of 
$P_{\Omega-1}$ is linear. Note that $A_i=(P_{\Omega-1})_i$ 
for $0\le i\le 5$. So we can write 
\begin{align*}
\phi(x)=f_1+f_0,\,\,
\phi(y)=g_2+g_1+g_0,\,\,
\phi(z)=h_3+h_2+h_1+h_0
\end{align*}
where $f_i,g_i,h_i$ are homogeneous polynomials of degree $i$ 
in $A_i=(P_{\Omega-1})_i$ for all possible $0\le i\le 3$. We 
apply $\phi$ to the following Poisson brackets:
\[\{x,y\}=2z+\lambda xy,\ \{y,z\}=6x^5+\lambda yz,\ \{z,x\}
=3y^2+\lambda xz\]
and examine the resulting equations at different degrees. First 
of all, the degree zero part yields that 
\[2h_0+\lambda f_0g_0=6f_0^5+\lambda g_0h_0
=3g_0^2+\lambda f_0h_0=0.\]
This is equivalent to the claim that $(f_0,g_0,h_0)$ is a 
singular point for $\Omega=0$. So $f_0=g_0=h_0=0$ since 
$\Omega$ only has an isolated singularity at the origin. Then, the degree 
one part further implies that $h_1=0$. Now since $\phi(\Omega)=1$ 
in $P_{\Omega-1}$, we get 
\[f_1^6+(g_2+g_1)^3+(h_3+h_2)^2+\lambda f_1(g_2+g_1)(h_3+h_2)=1.\]
The degree three part of the above equation yields that $g_1^3=0$ 
and thus $g_1=0$. Furthermore, the degree four part implies that 
$h_2^2=0$ and $h_2=0$. Hence, we have $\phi$ as in 
\eqref{E3.14.1}. Similar to the proof of part (1), indeed, we 
have $\phi(x)=\zeta x,\phi(y)=\zeta^2 y,\phi(z)=\zeta^3 z$ for 
some $\zeta\in \kk^\times$. Finally, $\phi(\Omega)=1$ implies 
that $\zeta^6=1$. 
\end{proof}

\begin{proof}[Proof of \Cref{zzthm3.9}] 
It suffices to prove the result when $\Omega$ is case (2) and  case (3), 
since case (1) was shown in \cite{MTU}. Though our argument 
works for all cases; here, we only provide the details for Case 
(2) as an illustration.  

Let $\phi$ be a Poisson algebra automorphism of $A_\Omega$. By 
the same argument as in \Cref{zzlem3.3}, the restriction of 
$\phi$ on the Poisson center $\kk[\Omega]$ is given by 
$\phi(\Omega)=\alpha \Omega$ with $\alpha\in \kk^\times$. So 
$\phi$ preserves the principal ideal generated by $\Omega$. Let 
$\phi'$ denote the induced Poisson algebra automorphism for 
$P_\Omega$. By \Cref{zzlem3.11}, $\phi'$ is graded. Moreover, 
since the equations $\{x,y\}=\Omega_z,\{y,z\}=\Omega_x,
\{z,x\}=\Omega_y$ are homogeneous of degree $<\deg(\Omega)$, 
we can lift $\phi'$ to a unique graded Poisson automorphism of 
$A_\Omega$, denoted by $\sigma$. It is clear that $\sigma'=\phi'$.
Let $\varphi=\phi\circ \sigma^{-1}$. Then it satisfies the 
equation $\varphi'={\rm id}_{P_\Omega}$. It remains to show 
that $\varphi={\rm id}_{A_\Omega}$. Note that we can write 
\begin{align*}
\varphi(x)=x+\Omega f,\,\,
\varphi(y)=y+\Omega g,\,\,
\varphi(z)=z+\Omega h
\end{align*}
for some polynomials $f,g,h\in A_\Omega$. An easy computation 
yields that $\varphi(\Omega)=\Omega+\Omega\,\alpha(f,g,h)$, 
where $\alpha(f,g,h)\in (A_\Omega)_{\ge 1}$. Since $\varphi$ 
induces an algebra automorphism of the Poisson center 
$\kk[\Omega]$ of $A_\Omega$, we must have $\alpha=0$ and 
$\varphi(\Omega)=\Omega$. 

Now we consider a $\kk$-linear basis $\mathbb B:
=\{1,x,y,z\}\cup \{b_s\}$ of the algebra $P_\Omega$ consisting 
of all possible monomials 
$\{x^{s_1}y^{s_2}z^{s_3}\,|\, s_1,s_2\ge 0, 0\le s_3\le 1\}$. 
We also treat $\mathbb B$ as a fixed subset of monomials in 
$A_\Omega$ by lifting. We claim that every polynomial $f$ in 
$A_\Omega$ can be written in the form  
\begin{align}
\label{E3.14.2}\tag{E3.14.2}
f=1f^1(\Omega)+xf^x(\Omega)+yf^y(\Omega)+zf^z(\Omega)
+\sum_{b_s} b_sf^{b_s}(\Omega)
\end{align}
where each $f^*(\Omega)\in \kk[\Omega]$. We prove this claim by 
induction on $\deg(f)$. It is trivial for $\deg(f)=0$. Suppose 
our claim holds for $\deg(f)\leq m$. For any polynomial $f$ of 
degree $m+1$, we can write 
\[f=1f^1+xf^x+yf^y+zf^z+\sum_{b_s}b_sf^{b_s}+g\Omega \]
for some scalars $f^*\in \kk$ and $\deg(g)=m-3$  by looking at 
the image of $f$ in $P_\Omega$. So, by the induction hypothesis, we 
can write $g$ in the required form. We get our claim by replacing 
$g$ above. Therefore, we can write
\begin{align}
\label{E3.14.3}\tag{E3.14.3}
\varphi(x)=1f^1(\Omega)+xf^x(\Omega)+yf^y(\Omega)+zf^z(\Omega)
+\sum_{b_s}b_sf^{b_s}(\Omega)
\end{align}
for some $f^*(\Omega)\in \kk[\Omega]$.

Now for each scalar $\xi\neq 0$, let $\pi_\xi: A_\Omega\to 
P_{\Omega-\xi}$ be the quotient map and write $\varphi_\xi'$ 
as the induced automorphism of $\varphi$ since 
$\varphi(\Omega)=\Omega$. Note that the image of $\mathbb B$ 
via $\pi_\xi$ is a $\kk$-basis of $P_{\Omega-\xi}$, which we 
continue to denote by $b_s$ etc. So we have
\begin{align}
\label{E3.14.4}\tag{E3.14.4}
\varphi_\xi'(x)=1f^1(\xi)+xf^x(\xi)+yf^y(\xi)+zf^z(\xi)
+\sum_{b_s}b_sf^{b_s}(\xi).
\end{align}
By \Cref{zzlem3.11}, $\varphi_\xi'$ is linear. Thus, 
$f^{b_s}(\xi)=0$ for all $\xi\neq 0$ and $b_{s}$. Hence 
$f^{b_s}(\Omega)=0$ and \eqref{E3.14.3} reduces to 
\begin{align*}
\varphi(x)=1f^1(\Omega)+xf^x(\Omega)+yf^y(\Omega)+zf^z(\Omega).
\end{align*}
Moreover, since $\varphi'(x)=x$ in $P_{\Omega}$, we have 
$f^x(\Omega)=1+p\Omega \neq 0$ for some $p\in \kk[\Omega]$. 
If $f^y(\Omega)\neq 0$, we can choose some nonzero $\xi_0$ 
such that $f^x(\xi_0),f^y(\xi_0)\neq 0$. But in this case, 
$\varphi_{\xi_0}'(x)$ in \eqref{E3.14.4} contains both 
terms of $x$ and $y$. This contradicts to the description 
of ${\rm Aut}_{grP}(P_{\Omega-\xi_0})$ by \Cref{zzlem3.13}. 
So $f^y(\Omega)=0$ and we get $f^1(\Omega)=f^z(\Omega)=0$ in 
the same fashion. This implies that $\varphi(x)=x(1+p\Omega\,)$. 
Similarly, we get $\varphi^{-1}(x)=x(1+\Omega\,h)$ for some 
$h\in \kk[\Omega]$. By using $\varphi(\Omega)=\Omega$, we obtain 
that
\[x=\varphi(\varphi^{-1}(x))=\varphi(x(1+\Omega\,h))
=x(1+\Omega\, p)(1+\Omega\,h )=x,\]
which implies that $p=h=0$. Therefore, $\varphi(x)=x$. We further get $\varphi(y)=y$ and $\varphi(z)=z$ by the same argument. 
Hence, $\varphi$ is the identity as required. 
\end{proof}

\begin{example}
\label{zzexa3.15}
Let $\Omega=z^2+x^3y$ in \Cref{tab:112} with $\deg x=\deg y=1$ 
and $\deg z=2$. Then the Poisson structure of $A_{\Omega}$ is 
determined by
$$\{x,y\}=2z, \quad
\{y,z\}=3x^2 y, \quad
\{z,x\}=x^3.$$
Consider the algebra automorphism of the polynomial ring 
$\Bbbk[x,y,z]$ defined by
\[\phi(x)=x,\,\, \phi(y)=y-x^3-2z, \,\,\,\textnormal{and}
\,\,\, \phi(z)=z+x^3.\] 
It is straightforward to check that $\phi$ is an ungraded 
Poisson automorphism of $A_\Omega$. 
\end{example}

We can prove that if $\Omega$ belongs to the Weyl type, meaning $Q(P_\Omega)\cong K_{Weyl}$, then $A_\Omega$ possesses ungraded Poisson automorphisms. However, we have not discovered any ungraded automorphisms for other irreducible $\Omega$, i.e., those satisfying $Q(P_\Omega)\cong K_q$. Therefore, we pose the following question.

\begin{question}
\label{zzque3.16}
Let $\Omega$ be a homogeneous polynomial of degree $|x|+|y|+|z|$. If $Q(P_\Omega)\cong K_q$ as Poisson fields, is every Poisson automorphism of $A_\Omega$ and $P_\Omega$ graded? 
\end{question}

\section{Rigidities}
\label{zzsec4}
In this section, we discuss several different rigidities about
the Poisson algebras related to $\Omega$. We intend to provide general methods for these rigidities and omit some details. Unless mentioned otherwise, we will retain \Cref{zzhyp0.1} for this whole section.

\subsection{Rigidity of graded twistings}
\label{zzsec4.1}
Note that the {\it rigidity of graded twistings}, denoted by 
$rgt$, was defined in \Cref{zzdef1.2}. The following lemma 
provides an easy way to compute Poisson derivations and $rgt$ 
of $A_\Omega$.

\begin{lemma}
\label{zzlem4.1}
Let $A$ be a connected graded unimodular Poisson algebra $A_\Omega$ as given in \Cref{zzhyp0.1}. For any derivation 
$\delta$ of $A$, we have the following:
\begin{enumerate}
\item[(1)] 
$\delta$ is a Poisson derivation if and only if 
$\divv(\delta)\pi_{\Omega}=\pi_{\delta(\Omega)}$ in $\mathfrak X^2(A)$.
\item[(2)] 
Suppose $\divv(\delta)=0$. Then $\delta$ is a Poisson derivation 
if and only if $\delta(\Omega)=0$.
\item[(3)]
We have 
\begin{align*}
rgt(A)&=1-\dim_\kk Gspd(A)\\
&=1-\dim_\kk Gpd(A)\\
&=1-\dim_\kk (PH^1(A))_0\\
&=-\dim_\kk\{\delta\in Gspd(A)\,|\, \divv(\delta)=0\}\\
&=-\dim_\kk\{\delta\in Gpd(A)\,|\, \divv(\delta)=0\}\\
&=-\dim_\kk\{\delta\in (Der(A))_0\,|\, \divv(\delta)=\delta(\Omega)=0\}.
\end{align*}
\end{enumerate}
\end{lemma}

\begin{proof}
Note that $\delta$ is a Poisson derivation if and only if 
$d_\pi^1(\delta)=0$. Recall that the cochain complexes 
\eqref{E1.3.1} and \eqref{E1.3.3} are isomorphic to each other. 
So $(d_{\pi}^1(\delta)(y,z), d_{\pi}^1(\delta)(z,x), 
d_{\pi}^1(\delta)(x,y))= \delta_\Omega^1(\delta(x), \delta(y), 
\delta(z))$ as vectors.
By \eqref{E1.3.5} and  \eqref{E3.1.1}, we conclude that 
\begin{equation}
\label{E4.1.1}\tag{E4.1.1}
d_\pi^1(\delta)=~\divv(\delta)\,\pi_{\Omega}-\pi_{\delta(\Omega)}
\end{equation}
in $\mathfrak X^2(A)$. Thus (1) follows immediately. 

For (2), suppose $\delta(\Omega)=0$. Then by \eqref{E4.1.1}, 
$\delta$ is a Poisson derivation. Conversely, suppose $\delta$ 
is a Poisson derivation. By \eqref{E4.1.1}, 
$\delta(\Omega)\in \kk=A_0$. We write $\delta=\sum \delta_i$, 
where each $\delta_i$ is a homogeneous derivation of degree $i$.  
Since $\deg(\Omega)=n$, we get $\delta_i(\Omega)=0$ for 
$i\neq -n$. Moreover, $\delta_{-n}(x)$, if not zero, has 
degree $a-n<0$. So $\delta_{-n}(x)=0$ as $A$ is 
$\mathbb N$-graded. Similarly, we get 
$\delta_{-n}(y)=\delta_{-n}(z)=0$. Therefore, 
$\delta_{-n}(\Omega)=0$ and so $\delta(\Omega)=0$, as desired.

For (3), consider the subspace 
$Gpd_0(A)=\{\delta\in Gpd(A)\,|\, \divv(\delta)=0\}$ of $Gpd(A)$. 
Let $\delta\in Gpd(A)$. By Lemma \cite[Lemma 1.2 (2)]{TWZ}, 
$\divv(\delta)\in \kk$. In particular, $\divv(E)=a+b+c\in 
\kk^\times$ for the Euler derivation $E$ of $A$. Hence 
\begin{equation}
\label{E4.1.2}\tag{E4.1.2}
\delta'=\delta-\frac{\divv(\delta)}{\divv(E)}\,E\in Gpd_0(A).  
\end{equation}
So $Gpd(A)=Gpd_0(A)\oplus \kk\, E$. Thus, the formulas of 
$rgt(A)$ can be derived from \cite[Lemma 4.1]{TWZ}, where 
the last equality follows from (2).
\end{proof}

When $A$ is generated in degree one, $rgt(A_{\Omega})$ is computed for each $\Omega$ from \Cref{tab:111} in \cite[Corollary 6.7]{TWZ}; and it follows that any graded unimodular structure $\pi_\Omega$ 
on $A$ is rigid (namely, $rgt(A_\Omega)=0$) if and only if $\Omega$ 
is irreducible. Moreover, one can easily obtain that ${\rm GK}(A_{sing})=0$ for irreducible $\Omega$ with isolated singularity, ${\rm GK}(A_{sing})=1$ for irreducible $\Omega$ without isolated singularity, and ${\rm GK}(A_{sing})$ is $1$ or $2$ for reducible $\Omega$. Now, we can generalize it to any weighted case.  

\begin{theorem}
\label{zzthm4.2}
Let $A_\Omega$ be a connected graded unimodular Poisson algebra given in \Cref{zzhyp0.1}. 
Then $A_\Omega$ is rigid if and only if the potential $\Omega$ 
is irreducible. In \Cref{tab:rgtanddim} in the appendix, we will 
list all $rgt(A)$ and $GKdim(A_{sing})$ of these $\Omega$ from 
\Cref{tab:112} to \Cref{tab:abc2}. 
\end{theorem}

\begin{proof}
We apply \Cref{zzthm3.4} (also see \Cref{zzlem1.3}) and 
\Cref{zzlem2.2}. Indeed, our method is a case-by-case 
verification. First of all, when $a=b<c$ and $\Omega=h(x,y)$ 
is reducible, we have $rgt(A)\neq 0$ according to \Cref{zzlem2.2}. 
For the rest of the proof, we only provide the details for the 
following two cases for illustration: irreducible 
$\Omega_1=z^2+xy^3+\lambda x^2y^2+x^3y$ and reducible 
$\Omega_2=x^2z+xy^3$ in \Cref{tab:112} (with $\deg(x)=\deg(y)=1$ 
and $\deg(z)=2$). Let $\phi$ be any graded Poisson derivation 
of degree zero. Replacing $\phi$ by $\phi+cE$ for some suitable 
scalar $c\in \kk$ if needed, we can always assume 
${\rm div}(\phi)=0$ by \eqref{E4.1.2}. Thus, we can write 
$\phi(x)=\alpha_1 x+\alpha_2y$, $\phi(y)=\alpha_3x+\alpha_4y$ and 
$\phi(z)=(-\alpha_1-\alpha_4)z+\alpha_5x^2+\alpha_6xy+\alpha_7y^2$ 
for some $\alpha_i\in \kk$. By \Cref{zzlem4.1} (2), we have 
$\phi(\Omega)=0$ for $\Omega=\Omega_1$ or $\Omega_2$. From 
$\phi(\Omega_1)=0$, one can show that $\alpha_i=0$ for 
$i=1, \cdots, 7$, which implies $\phi=0$. As a result, 
$rgt(A_{\Omega_1})=0$. From $\phi(\Omega_2)=0$, we obtain 
$\alpha_7+3\alpha_3=0$ and $\alpha_i=0$ for any 
$i\neq 3, 7$. This implies that 
$Gpd(A_{\Omega_2})={\rm span}\{E,\phi\}$, where $\phi(x)=0$, 
$\phi(y)=x$, $\phi(z)=-3y^2$. Hence we get 
$rgt(A_{\Omega_2})=-1$. Finally, a standard Gr\"obner basis 
computation yields all possible $\GKdim(A_{sing})$. We skip 
the details here. 
\end{proof}

\subsection{Rigidity of gradings}
\label{zzsec4.2}
In this and the following subsections, we use Poisson valuations to establish results about the rigidity of gradings and filtrations. We believe that these kinds of rigidities deserve more attention. 

\begin{theorem}
\label{zzthm4.3}
Let $A_\Omega$ be a connected graded Poisson algebra with 
$\Omega$ having an isolated singularity. Then
$P_\Omega$ has a unique connected grading such that it is 
Poisson graded with nonzero degree one part. 
\end{theorem}

\begin{proof}
We only prove the result for $\Omega:=\Omega_1$ in \Cref{zzlem3.10}. 
The argument for $\Omega_2$ and $\Omega_3$ is analogous.
By \cite[Theorem 3.8]{HTWZ1}, $Q(P_\Omega)$ only has two 
faithful $0$-valuations, which are denoted by $\{\nu^{\pm Id}\}$ as 
before. They correspond to connected gradings on $P_\Omega$ with 
$x,y,z$ in degree 1 according to \eqref{E3.10.1}. We assume 
$\mathbb G^{\pm Id}$ to be the $0$-filtrations on $Q(P_\Omega)$ 
associated with any new grading, and denote by $\mu^{\pm Id}$ the 
corresponding faithful $0$-valuations. It is clear that 
$\mu^{-Id}(f)<0$ for some $f\in P_\Omega$. So we have 
$\mu^{-Id}=\nu^{-Id}$ and $\mu(x)=-1,\mu(y)=-1,\mu(z)=-1$. 
So we can write $x=x_0+f$, $y=y_0+g$ , and $z=z_0+h$, where 
$x_0,y_0,z_0$ are homogeneous of new degree $1$ and $f,g,h\in \kk$. 
Since every non-trivial linear combination of $x,y,z$ has $\mu^{-Id}$-value 
$-1$ (as $\mu^{-Id}=\nu^{-Id}$), $x_0,y_0,z_0$ are linearly 
independent. Since $\Omega=x^3+y^3+z^3+\lambda xyz=0$ is a sum 
of homogeneous relations, we get
\[
(3f^2+\lambda gh)x_0+(3g^2+\lambda fh)y_0+(3h^2+\lambda fg)z_0=0.
\]
This implies that 
$3f^2+\lambda gh=3g^2+\lambda fh=3h^2+\lambda fg=0$, or 
equivalently $(f,g,h)$ is a singular point of $\Omega=0$. 
Since $\Omega$ has an isolated singularity at the origin, we get $f=g=h=0$ 
and $x=x_0$, $y=y_0$, $z=z_0$ are homogeneous of degree $1$ 
in this new grading. Since $P_\Omega$ is generated by $x,y,z$, 
the new grading agrees with the given grading. 
\end{proof}

\subsection{Rigidity of filtrations}
\label{zzsec4.3}

\begin{theorem}
\label{zzthm4.4}
Assume that $\Omega$ has an 
isolated singularity. Then $P_{\Omega-\xi}$, with 
$\xi\neq 0$, has a unique filtration $\mathbb F$ such 
that the associated graded ring 
$\gr_\mathbb F(P_{\Omega-\xi})$ is a connected graded Poisson 
domain with nonzero degree one part. 
\end{theorem}

\begin{proof}
Again, we only prove the result for $\Omega:=\Omega_1$ in 
\Cref{zzlem3.10}. The argument for $\Omega_2$ and $\Omega_3$ 
is analogous. The result follows from 
\cite[Theorem 3.11]{HTWZ1} that $Q(P_{\Omega-\xi})$ has 
only one faithful $0$-valuation $\nu^{-Id}$ and \Cref{zzpro1.6}.
\end{proof}

By \Cref{zzthm4.2}, we have the following
\begin{equation}
\label{E4.4.1}\tag{E4.4.1}
{\text{$\Omega$ being irreducible}}
\Leftrightarrow
{\text{$rgt(A_{\Omega})=0$.}}
\end{equation}

Now \Cref{zzthm4.3} and \Cref{zzthm4.4} can be summarized as
\begin{equation}
\label{E4.4.2}\tag{E4.4.2}
{\text{$P_{\Omega}$ having a unique grading}}
\Leftarrow 
{\text{$\Omega$ having isolated singularity}}
\Rightarrow
{\text{$P_{\Omega-1}$ having a unique filtration.}}
\end{equation}

There is another diagram for balanced irreducible
potentials $\Omega$, see \eqref{E5.11.1}.

\section{$K_1$-sealedness and $uPH^2$-vacancy}
\label{zzsec5}
In this section, we introduce two technical concepts --
$K_1$-sealedness [\Cref{zzdef5.2}(2)] and $uPH^2$-vacancy
[\Cref{zzdef5.8}(3)]. Together with $H$-ozoneness [\Cref{zzdef1.1}(3)], they will play an important role in computing Poisson cohomology in the next section. Note that 
the $uPH^2$-vacancy of $A_{\Omega}$ is independent of the
choices of the graded generators $(x,y,z)$; however, the
$K_1$-sealedness of $\Omega$ may be dependent on 
the choices of the graded generators $(x,y,z)$. In this section, we assume that  $n$ is not necessarily equal to $a+b+c$.

\subsection{$K_1$-sealedness}
\label{zzsec5.1}
Let $\Omega$ be a homogeneous element of degree $n>0$ in the 
weighted polynomial ring $A:=\Bbbk[x,y,z]$ with $({\rm deg}(x), {\rm deg}(y), {\rm deg}(z))=(a, b, c)\in \mathbb{N}_{+}^{3}$. Recall that the 
{\it Koszul complex} $K_\bullet(\overrightarrow{\nabla}\Omega)$ 
given by the sequence $\overrightarrow{\nabla}\Omega:
=(\Omega_x,\Omega_y,\Omega_z)$ in $A$ is:
\begin{equation}
\label{E5.0.1}\tag{E5.0.1}
\begin{tabular}{ccccc}
&\,$A[-2n+b+c]$ & &$A[a-n]$& \vspace*{-2mm}\\
$0\to A[-3n+(a+b+c)]\xrightarrow{\overrightarrow{\nabla}\Omega}$
& \hspace*{-2.5mm}$\oplus A[-2n+a+c]$
&\hspace*{-2.5mm}$\xrightarrow{\overrightarrow{\nabla}\Omega\times }$
& \hspace*{-2.5mm}$ \oplus A[b-n]$
&\hspace*{-2.5mm}$\xrightarrow{\overrightarrow{\nabla}\Omega\cdot }A
\to A/(\Omega_x,\Omega_y,\Omega_z)\hspace*{-1mm}\to \hspace*{-.8mm}0$.\\
&\hspace*{-2.5mm}$\oplus A[-2n+a+b]$& & \hspace*{-2.5mm}$\oplus A[c-n]$
&
\end{tabular}
\end{equation}

Note that \eqref{E5.0.1} is a complex of graded vector spaces, 
where the differentials are graded maps of degree zero. A 1-cycle 
in $\ker(\overrightarrow{\nabla}\Omega \cdot)$ is called 
{\it sealed} if $\overrightarrow{\nabla}\cdot \overrightarrow{f}=0$ 
when further considered as an element in $A_{sing}$. Let 
$s_1(\Omega)$ be the subspace of 
$\ker(\overrightarrow{\nabla}\Omega \cdot)$ 
consisting of all sealed 1-cycles in the above complex. The 
following follows from commutative algebra. 

\begin{lemma}
\label{zzlem5.1}
Retain the above notation.
\begin{enumerate}
\item[(1)]
If ${\rm gcd}(\Omega_x,\Omega_y,\Omega_z)=1$, then the Koszul 
complex \eqref{E5.0.1} is exact everywhere except possibly for 
the position at $K_1(\overrightarrow{\nabla}\Omega)$.
\item[(2)]
If $\GKdim A_{sing}\leq 1$, then 
${\rm gcd}(\Omega_x,\Omega_y,\Omega_z)=1$.
\item[(3)]
If $\Omega$ is irreducible {\rm{(}}and weighted homogeneous{\rm{)}}, 
then ${\rm gcd}(\Omega_x,\Omega_y,\Omega_z)=1$.
\item[(4)]
$\im(\overrightarrow{\nabla}
\Omega\times)\subseteq s_1(\Omega)
\subseteq \ker(\overrightarrow{\nabla}
\Omega\cdot)$
\end{enumerate}
\end{lemma}

\begin{proof}
(1) It follows from \cite[Remarks 3.6 \& 3.7]{Pi1}. 

(2,3) These are clear. 

(4) For any $\overrightarrow{g}\in A^{\oplus 3}$, by a 
computation, $\overrightarrow{\nabla}\cdot(
\overrightarrow{\nabla}\Omega\times \overrightarrow{g})=0$ 
in $A_{sing}$. The second inclusion follows from the definition.
\end{proof}

\begin{definition}
\label{zzdef5.2}
Let $\Omega\in A$ be a potential, and we consider the 
Koszul complex \eqref{E5.0.1}.
\begin{enumerate}
\item[(1)]
The {\it sealed first Koszul homology} of $(A, \Omega)$ is defined
to be
$$sK_1(A,\Omega):=s_1(\Omega)/\im(\overrightarrow{\nabla}
\Omega\times).$$
\item[(2)]
We say $\Omega$ is {\it $K_1$-sealed} if $sK_1(A,\Omega)=0$. That 
is, for any $\overrightarrow{f}\in A^{\oplus 3}$, if 
$\overrightarrow{f}\cdot \overrightarrow{\nabla}\Omega=0$ in $A$ and 
$\overrightarrow{\nabla}\cdot \overrightarrow{f}=0$ when 
considered as an element in $A_{sing}$, then 
$\overrightarrow{f}
=\overrightarrow{\nabla}\Omega\times \overrightarrow{g}$ 
for some $\overrightarrow{g}\in A^{\oplus 3}$.
\end{enumerate}
\end{definition}

 The property of being $K_1$-sealed was implicitly used by Luo in her thesis \cite{Luo}. It involves computing the Poisson homology using the first homology of the corresponding Koszul complex in certain exceptional cases. It is unclear if the $K_1$-sealedness of $\Omega$ depends on the choice of graded generators $(x,y,z)$. We assume a fixed set of $(x,y,z)$ when discussing $K_1$-sealedness. 

By definition, the $K_1$-sealedness for $\Omega$ can be 
reflected via the homology of the Koszul complex 
$K_\bullet(\overrightarrow{\nabla}\Omega)$ in the following 
way: for any 1-cycle $\overrightarrow{f}\in Z_1
(K_\bullet(\overrightarrow{\nabla}\Omega))$, if 
$\overrightarrow{\nabla}\cdot \overrightarrow{f}=0$ in 
$A_{sing}$, then $\overrightarrow{f}$ belongs to the 
1-boundary, namely $\overrightarrow{f}=0$ in 
$H_1(K_\bullet(\overrightarrow{\nabla}\Omega))$.
If $\Omega$ has an isolated singularity at the origin,
then $H_1(K_\bullet(\overrightarrow{\nabla}\Omega))=0$
\cite[Proposition 3.5]{Pi1}. Hence, such an 
$\Omega$ is always $K_1$-sealed.

In the rest of this subsection, we will show that $K_1$-sealedness holds for some other families of $\Omega$ that do not have isolated singularity. 

\begin{lemma}
\label{zzlem5.3}
 Let $\Omega=xyz+g(x,y)$ for some $g(x,y)\in \kk[x,y]$ satisfying the following conditions. 
 \begin{itemize}
     \item[(1)] $g(x,y)$ is homogeneous with respect to some new grading $\deg_{new}(x)=a'$ and $\deg_{new}(y)=b'$ for $a',b'\ge 3$.
     \item[(2)] $g(x,y)$ contains both terms $x^{b'}$ and $y^{a'}$.
     \item[(3)] $xy\mid g_xg_y$.
     \end{itemize}
Then $\Omega$ is $K_1$-sealed. In this case, by choosing such a new 
grading together with $\deg_{new}(z)=a'b'-a'-b'=:c'$, we have $h_{H_1(K_\bullet)}(t)=\frac{t^{c'+a'b'}}{1-t^{c'}}$. 
\end{lemma}
For example, $\Omega=xyz+x^{b^\prime}+y^{a^\prime}$ with $a^\prime, b^\prime \geq 3$ is $K_1$-sealed. In applications/examples in the next section, we have 
$(a,b,c)=(a'/g, b'/g, c'/g)$ for $g=\gcd(a',b',c')$. 

\begin{proof}[Proof of \Cref{zzlem5.3}]
We assign a new grading, denoted by $\deg_{new}$, on $A$ by 
choosing ${\rm deg}_{new}(x):=a'$, ${\rm deg}_{new}(y):=b'$, and 
${\rm deg}_{new}(z):=a'b'-a'-b'$. It is obvious that $\Omega$ 
becomes a homogeneous potential of degree $n':=a'b'$ under this 
new grading. It is easy to check that $\GKdim A_{sing}=1$ and a 
$\kk$-basis of $A_{sing}$ will be explicitly constructed later 
on. By Lemma \ref{zzlem5.1}, the Koszul complex 
$K_\bullet(\Omega_x,\Omega_y,\Omega_z)$ given in 
\eqref{E5.0.1} is exact everywhere except at $K_1$. We claim 
that $H_1(K_\bullet)$ is spanned by the images of the elements 
\[\overrightarrow{\varphi_l}:
=z^l\left(xz,\, -g_x,\, \frac{g_xg_y}{xy}-z^2\right)
\]
in $\left(A[a'-n']\oplus A[b'-n']\oplus A[c'-n']\right)_{c'(l+1)+n'}$ 
for all $l\ge 0$. It is easy to check that 
$\overrightarrow{\varphi_l}\cdot \overrightarrow{\nabla}\Omega=0$. 
Suppose we have $\overrightarrow{\varphi_l}
=\overrightarrow{f}\times \overrightarrow{\nabla}\Omega$ 
for some $\overrightarrow{f}=(f_1,f_2,f_3)\in 
\left(A[-a'-n']\oplus A[-b'-n']\oplus A[-c'-n']\right)_{c'(l+1)+n'}$. 
Then the third component of $\overrightarrow{\varphi_l}
=\overrightarrow{f}\times \overrightarrow{\nabla}\Omega$ is equal to
\[
f_1(xz+g_y)-f_2(yz+g_x)=\frac{g_xg_y}{xy}\,z^l-z^{l+2}.
\]
However, this is impossible since the monomial $z^{l+2}$ cannot 
appear on the left side. Since each homogeneous element 
$\overrightarrow{\varphi_l}$ has a different degree for distinct 
$l$, their images in $H_1(K_\bullet)$ must be linearly 
independent. To prove our claim, it suffices to match the Hilbert 
series of $H_1(K_\bullet)$ with $\frac{t^{c'+n'}}{1-t^{c'}}$, which 
is the one associated with the $\kk$-subspace spanned by 
\{$\overrightarrow{\varphi_l} \mid l\geq 0\}$. By our assumption (2), we can apply the diamond lemma to obtain a $\kk$-basis of 
$A_{sing}=A/(\Omega_x,\Omega_y,\Omega_y)=\kk[x,y,z]/(xy,xz+g_y,yz+g_x)$ given by 
\begin{align}
\label{E5.3.1}\tag{E5.3.1}
\{z^i\,|\, i\ge 0\}\ \cup\ \{x,\ldots,x^{b'-1}\}\ \cup\ \{y,\ldots,y^{a'-1}\}.\end{align}
Thus, $A_{sing}$ has a Hilbert series given below  
\[
h_{A_{sing}}(t)=\frac{1}{1-t^{c'}}
+\frac{t^{a'}-t^{a'b'}}{1-t^{a'}}
+\frac{t^{b'}-t^{b'a'}}{1-t^{b'}}.
\]
A calculation for the Hilbert series of the terms in  
\eqref{E5.0.1} yields the following
\[
h_{H_1(K_\bullet)}(t)=h_A(t)(t^{a'+b'}-1)(t^{b'+c'}-1)
(t^{a'+c'}-1)+h_{A_{sing}}(t)=\frac{t^{c'+n'}}{1-t^{c'}}.
\]
Therefore, we have proved the claim. 

Next, we show that $\Omega$ is $K_1$-sealed. Taking an 
arbitrary homogeneous element $\overrightarrow{f}\in 
Z_1(K_\bullet(\overrightarrow{\nabla}\Omega))$ satisfying 
$\overrightarrow{\nabla}\cdot \overrightarrow{f}=0$ in 
$A_{sing}$,  up to a boundary, we can take 
$\overrightarrow{f}=\lambda \overrightarrow{\varphi_l}$ for 
some $\lambda \in \kk$. One can check that
\[
\overrightarrow{\nabla}\cdot \overrightarrow{f}=
\lambda (\overrightarrow{\nabla}\cdot \overrightarrow{\varphi_l})
=\lambda\left(l\frac{g_xg_y}{xy}\,z^{l-1}-g_{xy}\,z^l-(l+1)z^{l+1}\right)=0\ \text{in $A_{sing}$}.
\]
By our assumption (2), one can check directly $y\nmid g_x$ and $x\nmid g_y$. So our assumption (3) implies $x\mid g_x$ and $y\mid g_y$. Thus, there is a surjection 
\[\pi: A\twoheadrightarrow A_{sing}=A/(\Omega_x,\Omega_y,\Omega_z)\twoheadrightarrow A/(x,y)\cong \kk[t]\]
via $\pi(x)=\pi(y)=0$ and $\pi(z)=t$. So we have
\begin{align*}
\pi(\overrightarrow{\nabla}\cdot \overrightarrow{f})&=\lambda\left(l\pi(g_x/x)\pi(g_y/y)\,\pi(z)^{l-1}-\pi(g_{xy})\,\pi(z)^l-(l+1)\pi(z)^{l+1}\right)\\
&=\lambda l\pi(g_x/x)\pi(g_y/y)t^{l-1}-\lambda \pi(g_{xy})\,t^l-\lambda (l+1)t^{l+1}
\end{align*}
which cannot be zero in $\kk[t]$ unless $\lambda=0$. Hence 
$\overrightarrow{f}\in B_1(K_\bullet(
\overrightarrow{\nabla}\Omega))$, and consequently,
$sK_1(A,\Omega)=0$. Therefore, $\Omega$ is $K_1$-sealed.              
\end{proof}

Another type of $K_1$-sealed $\Omega$ is given by $\Omega=z^n+g(x,y)$ for some $g(x,y)$ is in $\kk[x,y]$. We need the following definition.

\begin{definition}
    A polynomial $g\in \kk[x,y]$ is called {\it special} if for any polynomial $f\in \kk[x,y]$, $d\mid (\frac{g_y}{d}f)_x-(\frac{g_x}{d}f)_y$ implies that $d\mid f$ where $d=\gcd(g_x,g_y)$. 
\end{definition}

\begin{lemma}\label{lem5.5}
    The following $g\in \kk[x,y]$ are special.
    \begin{itemize}
\item[(1)] $g=x^ky^l$ with $k,l\ge 1$ and $\gcd(k,l)=1$.
\item[(2)] $g=x^ky^l+x^m$ with $l\ge 1$ and $m\ge k\ge 1$. 
\item[(3)] $g=x^ky^l+y^m$ with $k\ge 1$ and $m\ge l\ge 1$. 
    \end{itemize}
\end{lemma}
\begin{proof}
(1) We have $d=\gcd(g_x,g_y)=x^{k-1}y^{l-1}$. Suppose there is some $f\in \kk[x,y]$ such that $x^{k-1}y^{l-1}\mid(\frac{g_y}{d}f)_x-(\frac{g_x}{d}f)_y=(lxf)_x-(kyf)_y=(l-k)f+lxf_x-kyf_y$. Without loss of generality, we can take $f=x^my^n$ to be a monomial. Thus, $x^{k-1}y^{l-1}\mid(l(m+1)-k(n+1))f$. If to the contrary, $x^{k-1}y^{l-1}\nmid f$, then we must have $m<k-1$ or $n<l-1$ and $l(m+1)=k(n+1)$. The latter implies that $l\mid n+1$ and $k\mid m+1$ since $\gcd(k,l)=1$. This yields a contradiction. 

(2) We have $g_x=x^{k-1}(ky^l+mx^{m-k})$ and $g_y=lx^ky^{l-1}$. So $d=\gcd(g_x,g_y)=x^{k-1}$. Suppose there is some $f\in \kk[x,y]$ such that 
\[x^{k-1}\mid(\frac{g_y}{d}f)_x-(\frac{g_x}{d}f)_y=(lxy^{l-1}f)_x-((ky^l+mx^{m-k})f)_y=(1-k)ly^{l-1}f+ly^{l-1}xf_x-(ky^l+mx^{m-k})f_y.
\]
We consider $f$ a polynomial in $x$ with coefficients in $\kk[y]$ and denote its lowest term as $hx^n$ for some $0\neq h\in \kk[y]$. Then the possible lowest term in $(\frac{g_y}{d}f)_x-(\frac{g_x}{d}f)_y$ is $x^n$ with its coefficient given by
\[
-l(k-1-n)y^{l-1}h-k(y^l+\delta_{m, k})h_y.
\]
If $n\geq k-1$, then $d\mid f$. Suppose $n<k-1$. Then it implies that the above coefficient is zero and so we get  $h_y/h=-\frac{k-1-n}{k}ly^{l-1}/(y^l+\delta_{m, k})$. Thus $h=\lambda/ (y^l+\delta_{m, k})^{(k-1-n)/k}$ for some non-zero scalar $\lambda$, which is not a polynomial. This gives a contradiction.

(3) The proof is similar to (2). 
\end{proof}

\begin{lemma}\label{lem5.6}
Let $\Omega=z^n+g(x,y)$ where $n\ge 2$ and $g(x,y)\in \kk[x,y]$. Then $\Omega$ is $K_1$-sealed if $g$ is special.
\end{lemma}
\begin{proof}
Let $\overrightarrow{f}=(f_1,f_2,f_3)\in A^{\oplus 3}$ satisfying $\overrightarrow{\nabla}\Omega\cdot \overrightarrow{f}=0$ in $A$ and $\overrightarrow{\nabla}\cdot \overrightarrow{f}=0$ in $A_{sing}$. Write $f_i=\sum_{j=0}^{n-1} h_{ij}z^j$ where $h_{ij}\in \kk[x,y]$ for $0\le j\le n-2$ and $h_{i(n-1)}\in \kk[x,y,z]$ for $i=1,2$. Thus
\begin{align*}
 \overrightarrow{\nabla}\Omega\cdot \overrightarrow{f}=   \sum_{i=0}^{n-2}(h_{1i}g_x+h_{2i}g_y)z^i+(h_{1(n-1)}g_x+h_{2(n-1)}g_y+nf_3)z^{n-1}=0.
\end{align*}
So we get $h_{1i}g_x+h_{2i}g_y=0$ in $\kk[x,y]$ for $0\le i\le n-2$ and $f_3=-(h_{1(n-1)}g_x+h_{2(n-1)}g_y)/n$. Set $d=\gcd(g_x,g_y)$, we can further write $h_{1i}=(g_y/d)l_i$ and $h_{2i}=-(g_x/d)l_i$ with $l_i\in \kk[x,y]$ for all $0\le i\le n-2$. Therefore,
\begin{align*}
\overrightarrow{\nabla}\cdot\overrightarrow{f}&=\overrightarrow{\nabla}\cdot\left(\sum_{i=0}^{n-2} l_i\left(\frac{g_y}{d},-\frac{g_x}{d},0\right)z^i+h_{1(n-1)}\left(z^{n-1},0,-\frac{g_x}{n}\right)+h_{2(n-1)}\left(0,z^{n-1},-\frac{g_y}{n}\right)\right)\\
&=\overrightarrow{\nabla}\cdot\left(\sum_{i=0}^{n-2} l_i\left(\frac{g_y}{d},-\frac{g_x}{d},0\right)z^i+\overrightarrow{\nabla}\Omega\times\left(\frac{h_{2(n-1)}}{n},-\frac{h_{1(n-1)}}{n},0\right)\right)\\
&=\sum_{i=0}^{n-2}\left((\frac{g_y}{d}l_i)_x-(\frac{g_x}{d}l_i)_y\right)z^i-\overrightarrow{\nabla}\Omega\, \cdot\, \left(\overrightarrow{\nabla}\times \left(\frac{h_{2(n-1)}}{n},-\frac{h_{1(n-1)}}{n},0\right)\right)\\
&=\sum_{i=0}^{n-2}\left((\frac{g_y}{d}l_i)_x-(\frac{g_x}{d}l_i)_y\right)z^i\  \text{in $A_{sing}$}\\
&=0\  \text{in $A_{sing}$}.
\end{align*}
Note that $A_{sing}=\kk[x,y,z]/(g_x,g_y,z^{n-1})=\oplus_{i=0}^{n-2} (\kk[x,y]/(g_x,g_y))\, z^i$. This implies that $(\frac{g_y}{d}l_i)_x=(\frac{g_x}{d}l_i)_y$ in $\kk[x,y]/(g_x,g_y)$, and hence $d\mid(\frac{g_y}{d}l_i)_x-(\frac{g_x}{d}l_i)_y$ for all $0\le i\le n-2$. Since $g$ is special, by definition, we can write $l_i=dm_i$ for $m_i\in \kk[x,y]$ for all $0\le i\le n-2$, and hence 
\begin{align*}
   \overrightarrow{f}&=\sum_{i=0}^{n-2} l_i\left(\frac{g_y}{d},-\frac{g_x}{d},0\right)z^i+\overrightarrow{\nabla}\Omega\times\left(\frac{h_{2(n-1)}}{n},-\frac{h_{1(n-1)}}{n},0\right)\\
   &=\sum_{i=0}^{n-2} \left(g_y,-g_x,0\right)m_iz^i+\overrightarrow{\nabla}\Omega\times\left(\frac{h_{2(n-1)}}{n},-\frac{h_{1(n-1)}}{n},0\right)\\
&=\overrightarrow{\nabla}\Omega\times\left(\frac{h_{2(n-1)}}{n},-\frac{h_{1(n-1)}}{n},\sum_{i=0}^{n-2}m_iz^i\right).
\end{align*}
This proves our result. 
\end{proof}

\subsection{$uPH^2$-vacancy}
In this subsection, we introduce another new concept: $uPH^2$-vacancy 
for $A_\Omega$ and we shall show that if $\Omega$ is $K_1$-sealed, 
then its corresponding Poisson algebra $A_{\Omega}$ is $uPH^2$-vacant. 
By definition,
$$\ker(d_{\pi_\Omega}^2)~=
~\{\delta\in \mathfrak X^2(A)\,|\, [\delta,\pi_\Omega]_S=0\}$$
where $d_{\pi_\Omega}^\bullet$ is the differential
in the cochain complex 
$(\mathfrak X^\bullet(A_\Omega),d_{\pi_\Omega}^\bullet)$
\eqref{E1.3.1}. Before stating the definitions, we need a lemma. 

\begin{definition-lemma}
\label{zzlem5.7}
Retain the above notation. Let \[M^2(A):=
\{f\pi_{\Omega}+\pi_g\,|\, f,g\in A\}\] which is 
a subspace of $\mathfrak X^2(A)$. Then 
$\im (d_{\pi_\Omega}^1) \subseteq M^2(A)\subseteq 
\ker(d_{\pi_\Omega}^2)$.
\end{definition-lemma}

\begin{proof}
The first inclusion follows directly from \eqref{E4.1.1}. 
That is, $d_\pi^1(\delta)=~\divv(\delta)\,
\pi_{\Omega}-\pi_{\delta(\Omega)}$ for any $\delta\in 
{\mathfrak X}^1(A)$. For the second one, it suffices to 
show that $u:=[f \pi_{\Omega}, \pi_{\Omega}]_S=0$ and 
$v:=[\pi_g,\pi_{\Omega}]_S=0$ for all $f,g\in A_{\Omega}$. 
Note that both $u$ and $v$ are in ${\mathfrak X}^3(A)$. 
Hence, it remains to show that $u(x,y,z)=0$ and 
$v(x,y,z)=0$. They follow from the definition of the 
Schouten bracket and \eqref{E3.1.1} via a straightforward 
computation. We omit the details. 
\end{proof}

\begin{definition} 
\label{zzdef5.8}
Let $A_\Omega$ be the unimodular Poisson algebra defined in 
Notation \ref{zznot3.2}.
\begin{enumerate}
\item[(1)]
The {\it upper division of the second Poisson cohomology} of 
$A_{\Omega}$ is defined to be
$$uPH^2(A_{\Omega}):=\ker(d_{\pi_\Omega}^2)/M^2(A).$$
\item[(2)]
The {\it lower division of the second Poisson cohomology} of 
$A_{\Omega}$ is defined to be
$$lPH^2(A_{\Omega}):=M^2(A)/\im (d_{\pi_\Omega}^1).$$
\item[(3)]
We say $A_{\Omega}$ is {\it $uPH^2$-vacant} if 
$uPH^2(A_{\Omega})=0$, or equivalently 
$lPH^2(A_{\Omega})=PH^2(A_{\Omega})$.
\end{enumerate}
\end{definition}

According to \Cref{zzdef3.1}, the definition of $uPH^2(A_{\Omega})$ (and $lPH^2(A_{\Omega})$ as well as $uPH^2$-vacancy) is invariant under the choice of graded generators $(x,y,z)$. Consequently, $uPH^2$-vacancy is independent of the choice of graded generators $(x,y,z)$. According to the following lemma, $lPH^2(A_{\Omega})$ is typically non-zero.

\begin{lemma}
\label{zzlem5.9}
Retain the above notation and assume (1), (2) and (3) of \Cref{zzhyp0.1}. 
If $h_{PH^0(A)}(t)=h_{PH^1(A)}(t)=\frac{1}{1-t^n}$, then 
the Hilbert series of $lPH^2(A_\Omega)$ is given by 
\begin{equation*}
   h_{lPH^2(A)}(t)~=~\frac{1}{ t^{a+b+c}}\left(\frac{(1-t^{n-a})
(1-t^{n-b})(1-t^{n-c})}{(1-t^n)(1-t^a)(1-t^b)(1-t^c)}-1\right). 
\end{equation*}
\end{lemma}

\begin{proof}
We claim that the following sequence of graded vector spaces 
\begin{equation}
\label{E5.9.1}\tag{E5.9.1}
0\to \kk[\Omega][a+b+c]\xrightarrow{\alpha} 
A[-n+a+b+c]\oplus A[a+b+c]\xrightarrow{\beta} 
M^2(A)\xrightarrow{}0
\end{equation}
is exact, where $\alpha(g)=(-dg/d\Omega,g)$ and 
$\beta(f,g)=f\pi_\Omega+\pi_g$. It is clear that 
$\alpha,\beta$ are graded maps, $\alpha$ is injective 
and $\beta$ is surjective, and $\beta\alpha=0$. Hence 
it remains to show that 
$\ker (\beta)=\im (\alpha)$. Suppose 
$\beta(f,g)=f\pi_\Omega+\pi_g=0$. One can check that
\begin{equation*}
d_\pi^0(g)=[\pi_\Omega,g]_S=-[\pi_g,\Omega]_S=[f\pi_\Omega,\Omega]_S=0.
\end{equation*}
By the assumption, we have $Z_P(A)=\kk[\Omega]$. Therefore, we have $g\in Z_P(A)=\kk[\Omega]$ and $f=-dg/d\Omega$. 
Note that we have the following exact sequences of graded vector 
spaces:
\begin{equation}
\label{E5.9.2}\tag{E5.9.2}
0\to PH^0(A)\to \mathfrak X^0(A)\to \im(d_\pi^0)\to 0,
\end{equation}
\begin{equation}
\label{E5.9.3}\tag{E5.9.3}
0\to  \im(d_\pi^0)\to Pd(A)[w]\to PH^1(A)[w]\to 0,
\end{equation}
\begin{equation}
\label{E5.9.4}\tag{E5.9.4}
0\to Pd(A)[w]\to \mathfrak X^1(A)[w]\to M^2(A)[2w]
\to lPH^2(A)[2w]\to 0.
\end{equation}
One can deduce that 
\begin{equation}
\label{E5.9.5}\tag{E5.9.5}
h_{lPH^2(A)}(t)=h_{M^2(A)}(t)-h_{\mathfrak X^1(A)}(t)t^w
+h_{\mathfrak X^0(A)}(t)t^{2w}
+h_{PH^1(A)}(t)t^w-h_{PH^0(A)}(t)t^{2w}.
\end{equation}
The assertion follows from \eqref{E5.9.1} via a 
direct computation.
\end{proof}

We also need the {\it de Rham complex} for $A$: 
\begin{equation}
\label{E5.9.6}\tag{E5.9.6}
0\to \kk \to \Omega^0_A\xrightarrow{d} 
\Omega^1_A\xrightarrow{d} \Omega^2_A\xrightarrow{d} 
\Omega^3_A\to 0.  
\end{equation}
We will utilize the following natural isomorphisms of 
graded vector spaces 
\begin{equation}
\label{E5.9.7}\tag{E5.9.7}
\left \{
\centering
\begin{aligned}
&\Omega^0_A\xrightarrow{\sim}A &&\\
&\Omega^1_A\xrightarrow{\sim} A[-a]\oplus A[-b]\oplus A[-c]  
&&  fdg \mapsto f\overrightarrow{\nabla} g \\
&\Omega^2_A\xrightarrow{\sim} A[-b-c]\oplus A[-a-c]\oplus A[-a-b]
&& fdg\wedge dh\mapsto f\overrightarrow{\nabla} g\times 
   \overrightarrow{\nabla} h\\
&\Omega^3_A\xrightarrow{\sim} A[-a-b-c]  
&& fdx\wedge dy\wedge dz\mapsto f\\
\end{aligned}\right.
\end{equation}
to \eqref{E5.9.6} for $f, g,h\in A$. As a result, we have 
the following complex of graded vector spaces:  
\begin{equation}
\label{E5.9.8}\tag{E5.9.8}
\centering
\begin{tabular}{ccccc}
&\,$A[-a]$ & &$A[-b-c]$& \vspace*{-2mm}\\
$0\to \kk\to A\xrightarrow{\overrightarrow{\nabla}}$
& \hspace*{-2.5mm}$\oplus A[-b]$
&\hspace*{-2.5mm}$\xrightarrow{\overrightarrow{\nabla}\times }$
& \hspace*{-2.5mm}$ \oplus A[-a-c]$
&\hspace*{-2.5mm}$\xrightarrow{\overrightarrow{\nabla}\cdot }
A[-a-b-c]\to 0.$\\
&\hspace*{-2.5mm}$\oplus A[-c]$& & \hspace*{-2.5mm}$
\oplus A[-a-b]$&
\end{tabular}
\end{equation}
In the sequel, we will use the well-known fact that the de Rham 
complex \eqref{E5.9.8} for $A$ is always exact. 

The following lemma shows that $A_{\Omega}$ is $uPH^2$-vacant if and only if it is $H$-ozone. It is important to note that the $uPH^2$-vacancy of $A_{\Omega}$ plays a significant role in computing the Hilbert series of Poisson cohomology groups. In this lemma, we will give two more equivalent 
descriptions in terms of a combination of the Koszul complex for 
the sequence 
$\overrightarrow{\nabla}\Omega=(\Omega_x, \Omega_y, \Omega_z)$ 
and the de Rham complex for $A$. 

\begin{lemma}
\label{zzlem5.10}
Let $A$ be $A_\Omega$ satisfying (1), (2) and (3) of \Cref{zzhyp0.1}
where $\Omega$ is a nonconstant homogeneous potential. Then 
the following are equivalent.
\begin{itemize}
\item[(1)] 
$A$ is $H$-ozone.
\item[(2)] 
For any $\overrightarrow{f}\in A^{\oplus 3}$, 
$\overrightarrow{\nabla}\cdot \overrightarrow{f}
=\overrightarrow{f}\cdot \overrightarrow{\nabla}\Omega=0$ 
implies that $\overrightarrow{f}=\overrightarrow{\nabla}g
\times \overrightarrow{\nabla}\Omega$ for some $g\in A$.
\item[(3)] 
For any $\overrightarrow{f}\in A^{\oplus 3}$, if 
$(\overrightarrow{\nabla}\times \overrightarrow{f})\cdot 
\overrightarrow{\nabla}\Omega=0$, then 
$\overrightarrow{f}=\overrightarrow{\nabla}m
+g\overrightarrow{\nabla}\Omega$ for some $m,g\in A$.
\item[(4)] 
$A$ is $uPH^2$-vacant.
\end{itemize}
\end{lemma}

\begin{proof}
(1)$\Leftrightarrow$(2) Assume (2) holds and suppose 
$\delta$ is an ozone derivation of $A$, i.e., 
$\delta(\Omega)=0$. Write $\overrightarrow{\delta}
=(\delta(x),\delta(y),\delta(z))\in A^{\oplus 3}$. 
So $0=\delta(\Omega)=\overrightarrow{\delta} \cdot 
\overrightarrow{\nabla}\Omega$. By \eqref{E4.1.1}, 
we have $\overrightarrow{\nabla}\cdot 
\overrightarrow{\delta}={\rm div}(\delta)=0$ as 
$\delta(\Omega)=0$. So (2) implies that 
$\overrightarrow{\delta}=\overrightarrow{\nabla}g\times 
\overrightarrow{\nabla}\Omega$ for some $g\in A$, which 
is equivalent to $\delta=\{-,g\}$. So (1) holds. 

Conversely, suppose there is some $\overrightarrow{f}
\in A^{\oplus 3}$ such that $\overrightarrow{\nabla}\cdot 
\overrightarrow{f}=\overrightarrow{f}\cdot 
\overrightarrow{\nabla}\Omega=0$. Consider the 
derivation $\delta$ of $A$ defined by 
$(\delta(x),\delta(y),\delta(z))=\overrightarrow{f}$.  
The assumptions on $\overrightarrow{f}$ yield that 
$\delta(\Omega)=\divv(\delta)=0$. So \Cref{zzlem4.1}(2) 
implies that $\delta$ is a Poisson derivation. Note that 
in the proof of \cite[Lemma 1]{MTU}, it is shown that the 
Poisson center $Z_P(A)$ is algebraic over $\kk[\Omega]$. 
Hence, $\delta$ vanishes on $Z_P(A)$, so it is ozone. 
So by (1), $\delta=\{-,g\}$ for some $g\in A$,  that is, 
$\overrightarrow{f}=\overrightarrow{\nabla}g\times 
\overrightarrow{\nabla}\Omega$.

(2)$\Leftrightarrow$(3) Assume (2) holds and suppose that 
$(\overrightarrow{\nabla}\times \overrightarrow{f})\cdot 
\overrightarrow{\nabla}\Omega=0$ for  
$\overrightarrow{f}\in A^{\oplus 3}$. Write 
$\overrightarrow{h}:=\overrightarrow{\nabla}\times 
\overrightarrow{f}\in A^{\oplus 3}$. It is easy to see that 
$\overrightarrow{\nabla}\cdot \overrightarrow{h}
=0$ and $\overrightarrow{h}\cdot \overrightarrow{\nabla}\Omega=0$. 
So $\overrightarrow{h}=\overrightarrow{\nabla}g\times 
\overrightarrow{\nabla}\Omega=\overrightarrow{\nabla}
\times(g\overrightarrow{\nabla}\Omega)$ for some $g\in A$. 
This implies that $\overrightarrow{\nabla}\times
(\overrightarrow{f}-g\overrightarrow{\nabla}\Omega)=0$. 
By the exactness of the de Rham complex for $A$ \eqref{E5.9.8}, we can 
write  $\overrightarrow{f}-g\overrightarrow{\nabla}\Omega
=\overrightarrow{\nabla}m$ for some $m\in A$. So (3) holds. 

Conversely, suppose $\overrightarrow{\nabla}\cdot 
\overrightarrow{f}=0$ and $\overrightarrow{f}\cdot
\overrightarrow{\nabla}\Omega=0$ for some $\overrightarrow{f}\in A^{\oplus 3}$. 
By the exactness of the de Rham complex for $A$ \eqref{E5.9.8}, 
we can write $\overrightarrow{f}=\overrightarrow{\nabla}
\times \overrightarrow{h}$ for some $\overrightarrow{h}\in A^{\oplus 3}$. 
Hence $(\overrightarrow{\nabla}\times \overrightarrow{h})
\cdot \overrightarrow{\nabla}\Omega=0$. So, by (3), we can 
write  $\overrightarrow{h}=\overrightarrow{\nabla}m+g
\overrightarrow{\nabla}\Omega$ for some $g, m\in A$. Thus 
$\overrightarrow{f}=\overrightarrow{\nabla}\times 
(\overrightarrow{\nabla}m+g\overrightarrow{\nabla}\Omega)
=\overrightarrow{\nabla}g\times \overrightarrow{\nabla}\Omega$ 
and (2) follows. 

(3)$\Leftrightarrow$(4) We use the identification 
$\mathfrak X^2(A)\xrightarrow{\sim} A^{\oplus 3}$ described in 
\eqref{E1.3.2} via 
$$\delta\mapsto \overrightarrow{\delta}
=(\delta(y,  z), \delta(z, x), \delta(x, y)),$$
which is different from the notation $\overrightarrow{\delta}$ 
used in the proof of part (1). 
One can check that $\pi_g$ and $f\pi_\Omega$ correspond to 
$\overrightarrow{\nabla}g$ and $f\overrightarrow{\nabla}\Omega$, 
respectively. Moreover by \eqref{E1.3.6}, we have 
$-[\delta,\pi_\Omega]_S=-(\overrightarrow{\nabla}\times 
\overrightarrow{\delta})\cdot \overrightarrow{\nabla}\Omega$. 
Then, (3) and (4) are equivalent by reinterpreting the conditions 
through the above identification. 
\end{proof}

The result below links the $K_1$-sealness of a potential $\Omega$ to the $H$-ozoness of its associated Poisson algebra, $A_\Omega$.  

\begin{lemma}
\label{zzlem5.11}
Let $\Omega$ be a homogeneous polynomial 
of positive degree $n$. If $\Omega$ is $K_1$-sealed, then 
$A_\Omega$ is H-ozone.
\end{lemma} 

\begin{proof}
It suffices to show that $A$ satisfies condition (3) in 
\Cref{zzlem5.10}. Suppose that $\deg(x)=a,\deg(y)=b$, and 
$\deg(z)=c$ for some positive integers $a,b,c$.  Note that 
$\Omega$ is homogeneous with a positive degree $n$, which 
is not necessarily equal to $a+b+c$. As a result, the 
Poisson bracket $\pi_\Omega$ on $A$ is homogeneous of 
degree $n-a-b-c$. The following diagram combines the Koszul complex \eqref{E5.0.1} and the de Rham 
complex \eqref{E5.9.8}.
\[
\xymatrix{
& \left(A[-a] \oplus A[-b] \oplus A[-c]\right)[-n]
\ar[d]^-{\overrightarrow{\nabla}\Omega\times}&\\
A[-a]\oplus A[-b]\oplus A[-c]\ar[r]^-{\overrightarrow{\nabla}\times}
& A[-b-c]\oplus A[-a-c]\oplus A[-a-b]
\ar[d]^-{\overrightarrow{\nabla}\Omega\cdot}
\ar[r]^-{\overrightarrow{\nabla}\cdot}& A[-a-b-c]&\\
&A[n-a-b-c].&
}
\]
Without loss of generality, let us assume that 
$\overrightarrow{f}\in \left(A[-a]\oplus A[-b]
\oplus A[-c]\right)_\ell$ for some $\ell\in \mathbb Z$ 
such that $(\overrightarrow{\nabla}\times 
\overrightarrow{f})\cdot \overrightarrow{\nabla}\Omega=0$. 
The condition $(\overrightarrow{\nabla}\times 
\overrightarrow{f})\cdot \overrightarrow{\nabla}\Omega=0$ 
implies that $\overrightarrow{\nabla}\times 
\overrightarrow{f}\in Z_1(K_\bullet(\overrightarrow{\nabla}\Omega))$. 
So we can write 
\[\overrightarrow{\nabla}\times \overrightarrow{f}
=\overrightarrow{f_1}\times \overrightarrow{\nabla}\Omega
+\overrightarrow{h}\]
for some $\overrightarrow{f_1}\in \left(A[-a]\oplus A[-b]
\oplus A[-c]\right)_{\ell-n}$ and 
$\overrightarrow{h}\in Z_1(K_\bullet(
\overrightarrow{\nabla}\Omega))\setminus 
B_1(K_\bullet(\overrightarrow{\nabla}\Omega))$ or zero. 
If $\overrightarrow{h}\ne 0$, then 
\begin{align*}
\overrightarrow{\nabla}\cdot  \overrightarrow{h}
&=  \overrightarrow{\nabla}\cdot(\overrightarrow{\nabla}\times 
\overrightarrow{f}-\overrightarrow{f_1}\times 
\overrightarrow{\nabla}\Omega)\\
&=-\overrightarrow{\nabla}\cdot(\overrightarrow{f_1}\times 
\overrightarrow{\nabla}\Omega)\\
&=-(\overrightarrow{\nabla}\times \overrightarrow{f_1})\cdot 
\overrightarrow{\nabla}\Omega\\
&=0\ \text{in $A_{sing}$}.
\end{align*}
The last equality follows from a direct computation.
Since $\Omega$ is $K_1$-sealed, 
$\overrightarrow{h}\in B_1(K_\bullet(\overrightarrow{\nabla}\Omega))$, 
yielding a contradiction. Thus $\overrightarrow{h}=
\overrightarrow{0}$, and consequently,
$(\overrightarrow{\nabla}\times 
\overrightarrow{f_1})\cdot \overrightarrow{\nabla}\Omega=0$
by the above calculation. Repeating this procedure with 
initial setting $f_0=f$, we get 
\[\overrightarrow{\nabla}\times \overrightarrow{f_{m-1}}
=\overrightarrow{f_m}\times \overrightarrow{\nabla}\Omega
\quad{\text{and}}\quad
(\overrightarrow{\nabla}\times 
\overrightarrow{f_m})\cdot \overrightarrow{\nabla}\Omega=0
\]
for some $\overrightarrow{f_m}\in
\left(A[-a]\oplus A[-b]\oplus A[-c]\right)_{\ell-mn}$ for 
$m\geq 1$. 
Since $A$ is connected graded, we have $\overrightarrow{f_p}=0$ 
for some $p\gg 0$, which implies that 
$\overrightarrow{\nabla}\times \overrightarrow{f_{p-1}}=0$. 
By the exactness of the de Rham complex for $A$, we have 
$\overrightarrow{f_{p-1}}=\overrightarrow{\nabla}G_{p-1}$ 
for some $G_{p-1}\in A$. To complete the proof, it suffices 
to prove that, if we can write $\overrightarrow{f_{q}}=\overrightarrow{\nabla}G_q
+H_q\overrightarrow{\nabla}\Omega$ for some $G_q,H_q\in A$, 
then so can we do the same for $\overrightarrow{f_{q-1}}$. An 
easy calculation yields that 
\begin{align*}
\overrightarrow{\nabla}\times \overrightarrow{f_{q-1}}
&=\overrightarrow{f_q}\times \overrightarrow{\nabla}\Omega\\
&=(\overrightarrow{\nabla}G_q+H_q\overrightarrow{\nabla}\Omega)
  \times \overrightarrow{\nabla}\Omega\\
&=\overrightarrow{\nabla}G_q\times \overrightarrow{\nabla}\Omega\\
&=\overrightarrow{\nabla}\times (G_q\overrightarrow{\nabla}\Omega).
\end{align*}
This implies that $\overrightarrow{\nabla}\times
(\overrightarrow{f_{q-1}}-G_q\overrightarrow{\nabla}\Omega)=0$. 
By the exactness of the de Rham complex for $A$, we can write 
$\overrightarrow{f_{q-1}}-G_q\overrightarrow{\nabla}\Omega
=\overrightarrow{\nabla}G_{q-1}$ for some $G_{q-1}\in A$. 
Now, our result is deduced by a downward induction.  
\end{proof}

By the above two lemmas (and Theorem \ref{zzthm0.6}),  we shall have 
\begin{equation}
\label{E5.11.1}\tag{E5.11.1}
{\text{$\Omega$ is $K_1$-sealed}}
\Rightarrow
{\text{$A_{\Omega}$ is $uPH^2$-vacant}}
\Leftrightarrow
{\text{$A_{\Omega}$ is $H$-ozone}}
\Leftrightarrow
{\text{$\Omega$ is irreducible and balanced}}
\end{equation}

However, it is not clear if
``${\text{$\Omega$ is $K_1$-sealed}}
\Leftarrow
{\text{$A_{\Omega}$ is $uPH^2$-vacant}}$''. 
So, we ask the following question.  

\begin{question}
\label{zzque5.12}
Is $K_1$-sealing condition equivalent to $H$-ozone for $\Omega$ with $|\Omega|=|x|+|y|+|z|$?
\end{question}

\section{Poisson cohomology of weighted Poisson algebras} 
\label{zzsec6}
In this last section, we investigate Poisson cohomology for a graded unimodular Poisson algebra $A$ in dimension three assuming \Cref{zzhyp0.1} unless mentioned otherwise. We aim to expand upon the results of Van den Bergh \cite{VdB} and Pichereau \cite{Pi1} to encompass a wider range of Poisson algebras. First, we show that $A_\Omega$ is $H$-ozone for all balanced irreducible $\Omega$. When $\Omega$ has an isolated singularity, $\Omega$ is $K_1$-sealed and $A_{\Omega}$ is thus $H$-ozone. By \Cref{zzthm3.4} and \Cref{zzlem1.3}, we can use the classification of $\Omega$ in 
\Cref{zzthm2.5} to further divide the rest of the irreducible ones (without isolated singularity) into the following four cases: 
\begin{itemize}
\item[(1)] Under a new grading, the associated graded Poisson algebra of $A_{\Omega}$ is $A_{\overline{\Omega}}$ for some $\overline{\Omega}=z^n+x^ky^l$ with $n\geq 2$ and $\gcd(k,l)=1$ up to a permutation of $x,y,z$; that is, $\Omega$ is one of these irreducible potentials not included in Cases (2)-(4) below.
\item[(2)] $\Omega=xyz+x^{b'}+y^{a'}$, where $a',b'\ge 3$.
\item[(3)] $\Omega$ is non-balanced and irreducible; that is,  $\Omega=z^2+y^3$, or $z^2+x^{2+\frac{2b}{a}}$, or $x^{\frac{b}{{\rm gcd}(a, b)}}+y^{\frac{a}{{\rm gcd}(a, b)}}$.
\item[(4)] $\Omega=z^2+y^3+x^2 y^2$ or $z^2+x^2y^2+x^{2+\frac{2b}{a}}$.
\end{itemize}

\subsection{Case (1)}
\label{zzsec6.1}
We first verify that $A_\Omega$ is $H$-ozone in the case 
$\Omega=z^n+x^ky^l$ with $n\geq 2$ and ${\rm gcd}(k, l)=1$. For any general 
$\Omega$ containing $z^n+x^ky^l$, we construct a 
$w$-filtration on $A_\Omega$ and consider the associated 
$w$-graded Poisson algebra instead. Applying a spectral sequence 
argument, we can establish that $A_\Omega$ is still 
$H$-ozone for any such general $\Omega$.

\begin{lemma}
\label{zzlem6.1}
Let $\Omega=z^n+x^ky^l$ for some positive integers $n, k,l$ such 
that $n\geq 2$ and $\gcd(k,l)=1$. Then $A_{\Omega}$ is H-ozone.
\end{lemma}

\begin{proof}
By \Cref{lem5.5}(1) and \Cref{lem5.6}, $\Omega$ is $K_1$-sealed. The assertion follows
from \Cref{zzlem5.11}.
\end{proof}

Now suppose $\Omega$ is as in Case (1). We will 
consider some Poisson $w$-filtration on $A_\Omega$, whose 
associated $w$-graded Poisson algebra is defined by a potential 
of the form $z^n+x^ky^l$. 

For the rest of this subsection, it is more convenient to use 
the filtration $\mathbb F=\{F_i\,|\, i\in \mathbb Z\}$ on 
$A$ consisting of an increasing chain of $\kk$-subspaces 
$F_i\subseteq F_{i+1}$ (which is different from the one given in
\Cref{zzdef1.5}). Accordingly, we have to modify other parts of
\Cref{zzdef1.5}. In particular, we will change $w$ to $-w$ 
when we use a Poisson $w$-filtration.

\begin{lemma}
\label{zzlem6.2}
Let $A$ be a connected graded polynomial Poisson algebra and 
$\mathbb F$ be a $w$-filtration of $A$ with $w$ being not necessarily zero. Suppose the following 
hold: 
\begin{itemize}
\item[(1)] 
The Poisson center $Z_P(A)=\kk[\chi]$ with $\deg(\chi)>0$.
\item[(2)] 
The Euler derivation $E$ of $A$ preserves the 
$w$-filtration $\mathbb F$.  
\item[(3)] 
The associated graded algebra ${\rm gr}_\mathbb FA$ 
is a  connected $w$-graded polynomial Poisson algebra. 
\item[(4)] 
The Poisson center $Z_P({\rm gr}_\mathbb FA)$ is 
$\kk[\overline{\chi}]$ with $s:=\deg_{new}(\overline{\chi})>0$ 
here $\deg_{new}$ is the new grading associated to the 
filtration ${\mathbb F}$.
\item[(5)] 
${\rm gr}_\mathbb FA$ has no non-zero Poisson derivation of degree $-s$. 
\item[(6)] 
${\rm gr}_\mathbb FA$ is $H$-ozone.
\end{itemize}
Then $A$ is both $H$-ozone and $PH^1$-minimal. 
\end{lemma}

\begin{proof}
Let $E$ be the Euler derivation of $A$. By (2), we have 
the induced graded Poisson derivation, denoted by $E^{ind}$, 
on ${\rm gr}_\mathbb FA$. So we can write 
$E^{ind}(\overline{\chi})=m\,\overline{\chi}$, where 
$m=\deg(\chi)$ is under the original grading of $A$. 
Denote $Z=Z_P(\gr_\mathbb F A)$. We claim that 
$PH^1(\gr_\mathbb F A)\cong ZE^{ind}$ by following the argument 
in \cite[Lemma 7.7]{TWZ}. Note that the induced Poisson 
bracket on ${\rm gr}_\mathbb FA$ is homogeneous of degree $w$. 
So, it suffices to consider all homogeneous Poisson derivations. 
Say $\phi$ is one of degree $i$. One can check that 
$\phi(\overline{\chi})\in Z$. By (5), we get 
$\phi(\overline{\chi})=a \overline{\chi}^n$ for some $a\in \kk$ 
and $n>0$ (and further $\phi(\overline{\chi})=0$ if $s\nmid i$). 
Write $\phi'=\phi-\frac{a}{m}\overline{\chi}^{n-1}E^{ind}$. Then 
$\phi'(\overline{\chi})=0$, whence $\phi'$ is ozone. Now (6) 
implies that $\phi=\frac{a}{m}\overline{\chi}^{n-1}E^{ind}+\phi'$ 
for some Hamiltonian derivation $\phi'$. This means that 
$Pd({\rm gr}_\mathbb FA)=ZE^{ind}+Hd({\rm gr}_\mathbb FA)$. It 
remains to show that $ZE^{ind}\cap Hd({\rm gr}_\mathbb FA)=0$. 
Let $\phi=fE^{ind}$ be Hamiltonian for some $f\in Z$. So 
$\phi(\overline{\chi})=fE^{ind}(\overline{\chi})=mf\overline{\chi}=0$. 
Hence, $f=0$ and $\phi=0$ since $Z$ is an integral domain. 
This proves our claim. 

Next, we use the $w$-filtration 
$\mathbb F=\{F_iA\mid i\in \mathbb Z\}$ of $A$ to filter 
the cochain complex  $(\mathfrak X^\bullet(A), d_\pi^\bullet)$ 
and compute $PH^1(A)$ by spectral sequence. For each 
$p,i\in \mathbb Z$ with $i\ge 0$, we define a $\kk$-subspace 
$F_p\mathfrak X^i(A)$ of $\mathfrak X^i(A)$ 
\begin{align*}
F_p\mathfrak X^i(A)=\{f\in \mathfrak X^i(A)\,|\, 
f(a_{1},\ldots,a_{i})\in F_{l_1+\cdots +l_i-p+iw}A\ 
\text{for any $a_{j}\in  F_{l_j}A$},
1\leq j\leq i\}. 
\end{align*}
Applying the differential formula \eqref{E1.0.1}, for any 
$f\in F_p\mathfrak X^i(A)$ we have
\begin{align*}
d_{\pi}^i(f)(a_0,\ldots,a_i)
&=\sum_{j=0}^i (-1)^j \{ a_j, f(a_0,\ldots,
   \widehat{a_j},\ldots,a_i)\}\\
&\quad +\sum_{0\leq j<k\leq i}
(-1)^{j+k}
f(\{a_j,a_k\}, a_0,\ldots,\widehat{a_j},\ldots,
   \widehat{a_k},\ldots,a_i)\\
&\in F_{l_0+\cdots+l_i-p+(i+1)w}A
\end{align*}
where $a_j\in F_{l_j}A$ for $0\leq j\leq i$. So 
$d_\pi^i: F_p\mathfrak X^i(A)\to F_p\mathfrak X^{i+1}(A)$. 
Then $\{F_p\mathfrak X^\bullet(A)\,|\, p\in \mathbb Z\}$ 
is a (decreasing) filtration on the cochain complex 
$(\mathfrak X^\bullet(A),d_\pi^\bullet)$, which is exhaustive 
and bounded below since ${\rm gr}_\mathbb FA$ is connected 
graded. Thus according to \cite[\S 5.4]{We}, we have a 
cohomology spectral sequence with 
\[E_0^{p,q}=F_p\mathfrak X^{p+q}(A)/F_{p+1}\mathfrak X^{p+q}
\cong \mathfrak X^{p+q}({\rm gr}_\mathbb FA)_{-p+(p+q)w}\]
and 
\[
E_1^{p,q}=PH^{p+q}({\rm gr}_\mathbb FA)_{-p+(p+q)w}
\Longrightarrow PH^{p+q}(A).
\] 
By (4), we know $PH^0({\rm gr}_\mathbb FA)=\kk[\overline{\chi}]$. 
Since every cocycle $(\overline{\chi})^n$ in $E_1$-page can be 
lifted to a Poisson central element $\chi^n$ in $PH^0(A)$, 
the elements in $Z=\kk[\overline{\chi}]$ are all permanent 
cocycles and survive to $E_\infty$-page. As a consequence, the 
differentials $d_r^{p,q}: E_r^{p,q}\to E_r^{p+r,q-r+1}$ are all 
zero whenever $p+q=0$. So $E_1^{p,q}=E_\infty^{p,q}$ when $p+q=0$. 
By our previous claim, we know $PH^1({\rm gr}_\mathbb FA)=ZE^{ind}$. 
Similarly, every cocycle $\overline{\chi}^nE^{ind}$ can be lifted 
to a Poisson derivation $\chi^nE$ in $PH^1(A)$. Hence 
$E_1^{p,q}=E_\infty^{p,q}$ when $p+q=1$. This implies that 
$PH^1(A)=\kk[\chi]E$ and hence $A$ is $PH^1$-minimal. Finally, 
$A$ is $H$-ozone by \cite[Lemma 7.5]{TWZ}. 
\end{proof}

By combining \Cref{zzlem6.1} and \Cref{zzlem6.2}, we obtain 
the following consequence.

\begin{proposition}
\label{zzpro6.3}
If $\Omega$ is of the form in Case (1), then $A_{\Omega}$ is 
$H$-ozone.
\end{proposition}

\begin{proof}
We first show that any such $A_\Omega$ satisfies the assumptions 
in \Cref{zzlem6.2}, and then the result follows. We use 
$\Omega=z^2+x^2y^2+x^4y$ as an illustration. Note that the 
original grading on $A$ is given by $\deg(x)=1,\deg(y)=2,\deg(z)=3$. 
Now we set $\deg_{new}(x)=3,\deg_{new}(y)=2,\deg_{new}(z)=7$ and 
consider the following algebra filtration 
$\mathbb F=\{F_i\,|\,i\in \mathbb N\}$, where $F_i$ are spanned by 
all monomials $x^jy^kz^l$ satisfying $3j+2k+7l\le i$. By 
\cite[Lemma 2.9]{HTWZ1}, it is easy to check that 
$\{F_i,F_j\}\subseteq F_{i+j+2}$ for all possible 
$i,j\in \mathbb N$. So $\mathbb F$ is a $w$-filtration for the 
Poisson algebra $A_\Omega$ with $w=2$. It is clear that 
${\rm gr}_\mathbb F A= \kk[\overline{x},\overline{y},\overline{z}]$ 
with the new grading $\deg_{new}(\overline{x})=3,
\deg_{new}(\overline{y})=2,\deg_{new}(\overline{z})=7$. 
One can further verify that ${\rm gr}_\mathbb F A$ is still 
unimodular with a homogeneous potential 
$\overline{\Omega}=\overline{z}^2+\overline{x}^4\overline{y}$. By 
\Cref{zzlem6.1}, $(\gr_\mathbb F A)_{\overline{\Omega}}$ is 
$H$-ozone. Moreover, it is routine to check that 
$(\gr_\mathbb F A)_{\overline{\Omega}}$ satisfies all the 
requirements in \Cref{zzlem6.2}. So $A_\Omega$ is $H$-ozone. 
\end{proof}

\subsection{Case (2)}
\label{zzsec6.2}
\begin{proposition}
\label{zzpro6.4}
If $\Omega$ is of the form in Case (2), then $A_{\Omega}$ is 
$H$-ozone.
\end{proposition}

\begin{proof} If $\Omega$ is of the form $xyz+x^{b'}+y^{a'}$, by \Cref{zzlem5.3}, $\Omega$ is $K_1$-sealed. The assertion follows from \Cref{zzlem5.11}.
\end{proof}

\subsection{Case (3)}
\label{zzsec6.3}
In this subsection, we show that $A_\Omega$ is not $H$-ozone for each non-balanced irreducible $\Omega$. 

\begin{lemma}
\label{zzlem6.5}
Let $\Omega$ be non-balanced and irreducible. 
Then $A$ is not $H$-ozone.
\end{lemma}

\begin{proof} By a choice of $(x,y,z)$, we may assume that
$\Omega=h(x,y)=\sum_{ai+bj=n} \alpha_{ij} x^{i} y^{j}$ 
where $n=a+b+c$. Since $\Omega$ is irreducible, $\alpha_{0k}$
and $\alpha_{l0}$ are nonzero where $kb=n=al$. Since
$c=n-a-b=kb -a-b$, the assumption $\gcd(a,b,c)=1$ implies
that $\gcd(a,b)=1$. Then the equation $kb=al$ implies that
$k=ag$ and $l=bg$ where $g=\gcd(k,l)$. If $\alpha_{ij}\neq 0$
for some $(i,j)$, then $ai+bj=n=abg$. Hence $a\mid j$
and $b\mid i$. This means that $h(x,y)=f(x^b, y^a)$
for some homogeneous polynomial $f(s,t)$. Since $\Omega$ is
irreducible, so is $f(s,t)$. The only possibility is 
when $f(s,t)$ is linear, or equivalently, $h(x,y)=x^b+y^a$.
Since $ab=n=a+b+c$, $a,b\geq 2$.

To show $A$ being not $H$-ozone, it suffices to show that 
$\Omega$ does not satisfy condition (2) in \Cref{zzlem5.10}. 
Considering $\overrightarrow{f}=(0,0,1)\in A^{\oplus3}$. It is 
easy to check that 
$\overrightarrow{\nabla}\cdot \overrightarrow{f}
=\overrightarrow{f}\cdot \overrightarrow{\nabla}\Omega=0$. 
Suppose 
$\overrightarrow{f}
=\overrightarrow{\nabla}g\times \overrightarrow{\nabla}\Omega$ 
for some $g\in A$. Write 
$\overrightarrow{\nabla}g=(g_1,g_2,g_3)\in A^{\oplus 3}$. Then we must 
have $1=ag_1y^{a-1}-bg_2x^{b-1}$, which is impossible since 
$a,b\ge 2$. 
\end{proof}

\subsection{Case (4)}
\label{zzsec6.4}
Finally, we deal with the two exceptional potentials: $\Omega=z^2+y^3+x^2 y^2$ or $z^2+x^2y^2+x^{2+\frac{2b}{a}}$.

\begin{proposition}
\label{zzpro6.6}
If $\Omega$ is of one of the above forms, then $A_{\Omega}$ is H-ozone.
\end{proposition}
\begin{proof}
By \Cref{lem5.5}(2)-(3) and \Cref{lem5.6}, $\Omega$ is $K_1$-sealed. The assertion follows
from \Cref{zzlem5.11}.
\end{proof}

\begin{corollary}
\label{zzcor6.7}
Let $\Omega$ be an irreducible potential in the classification 
that is neither $x^{k}+y^{l}$ nor $x^{k}+z^{l}$ 
nor $y^{k}+z^{l}$. Then, $\Omega$ is balanced.
\end{corollary}

\begin{proof}
By \Cref{zzpro6.3}, \Cref{zzpro6.4} and \Cref{zzpro6.6}, for such an 
$\Omega$, $A_{\Omega}$ is $H$-ozone. The assertion follows
from \Cref{zzlem6.5}.
\end{proof}

\subsection{Main results on Poisson cohomology}
\label{zzsec6.5}
We compute the Hilbert series of the Poisson cohomology groups 
of $A_\Omega$ when it is $H$-ozone. Our result shows that 
their Hilbert series only depends on the grading of $A_\Omega$.

\begin{theorem}
\label{zzthm6.8}
Let $A=\kk[x,y,z]$ be a connected graded algebra such that 
$\deg(x)=a,\deg(y)=b,\deg(z)=c$ with unimodular Poisson structure given by some homogeneous polynomial $\Omega$ of 
degree $n>0$. Here $n$ may not be necessarily equal to $a+b+c$. Suppose the following statements hold:
\begin{itemize}
\item[(a)] $Z_P(A)=\kk[\Omega]$. 
\item[(b)] $A$ is $H$-ozone. 
\item[(c)] $A$ has a degree zero Poisson derivation that is not ozone. 
\item[(d)] $A$ has no non-zero Poisson derivation of degree $-n$. 
\end{itemize}
Then, the Hilbert series of the Poisson cohomology groups of $A$ is 
given by  
\begin{enumerate}
\item[(1)] 
$h_{PH^0(A)}(t)=\frac{1}{1-t^n}$.
\item[(2)] 
$h_{PH^1(A)}(t)=\frac{1}{1-t^n}$.
\item[(3)]
$h_{PH^2(A)}(t)=\frac{1}{t^{a+b+c}}\left(\frac{(1-t^{n-a})
(1-t^{n-b})(1-t^{n-c})}{(1-t^n)(1-t^a)(1-t^b)(1-t^c)}-1\right)$.
\item[(4)]
$h_{PH^3(A)}(t)=\frac{(1-t^{n-a})(1-t^{n-b})
(1-t^{n-c})}{t^{a+b+c}(1-t^n)(1-t^a)(1-t^b)(1-t^c)}$.
\end{enumerate}
\end{theorem}

\begin{proof}
Consider the cochain complex 
$(\mathfrak X^\bullet(A),d_\pi^\bullet)$ in \eqref{E1.3.1} for 
computing $PH^\bullet(A)$. 

(1) By (a), it is clear that 
$h_{PH^0(A)}(t)=h_{Z_P(A)}(t)=\frac{1}{1-t^n}$. 

(2) By (c), we have a degree zero Poisson derivation, denoted by 
$\delta$, such that $\delta(\Omega)\neq 0$. Without loss of 
generality, we can assume $\delta(\Omega)=\Omega$ as $\delta$ has
degree zero. By (d) and a similar argument of \cite[Lemma 7.7]{TWZ} 
(also see the proof of \Cref{zzlem6.2}), one can see that  
$PH^1(A)\cong Z_P(A)\delta$. So 
$h_{PH^1(A)}(t)=h_{Z_P(A)}(t)=\frac{1}{1-t^n}$.

(3) By (b) and \Cref{zzlem5.10}, $A$ is $uPH^2(A)$-vacant, namely,  
$lPH^2(A)=PH^2(A)$. Then the result follows from \Cref{zzlem5.9} 
and parts (1, 2).

(4) Notice $\ker(d_\pi^2)=M^2(A)$ since $A$ is $uPH^2(A)$-vacant.  Considering 
\eqref{E5.9.2}--\eqref{E5.9.3} and the following exact sequence
\begin{equation}
\label{E6.7.1}\tag{E6.7.1}
0\to Pd(A)[w]\to \mathfrak X^1(A)[w]\to \ker(d_\pi^2)[2w]
\to PH^2(A)[2w]\to 0,
\end{equation}
we obtain
\begin{align}
\label{E6.7.2}\tag{E6.7.2}
h_{PH^3(A)}(t)=h_{\ker(d_\pi^2)}(t)t^w
 -h_{\mathfrak X^1(A)}(t)t^{2w}+h_{\mathfrak X^0(A)}(t)t^{3w}
 +\frac{1}{t^{a+b+c}}\frac{(1-t^{w+a})(1-t^{w+b})(1-t^{w+c})}
{(1-t^a)(1-t^b)(1-t^c)}
\end{align}
via \eqref{E1.3.7}. The rest can be deduced from \eqref{E5.9.1} via another direct computation. 
\end{proof}

As a consequence, we obtain the Hilbert series of Poisson 
cohomology for any connected graded unimodular Poisson algebra 
$A_\Omega$ for any balanced irreducible $\Omega$ (but not 
necessarily having isolated singularity).

\begin{corollary}
\label{zzcor6.9}
Assume \Cref{zzhyp0.1}. If $\Omega$ is an irreducible potential in the classification that is neither $x^{k}+y^{l}$ nor $x^{k}+z^{l}$ 
nor $y^{k}+z^{l}$, then the Hilbert series of Poisson cohomology of $A$ are given as follows.
\begin{enumerate}
\item[(1)] 
$h_{PH^0(A)(t)}=\frac{1}{1-t^n}$.
\item[(2)]
$h_{PH^1(A)}(t)=\frac{1}{1-t^n}$.
\item[(3)]
$h_{PH^2(A)}(t)=\frac{1}{t^n}\left(\frac{(1-t^{a+b})(1-t^{a+c})
(1-t^{b+c})}{(1-t^n)(1-t^{a})(1-t^{b})(1-t^{c})}-1\right)$.
\item[(4)]
$h_{PH^3(A)}(t)=\frac{(1-t^{a+b})(1-t^{a+c})(1-t^{b+c})}
{t^n(1-t^n)(1-t^{a})(1-t^{b})(1-t^{c})}$.
 \end{enumerate}
\end{corollary}

\begin{proof}
It suffices to check that $A_\Omega$ satisfies all the 
requirements in \Cref{zzthm6.8} when $\Omega$ is an irreducible potential mentioned above. (a), (c) and (d) are obvious. So it remains to show all such $A_\Omega$ are $H$-ozone. If $\Omega$ has isolated singularity, $A_{\Omega}$ is $H$-ozone. As discussed before \Cref{zzlem6.1}, the rest of the irreducible $\Omega$ classified in \Cref{zzthm2.5} are divided into 
four classes. Note that \Cref{zzpro6.3} covers class (1), 
\Cref{zzpro6.4} covers class (2), \Cref{zzpro6.6} covers class (4), and \Cref{zzlem6.5} covers class (3). So our result follows immediately by letting $n=a+b+c$.
\end{proof}

The Hilbert series of the Poisson cohomology only depends on the weights of $x, y, z$ when $\Omega$ is balanced. Thus, these connected graded unimodular Poisson algebras exhibit the same homological behavior, making it impossible to distinguish irreducible potentials solely by using the Hilbert series. 
It would be interesting to see if additional structures in Poisson cohomology can distinguish different irreducible potentials. See the next question.  

\begin{question}
\label{zzque6.10}
Can we use $\kk[\Omega]$-module structures on $PH^\bullet(A)$ to distinguish between singular and smooth curves and identify types of singularity for $\Omega=0$?  
\end{question}

Finally, we generalize \cite[Theorem 0.6]{TWZ} by removing the restriction that generators are in degree one and including more equivalent conditions involving the second Poisson cohomological group. The result is stated as  \Cref{zzthm0.6}.

\begin{proof}[Proof of \Cref{zzthm0.6}]
(1) $\Rightarrow$ (2): Since $rgt(A)=0$, every degree zero
semi-Poisson derivation $\delta$ for $A$ is of the form $\alpha E$ 
for some $\alpha\in \kk$ where $E$ is the Euler derivation. 
Thus $E\wedge \delta=0$. By \eqref{E1.1.2}, $\pi_{new}=\pi$. 
So $A=A^{\delta}$. The assertion follows. In particular, $A$ is unimodular by \cite[Corollary 0.3]{TWZ}.

(2) $\Rightarrow$ (1): By \cite[Corollary 0.3]{TWZ}, there 
is a Poisson derivation $\delta$ of degree zero such that 
$A^{\delta}$ is unimodular. Since $A^{\delta}\cong A$ for 
all $\delta$, $A$ is unimodular. Suppose, to the contrary 
that $A$ is not rigid, that is $rgt(A)\neq 0$. Then, there is a Poisson derivation 
$\delta$ of degree zero not in $\kk\,E$. Thus, by \cite[Theorem 0.2]{TWZ}, the modular derivation of $A^{\delta}$ is
$${\mathbf n}=0+n\delta-\divv(\delta) E$$
which cannot be zero as $\divv(\delta)\in \Bbbk$ 
\cite[Lemma 1.2(3)]{TWZ}. Therefore, $A^{\delta}$ is not 
isomorphic to $A$, yielding a contradiction.

(1) $\Leftrightarrow$ (8): In the proof of (1)$\Leftrightarrow$(2), 
(1) implies that $A$ is unimodular with a potential $\Omega$. By \Cref{zzthm4.2}, $\Omega$ 
is irreducible. According to \Cref{zzthm2.5} and \Cref{zzcor6.7}, up to graded 
isomorphisms of $A$, the non-balanced potentials are 
$\Omega=z^2+y^3$, $z^2+x^{2+\frac{2b}{a}}$,  or $\Omega=f(x,y)$ 
in \Cref{zzlem2.2}. It is easy to check that these corresponding 
Poisson algebras have Poisson derivations of negative degree. For example, one can define a Poisson derivation $\delta$ such that $\delta(z)=\delta(y)=0, \delta(x)=1$ when $\Omega=z^2+y^3$, and $\delta(z)=\delta(x)=0, \delta(y)=1$ when $\Omega=z^2+x^{2+\frac{2b}{a}}$, and $\delta(x)=\delta(y)=0, \delta(z)=1$ when $\Omega=f(x,y)$. So we get (8). The other direction follows from the calculation 
of $PH^1(A)$ in \Cref{zzcor6.9} and \cite[Remark 5.2]{TWZ}. 

(3) $\Leftrightarrow$ (5): This follows from \cite[Proposition 7.4]{TWZ}.

(5) $\Leftrightarrow$ (6): Under the hypothesis (5), $A$ is 
$PH^1$-minimal. One implication follows by \cite[Lemma 7.5]{TWZ}
and the other direction is clear.

(6) $\Rightarrow$ (7): See the proof of \cite[Lemma 7.5]{TWZ}.

(7) $\Rightarrow$ (1): Note that $A$ is unimodular in this case. From the classification of $\Omega$ in 
\Cref{zzthm2.5}, it is easy to see that $Z_P(A)=\kk[\chi]$ with $\deg(\chi)>0$. Thus, $rgt(A)=0$ from 
\cite[Lemma 7.7(2)]{TWZ}. In particular, $\Omega$ is irreducible and $Z_{P}(A)=\Bbbk[\Omega]$. Moreover, suppose
 $\delta\in Pd(A)_{<0}$. If $\delta(\Omega)=0$, then $\delta=H_f$ 
for some $f\in A$. This is impossible since 
$\deg(f)=\deg(\delta)<0$ and  $A$ is connected graded. If 
$\delta(\Omega)\in \kk^\times$, then $\deg(\delta)=-n$. Then 
$\deg(\delta(x))=-b-c,\deg(\delta(y))=-a-c,\deg(\delta(z))=-a-b$, 
which are all negative. Hence $\delta=0$. So $Pd(A)_{<0}=0$.  

(8) $\Rightarrow$ (4,5,10,11): This assertion follows from 
\Cref{zzcor6.9}.

(4) $\Rightarrow$ (1): The assertion $rgt(A)=0$ follows from 
\cite[Remark 5.2]{TWZ}. Moreover, since $A$ is connected graded, 
each non-zero Hamiltonian derivation is of positive degree. 
As a result, we get $Pd(A)_{<0}=PH^1(A)_{<0}=0$.  

(5) $\Leftrightarrow$ (9): It follows from \eqref{E1.3.7} that 
$h_{PH^3(A)}(t)-h_{PH^2(A)}(t)
=h_{PH^0(A)}(t)-h_{PH^1(A)}(t)+t^{-n}$.
We know that $PH^0(A)=Z_P(A)$. So, the assertion follows from the 
fact that $h_{PH^1(A)}(t)=h_{PH^0(A)}(t)=h_Z(t)$ if and only 
if $h_{PH^3(A)}(t)-h_{PH^2(A)}(t)=t^{-n}$.

(10) $\Rightarrow$ (7): Note that the subspace 
$M^2(A)=\{f\pi_\Omega+\pi_g\,|\, f,g\in A\}$ of $\mathfrak X^2(A)$ 
lies in $\ker(d_\pi^2)$ [Definition-Lemma \ref{zzlem5.7}]. 
For any two locally finite graded vector spaces $M$ and $N$, we 
use the notation $h_M(t)\ge h_N(t)$ to mean $\dim M_i\ge \dim N_i$ 
for each $i\in \mathbb Z$. By \eqref{E5.9.1}--\eqref{E5.9.3} and 
\eqref{E6.7.1} with $w=0$, we get 
\begin{align*}
h_{PH^2(A)}(t)
&=h_{\ker(d_\pi^2)}(t)-h_{\mathfrak X^1(A)}(t)
+h_{\mathfrak X^0(A)}(t)+(h_{PH^1(A)}(t)-h_{PH^0(A)}(t))\\
&\ge h_{M^2(A)}(t)-h_{\mathfrak X^1(A)}(t)
+h_{\mathfrak X^0(A)}(t)+( h_{Z_P(A)E}(t)-h_{PH^0(A)}(t))\\
&\ge h_{M^2(A)}(t)-h_{\mathfrak X^1(A)}(t)+h_{\mathfrak X^0(A)}(t)\\
&=\frac{1}{t^n}\left(\frac{(1-t^{a+b})(1-t^{a+c})(1-t^{b+c})}
{(1-t^n)(1-t^{a})(1-t^{b})(1-t^{c})}-1\right)
\end{align*}
where the last equality follows from the exact sequence 
\eqref{E5.9.1} together with a direct computation. By the 
assumption, one can obtain $M^2(A)=\ker(d_\pi^2)$ and so $A$ is 
$uPH^2$-vacant. This is equivalent to $A$ being $H$-ozone by 
\Cref{zzlem5.10}. 

(11) $\Rightarrow$ (7): The argument is similar to the proof 
of  (10) $\Rightarrow$ (7) by using \eqref{E6.7.2}.

(7) $\Leftrightarrow$ (12): It follows from \Cref{zzlem5.10}.
\end{proof}
\pagebreak

\appendix

\section{Classification of Weight Polynomials $\Omega$ of degree $|x|+|y|+|z|$ in $\kk[x,y,z]$} 
\label{appendix}


The classification of $\Omega$ when $\deg(x)=\deg(y)=\deg(z)=1$
is well-known; see \cite{BM, DH, DML, KM, LX} and the first table below.

\begin{table}[H]
\caption{\label{tab:111} (a, b, c)=(1, 1, 1).}
\begin{center}
\begin{tabular}{|cc|cc|}
 \hline
 {\rm Irreducible} $\Omega$ &&& {\rm Reducible} $\Omega$\\ [0.2em] 
 \hline
$x^3+y^2z$~~\raisebox{.5pt}{\textcircled{\raisebox{-.5pt} {w}}} , 
$x^3+x^2z+y^2z$~~\raisebox{.5pt}{\textcircled{\raisebox{-.5pt} {q}}}
&&& $ x^3, x^2y, xyz$ \\[0.2em]
$x^3+y^3+z^3+3\lambda xyz$ $(\lambda^3 \ne -1 ,\raisebox{.5pt}{\textcircled{\raisebox{-.5pt} {i}}})$ &&& $xy(x+y), xyz+x^3, xy^2+x^2z$ \\ [0.2em] 
 \hline
\end{tabular}
\end{center}
\end{table}

\begin{table}[H]
\caption{\label{tab:112} (a, b, c)=(1, 1, 2).}
\begin{center}
\begin{tabular}{|cc|cc|}
 \hline
 {\rm Irreducible} $\Omega$ &&& {\rm Reducible} $\Omega$\\ [0.2em] 
 \hline
$z^2+x^3y$~~\raisebox{.5pt}{\textcircled{\raisebox{-.5pt} {w}}} , 
$z^2+ x^2y^2+ x^3y$~~\raisebox{.5pt}{\textcircled{\raisebox{-.5pt} {q}}}
&&& $ x^4, x^3y, x^2y^2, x^2y^2+x^3y$ \\[0.2em]
$z^2+xy^3+\lambda x^2y^2+x^3y$ &&& $z^2+x^2y^2, z^2, z^2+x^4$, \\ [0.2em] 
$(\lambda \ne \pm 2,\raisebox{.5pt}{\textcircled{\raisebox{-.5pt} {i}}})$,\,
$(\lambda=\pm 2, \raisebox{.5pt}{\textcircled{\raisebox{-.5pt} {q}}} )$
&&& $xy^3+\lambda x^2y^2+x^3y$ \\[0.2em]
 
$x^2z+y^4$~~\raisebox{.5pt}{\textcircled{\raisebox{-.5pt} {w}}} , 
$x^2z+xy^3+y^4$ ~~\raisebox{.5pt}{\textcircled{\raisebox{-.5pt} {w}}} 
&&& $x^2z, x^2z+xy^3$\\[0.2em]
 
 $xyz+x^4+y^4$~~\raisebox{.5pt}{\textcircled{\raisebox{-.5pt} {q}}} 
&&&$xyz, xyz+x^4$\\ [0.2em]
 \hline
\end{tabular}
\end{center}
\end{table}

\begin{table}[H]
\caption{\label{tab:a=b<c} $(a=b<c)\neq (1, 1, 2)$}
    \begin{center}
\begin{tabular}{|c|c|c|}
 \hline
 $k\in \mathbb N$ & {\rm Irreducible} $\Omega$ & {\rm Reducible} $\Omega$ \\ [0.4em] 
 \hline
$c\ne ka$ &  & $xyz, x^2z$\\ [0.4em] 
\hline
&  $x^2z+xy^{k+1}+y^{k+2}$ ~~\raisebox{.5pt}{\textcircled{\raisebox{-.5pt} {w}}}  
& $x^2z, x^2z+xy^{k+1}$  \\[0.4em]
$c=ka$ &$x^2z+y^{k+2}$ ~~\raisebox{.5pt}{\textcircled{\raisebox{-.5pt} {w}}} 
& $xyz, xyz+x^{k+2}$\\[0.4em]
$k\ne 2$ &$xyz+x^{k+2}+y^{k+2}$ ~~\raisebox{.5pt}{\textcircled{\raisebox{-.5pt} {q}}}  
& $h(x,y)$ {\rm of degree} $(k+2)a$\\[0.4em]
 \hline 
\end{tabular}
\end{center}
\end{table}

\begin{table}[H]
\caption{\label{tab:a<b=c} $(a<b=c)$}
    \begin{center}
\begin{tabular}{|c|}
 \hline
 {\rm Reducible} $\Omega$ \\ [0.4em] 
 \hline
 $xyz, \,xy^2,\,xz^2+x^{1+\frac{2b}{a}},\, xyz+x^{1+\frac{2b}{a}}$,\, 
$x^{1+\frac{b}{a}}y+xz^2$,\, $x^{1+\frac{2b}{a}}, x^{1+\frac{b}{a}}y$ \\[0.4em]
 \hline 
\end{tabular}
\end{center}
\end{table}

\begin{table}[H]
\caption{\label{tab:123} (a, b, c)=(1, 2, 3).}
\begin{center}
\begin{tabular}{|cc|cc|}
 \hline
 {\rm Irreducible }$\Omega$ &&& {\rm Reducible} $\Omega$ \\ [0.4em] 
 \hline
$z^2+y^3$ ~~\raisebox{.7pt}{\textcircled{\raisebox{-.01pt} {w}}} 
&&&   $z^2$, $y^3, x^6$, $y^3+x^2y^2$\\ [0.4em] 
$z^2+y^3+x^2y^2$ ~~\raisebox{.5pt}{\textcircled{\raisebox{-.5pt} {q}}}  
&&& $x^2y^2$, $x^2y^2+x^4y$ \\[0.4em]
$z^2+y^3+\lambda x^2y^2+x^4y$ &&&  $y^3+\lambda x^2y^2+x^4y$\\[0.4em]
$(\lambda \ne \pm 2,\raisebox{.5pt}{\textcircled{\raisebox{-.5pt} {i}}})$, 
$(\lambda=\pm 2, \raisebox{.5pt}{\textcircled{\raisebox{-.5pt} {q}}})$ 
&&&$x^4y$ \\[0.4em]
$z^2+x^4y$~~\raisebox{.5pt}{\textcircled{\raisebox{-.5pt} {w}}}   
&&&  $z^2+x^6$\\[0.4em]
$z^2+x^2y^2+x^4y$ ~~\raisebox{.5pt}{\textcircled{\raisebox{-.5pt} {q}}} 
 &&& $z^2+x^2y^2$\\[0.4em]
 
$x^3z+y^3+ x^2y^2$  ~~\raisebox{.5pt}{\textcircled{\raisebox{-.5pt} {w}}} , 
$x^3z+y^3$~~\raisebox{.5pt}{\textcircled{\raisebox{-.5pt} {w}}}  
&&& $x^3z$, $x^3z+x^2y^2$ \\[0.4em]
 
$xyz+x^6+y^3$~~\raisebox{.5pt}{\textcircled{\raisebox{-.5pt} {q}}}  
&&& $xyz, xyz+x^6$, $xyz+y^3$\\[0.4em] 
\hline
\end{tabular}
\end{center}
\end{table}

\begin{table}[H]
\caption{\label{tab:abc2} $(a<b<c)\neq (1, 2, 3)$}
\begin{tabular}{|c|c|c|}
 \hline
 $m,n\in \mathbb N\cup \{-1\}$ & {\rm Irreducible }$\Omega$ & {\rm Reducible} $\Omega$  \\ [0.4em] 
 \hline
$c\ne ma+nb$ &  & $xyz$, \, $ x^{1+\frac{b}{a}}z$ \\ [0.4em] 

\hline

$c=a+b$ & $z^2+x^{2+\frac{2b}{a}}$~~\raisebox{.5pt}{\textcircled{\raisebox{-.5pt} {w}}}   
& $z^2, z^2+x^2y^2$  \\[0.4em]
$a\nmid b$& $z^2+x^2y^2+x^{2+\frac{2b}{a}}$~~\raisebox{.5pt}{\textcircled{\raisebox{-.5pt} {q}}}  
& $xyz, xyz+x^{2+\frac{2b}{a}}$ \\[0.4em]
$\quad\quad\quad$ & $\quad\quad\quad$ & $h(x,y)$ of degree $(2a+2b)$\\[0.4em]
  
\hline

& &$ z^2, z^2+x^{2+\frac{2b}{a}}$  \\[0.4em]
$c=a+b$ & $z^2+x^2y^2+x^{2+\frac{b}{a}}y$~~\raisebox{.5pt}{\textcircled{\raisebox{-.5pt} {q}}} 
& $z^2+x^2y^2, x^{1+\frac{b}{a}}z$  \\[0.4em]
$b\ne 2a,\, a\mid b$ & $ z^2+x^{2+\frac{b}{a}}y$~~\raisebox{.5pt}{\textcircled{\raisebox{-.5pt} {w}}} 
& $ x^{1+\frac{b}{a}}z+x^2y^2, xyz $ \\ [0.4em]
&  & $xyz+x^{2+\frac{2b}{a}}$ \\[0.4em]
$\quad\quad\quad$ & $\quad\quad\quad$ & $h(x,y)$ of degree $(2a+2b)$\\[0.4em]
\hline

$c=ma+nb$  &   & $ xyz+x^{m+1+k}$ \\[0.4em]
$c\ne a+b$,\,  $a\nmid b$ 
& $ xyz+x^{m+1+k}+y^{n+1+l}$~~\raisebox{.5pt}{\textcircled{\raisebox{-.5pt} {q}}}  
& $xyz+y^{n+1+l}$  \\[0.4em]
$\frac{m+1}{k}=\frac{l}{n+1}$& $x^{\frac{b}{\gcd(a,b)}}+y^{\frac{a}{\gcd(a,b)}}$~~\raisebox{.5pt}{\textcircled{\raisebox{-.5pt} {w}}} & $xyz$ 
\\[0.4em]
$\quad\quad\quad$ & $\quad\quad\quad$ & $h(x,y)$ of degree $(m+1)a+(n+1)b$\\[0.4em]
\hline
$c=ma+nb$ 
& $ xyz+x^{m+1+k}+y^{n+1+l}$ ~~\raisebox{.5pt}{\textcircled{\raisebox{-.5pt} {q}}}  
& $ x^{1+\frac{b}{a}}z, xyz+x^{m+1+k}$ \\[0.4em]
$c\ne a+b$,\, $a\mid b$  &$x^{1+\frac{b}{a}}z+y^{1+n+l}$ ~~\raisebox{.5pt}{\textcircled{\raisebox{-.5pt} {w}}} 
&  $xyz, xyz+y^{n+1+l}$ \\[0.4em]
$\frac{m+1}{k}=\frac{l}{n+1}$ & $x^{1+\frac{b}{a}}z+ x^{\frac{b}{a}}y^{n+l}+y^{1+n+l}$  ~~\raisebox{.5pt}{\textcircled{\raisebox{-.5pt} {w}}} & $x^{1+\frac{b}{a}}z+x^{\frac{b}{a}}y^{n+l}$ \\[0.4em]
$\quad\quad\quad$ & $\quad\quad\quad$ & $h(x,y)$ of degree $(m+1)a+(n+1)b$\\[0.4em]
\hline
\end{tabular}
\end{table}

\begin{table}[H]
\caption{\label{tab:rgtanddim} $rgt:=rgt(A)$ and $GK:=GKdim(A_{ sing})$ for $\Omega$ in Tables $3$-$7$} 
\begin{center}
\begin{tabular}{|c|c|c|c|c|c|}
 \hline
$\Omega$ & $rgt$ & $GK$ & $\Omega$& $rgt$ & $GK$\\
 \hline
$z^2+x^3y$ & $0$ & $1$ & $xz^2+x^{1+\frac{2b}{a}}$& $\substack{-2(a\mid b)\\ -1 (a\nmid b)}$ & $2$\\
\hline
$z^2+ x^2y^2+ x^3y$ & $0$ & $1$ &$xyz+x^{1+\frac{2b}{a}}$ & $-1$ & $1$\\
\hline 
$z^2+xy^3+\lambda x^2y^2+x^3y$
&$0$ & $\substack{0\, (\lambda \ne \pm 2) \\ 1\, (\lambda=\pm 2)}$ &$x^{1+\frac{b}{a}}y+xz^2$ & $-1$ & $1$\\
\hline
$x^2z+y^4$  &$0$&  $1$  & $x^{1+\frac{2b}{a}}$& $\substack{-5(a\mid b)\\ -3(a\nmid b)}$& $2$\\
\hline
$x^2z+xy^3+y^4$ 
 &$0$ &  $1$ &$x^{1+\frac{b}{a}}y\,\, (\ast)$ & $-3$ & $2$\\
\hline
$xyz+x^4+y^4$  & $0$& $1$  &$x^2z+xy^{k+1}+y^{k+2}$ & $0$& $1$\\
\hline 
$ x^4$ & $-5$ & $2$  &$x^2z+y^{k+2}$  &$0$ & $1$\\
\hline
 $x^3y$& $-4$ &  $2$ & $xyz+x^{k+2}+y^{k+2}$  & $0$  & $1$\\
\hline
$x^2y^2$& $-4$& $2$ & $xyz$& $-2$ & $1$\\
\hline
$x^2y^2+x^3y$& $-3$ & $2$ &$x^2z$ & $-2$ & $2$ \\
\hline
 $z^2+x^2y^2$& $-1$ & $1$ & $x^2z$& $-2$ & $2$\\ 
 
 \hline
 $z^2$ &$-3$ & $2$ &$x^2z+xy^{k+1}$ & $-1$ & $1$\\

\hline
$z^2+x^4$ & $-1$ & $1$ &$xyz$ & $-2$ & $1$\\
\hline
 $xy^3+\lambda x^2y^2+x^3y$ & $-3$ & $1$ &$xyz+x^{k+2}\,\, (\ast)$ &$-1$ & $1$\\
\hline
 $x^2z$ & $-2$ & $2$ & $z^2+x^{2+\frac{2b}{a}} $ & $0 (a\nmid b)$& $1$\\
\hline
$x^2z+xy^3$ & $-1$ &  $1$ & $z^2+x^2y^2+x^{2+\frac{2b}{a}}$ & $0$& $1$\\
\hline 

$xyz$ & $-2$ & $1$  &$z^2+x^2y^2+x^{2+\frac{b}{a}}y$  & $0$ & $1$\\
\hline
$xyz+x^4\,\, (\ast)$ & $-1$ & $1$ & $z^2+x^{2+\frac{b}{a}}y$ & $0$& $1$\\
 \hline
 $z^2+y^3$ &$0$ & $1$  & $xyz+x^{m+1+k}+y^{n+1+l}\, (a\nmid b)$ & $0$& $1$\\
 \hline
 $z^2+y^3+x^2y^2$  &$0$ & $1$  &$xyz+x^{m+1+k}+y^{n+1+l}\, (a\mid b)$ & $0$ & $1$\\
 \hline
 $z^2+y^3+\lambda x^2y^2+x^4y$ &$0$ &  $\substack{0 (\lambda \ne \pm 2)\\ 1 (\lambda=\pm 2)}$ & $x^{1+\frac{b}{a}}z+y^{1+n+l}$ & $0$ & $1$\\
 \hline
 $z^2+x^4y$ &$0$ & $1$ & $x^{1+\frac{b}{a}}z+x^{\frac{b}{a}}y^{n+l}+y^{1+n+l}$ & $0$ & $1$\\
 \hline
 $z^2+x^2y^2+x^4y$ &  $0$ & $1$ &$xyz$ &$-2$ & $1$\\
 \hline
 $x^3z+y^3+x^2y^2$ &$0$ & $1$  &$x^{1+\frac{b}{a}}z$ & $-2$& $2$\\
 \hline

 $x^3z+y^3$  &$0$ & $1$   & $z^2$& $-1(a\nmid b)$& $2$\\
 \hline 
 $xyz+x^6+y^3$ &$0$ & $1$  &$z^2+x^2y^2$ &$-1$ & $1$\\
 \hline
 $z^2$ &$-2$ & $2$ &$xyz$ &$-2$ & $1$\\
 \hline
 $y^3$ & $-3$ & $2$ &$xyz+x^{2+\frac{2b}{a}}$ & $-1$ & $1$ \\
 \hline 
 $x^6$ & $-4$ & $2$  &$z^2$ & $-2 (a\mid b)$ & $2$\\
 \hline 
 $y^3+x^2y^2$ & $-2$ & $2$ &$z^2+x^{2+\frac{2b}{a}}$ & $-1(a\mid b)$ & $1$\\
 \hline
 $x^2y^2$ & $-3$ & $2$ & $z^2+x^2y^2$ & $-1$ & $1$\\
 \hline 
 $x^2y^2+x^4y$ & $-2$ & $2$  & $x^{1+\frac{b}{a}}z$ &  $-2$& $2$\\
 \hline 
 $y^3+\lambda x^2y^2+x^4y$ & $-2$ & $\substack{1\,(\lambda\ne \pm 2)\\  2\,(\lambda =\pm 2)}$& $x^{1+\frac{b}{a}}z+x^2y^2$ &$-1$ & $2$\\
 \hline
 $x^4y$ & $-3$ & $2$ & $xyz$ & $-2$& $1$\\
 \hline 
 $z^2+x^6$ & $-1$ & $1$ &$xyz+x^{2+\frac{2b}{a}}$ &$-1$ & $1$\\
 \hline 
 $z^2+x^2y^2$ & $-1$ & $1$  & $xyz+x^{m+1+k}$& $-1$& $1$\\
 \hline 
 $x^3z$ & $-2$ &  $2$  &$xyz+y^{n+1+l}$ & $-1$& $1$\\
 \hline 
 $x^3z+x^2y^2$ &$-1$ & $2$ &$xyz$ & $-2$ & $1$\\
 \hline 
 $xyz$ & $-2$ & $1$ &$x^{1+\frac{b}{a}}z$ & $-2$& $2$\\
 \hline
 $xyz+x^6$ & $-1$ & $1$ &$xyz+x^{m+1+k}$ & $-1$& $1$\\
 \hline
 $xyz+y^3\,\, (\ast)$ & $-1$ & $1$ &$xyz$ & $-2$& $1$\\
 \hline
 $xyz$ & $-2$   & $1$ & $xyz+y^{n+1+l}$ &$-1$ & $1$\\
 \hline 
 $xy^2$ & $-2$ & $2$ &$x^{1+\frac{b}{a}}z+x^{\frac{b}{a}}y^{n+l}$ & $-1$& $2$\\
 \hline
 $x^{b/\gcd(a,b)}+y^{a/\gcd(a,b)}$ & $0$   & $1$ & $\text{Reducible}\, h(x,y)$ &$\leq-1$ & $\{1, 2\}$\\
 \hline
\end{tabular}
\end{center}
A `$\ast$' indicates the end of a table.
\end{table}

\newpage
\subsection*{Acknowledgments} 
The authors would like to thank the referee for carefully reading the paper and providing many useful suggestions. Wang was partially supported by the Simons Collaboration Grant 
\#688403 and Air Force Office of Scientific Research grant 
FA9550-22-1-0272. Zhang was partially supported by the US 
National Science Foundation (No. DMS-2001015 and DMS-2302087). Part of this research work was done during 
the third author's visit to the Department of Mathematics at 
Rice University in November 2022, and the first and third 
authors' visit to the Department of Mathematics at the University of Washington in January 2023. They wish to thank Rice 
University and University of Washington for their hospitality.

\providecommand{\bysame}{\leavevmode\hbox to3em{\hrulefill}\thinspace}
\providecommand{\MR}{\relax\ifhmode\unskip\space\fi MR }
\providecommand{\MRhref}[2]{%

\href{http://www.ams.org/mathscinet-getitem?mr=#1}{#2} }
\providecommand{\href}[2]{#2}

\end{document}